\documentclass[twocolumn,10pt]{IEEEtran}
\IEEEoverridecommandlockouts

\usepackage[ruled, vlined, algo2e]{algorithm2e}

\usepackage{algorithm}% http://ctan.org/pkg/algorithm

\usepackage{graphicx}
\usepackage{epstopdf}
\graphicspath{{../}}
\usepackage{multirow}
\usepackage{float}
\usepackage{amsmath,amsthm}
\usepackage{amssymb}% http://ctan.org/pkg/amssymb
\usepackage{pifont}% http://ctan.org/pkg/pifont
\newcommand{\cmark}{\ding{51}}%
\usepackage{array}
\usepackage{mdwmath}
\usepackage{mdwtab}
\usepackage{eqparbox}
\usepackage{url}
\usepackage{nomencl}
\usepackage{cite}
\usepackage{makeidx}
\usepackage{ifthen}
\usepackage{subfigure}
\usepackage{caption}
\usepackage{pdflscape}
\usepackage{hyperref}
\usepackage{subcaption}
\usepackage{xcolor}
\usepackage[scr=boondoxo]{mathalpha}

\def\bz{\mathbf{z}}
\def\mbi{\mathbb{I}}

\def\Nikin{\mathcal{N}_{ik}^{\mathrm{in}}}
\def\Nikout{\mathcal{N}_{ik}^{\mathrm{out}}}

\def\a{{\alpha}}

\def\1{\mathbf{1}}
\def\0{\mathbf{0}}
\def\la{\langle}
\def\ra{\rangle}

\newcommand{\SAB}{$\mathcal{S\!\!-\!\!AB}$}
\newcommand{\SABTV}{$\mathcal{S\!AB}-$\textnormal{TV}}
\newcommand{\DSGTmTV}{$\textnormal{DSGT}m-$\textnormal{TV}}
\newcommand{\DGTmTV}{$\textnormal{DGT}m-$\textnormal{TV}}

\usepackage{amsmath,mathtools, amssymb, bbm, xspace}
\usepackage{booktabs,ragged2e}
\usepackage{multirow}
\usepackage[flushleft]{threeparttable}
\usepackage{graphicx}

\def\diag{\mathop{\hbox{\rm diag}}}
  
\def\exp{\mathop{\hbox{\rm exp}}}

\def\spose#1{\hbox to 0pt{#1\hss}}

\def\text #1{\hbox{\quad#1\quad}}

\def\bbG{\mathbb{G}}

\def\nthinsp{\mskip -2   mu}

\def\F{\mathcal{F}}
\def\G{_{\scriptscriptstyle G}}

\def\superstar{^{\raise 0.5pt\hbox{$\nthinsp *$}}}
\def\SUPERSTAR{^{\raise 0.5pt\hbox{$*$}}}
\def\lamstarT {\lambda^{\raise 0.5pt\hbox{$\nthinsp *$}T}}

\def\hbar{\skew{4.2}\bar h}

		\def\bkE{{\rm I\kern-.17em E}}
		\def\bk1{{\rm 1\kern-.17em l}}
		\def\bkD{{\rm I\kern-.17em D}}
		\def\bkR{{\rm I\kern-.17em R}}
		\def\bkP{{\rm I\kern-.17em P}}
		\def\bkY{{\bf \kern-.17em Y}}
		\def\bkZ{{\bf \kern-.17em Z}}

		\def\bc{\begin{center}}
		\def\be{\begin{enumerate}}
		\def\bi{\begin{itemize}}
		\def\ec{\end{center}}
		\def\ee{\end{enumerate}}
		\def\ei{\end{itemize}}
		\def\es{\end{small}}
		\def\eS{\end{slide}}

	\def\cp2problem#1#2#3#4{\fbox
		 {\begin{tabular*}{0.9\textwidth}
			{@{}l@{\extracolsep{\fill}}l@{\extracolsep{6pt}}l@{\extracolsep{\fill}}c@{}}
				#1 & & $#4 $ 
			\end{tabular*}}}

		\def\bkE{{\rm I\kern-.17em E}}
		\def\bk1{{\rm 1\kern-.17em l}}
		\def\bkD{{\rm I\kern-.17em D}}
		\def\bkR{{\rm I\kern-.17em R}}
		\def\bkP{{\rm I\kern-.17em P}}
		
		\def\bkZ{{\bf{Z}}}

\newcommand {\beeq}[1]{\begin{equation}\label{#1}}
\newcommand {\eeeq}{\end{equation}}
\newcommand {\bea}{\begin{eqnarray}}
\newcommand {\eea}{\end{eqnarray}}

\def\texitem#1{\par\smallskip\noindent\hangindent 25pt
               \hbox to 25pt {\hss #1 ~}\ignorespaces}

\newcommand{\argmin}{\operatornamewithlimits{argmin}}

\newcommand{\beq}{\begin{equation}}
\newcommand{\eeq}{\end{equation}}
\newcommand{\beqn}{\begin{eqnarray}}
\newcommand{\eeqn}{\end{eqnarray}}
\newcommand{\beqno}{\begin{eqnarray*}}
\newcommand{\eeqno}{\end{eqnarray*}}
\newcommand{\bma}{\begin{displaymath}}
\newcommand{\ema}{\end{displaymath}}
\newcommand{\bnu}{\begin{enumerate}}
\newcommand{\enu}{\end{enumerate}}
\newcommand{\bce}{\begin{center}}
\newcommand{\ece}{\end{center}}
\newcommand{\btb}{\begin{tabular}}
\newcommand{\etb}{\end{tabular}}

\def\G{{\mathbb{G}}}

\def\ba{{\mathbf{a}}}
\def\bb{{\mathbf{b}}}

\def\bx{{\mathbf{x}}}
\def\by{{\mathbf{y}}}
\def\bz{{\mathbf{z}}}
\def\bu{{\mathbf{u}}}
\def\bv{{\mathbf{v}}}

\def\b1{{\mathbf{1}}}
\def\cNinik{{\mathcal{N}^{\rm in}_{ik}}}
\def\cNoutik{{\mathcal{N}^{\rm out}_{ik}}}
\def\cNini{{\mathcal{N}^{\rm in}_{i}}}
\def\cNouti{{\mathcal{N}^{\rm out}_{i}}}
\def\diag{{\rm diag}}

\usepackage[normalem]{ulem} % to use \sout
\newtheorem{theorem}{Theorem}

\newtheorem{proposition}{Proposition}
\newtheorem{lemma}{Lemma}
\newtheorem{remark}{Remark}

\newtheorem{assumption}{Assumption}

\newcommand{\bg}{{\mathbf{g}}}

\newcommand{\g}{{\mathbf{g}}}

\newcommand{\Diag}{{\mathrm{Diag}}}

\newcommand{\bs}{{\mathbf{s}}}

\newcommand{\E}{{\mathbb{E}}}

\newcommand{\one}{{\mathbf{1}}}
\newcommand{\zero}{{\mathbf{0}}}

\newcommand{\T}{{\mathsf{T}}}

\definecolor{myBlue}{rgb}{0.80,0.85,1.00}
\definecolor{myYellow}{rgb}{0.951,1.000,0.547}

\def\la{{\langle}}
\def\ra{{\rangle}}

\def\bit{\begin{itemize}}
\def\eit{\end{itemize}}
\def\BEAS{\begin{eqnarray*}}
\def\EEAS{\end{eqnarray*}}

\def\re{{\mathbb R}}

\def\a{\alpha}
\def\b{\beta}

\def\g{\gamma}

\def\argmin{\mathop{\rm argmin}}

\newcommand{\balpha}{{\boldsymbol{\alpha}}}
\newcommand{\bbeta}{{\boldsymbol{\beta}}}
\newcommand{\Da}{{\boldsymbol{D_\alpha}}}
\newcommand{\Db}{{\boldsymbol{D_\beta}}}
\newcommand{\bxi}{{\boldsymbol{\xi}}}

\newcommand{\bdelta}{{\boldsymbol{\delta}}}
\newcommand{\bby}{{\bar y}}

\makeatletter
\def\munderbar#1{\underline{\sbox\tw@{$#1$}\dp\tw@\z@\box\tw@}}

\makeatother

\begin{document}

\title{Decentralized Federated Learning with Gradient Tracking over Time-Varying Directed Networks
% A Fast Decentralized Stochastic Optimization Algorithm over Time-Varying Directed Networks 
% Distributed Stochastic Optimization with Gradient Tracking, Acceleration and Uncoordinated Step-Sizes over Time-Varying Directed Networks
% Decentralized Federated Learning in Dynamic Networks: Accelerated Gradient Tracking with Uncoordinated Step-Sizes
% Accelerated Gradient Tracking for Decentralized Federated Learning in Time-Varying Directed Networks
}

\author{
Duong Thuy Anh Nguyen$^{*}$, \IEEEmembership{ Graduate Student Member, IEEE}, Su Wang,
%\IEEEmembership{, IEEE}, 
Duong Tung Nguyen, \IEEEmembership{Member, IEEE}, Angelia Nedi\'c, \IEEEmembership{Member, IEEE}, and H. Vincent Poor, \IEEEmembership{Life Fellow, IEEE}
\vspace{-0.1cm}
\thanks{Duong Thuy Anh Nguyen, Duong Tung Nguyen, and Angelia Nedi\'c are with the School of Electrical, Computer and Energy Engineering, Arizona State University, Tempe, AZ, USA. 
Email: \{dtnguy52,~duongnt,~Angelia.Nedich\}@asu.edu.
Su Wang and H. Vincent Poor are with the Department of Electrical and Computer Engineering, Princeton University, Princeton, NJ, USA. Email: \{hw5731,~poor\}@princeton.edu.
This material is based in part upon work supported by the NSF award CCF-2106336, the ONR award N00014-21-1-2242, the ARPA-H award SP4701-23-C-0074, and a grant from Princeton Language and Intelligence.\\
\indent \textit{Corresponding Author}: Duong Thuy Anh Nguyen.\\
\indent A preliminary version of this article was presented at the 2023 Asilomar Conference on Signals, Systems, and Computers (ACSSC) \cite{Nguyen2023SABTV}.\\
\indent $^{*}$Part of this work was done while at Princeton University.
}
}

% \author{\IEEEauthorblockN{Duong Thuy Anh Nguyen, Duong Tung Nguyen, and Angelia Nedi\'c}
% \IEEEauthorblockA{School of Electrical, Computer and Energy Engineering, Arizona State University, AZ, USA \\
% Email: dtnguy52@asu.edu, duongnt@asu.edu, Angelia.Nedich@asu.edu}
% }

\maketitle

\IEEEpubid{\begin{minipage}{\textwidth}\ \\[9pt]
\vspace{-1cm}
\!\!\!\!\!\!\!\!\!\!\!\!\!\! This work has been submitted to the IEEE for possible publication. Copyright may be transferred without notice, after which this version may no longer be accessible. \!\!\!\!\!\!\!\!\!\!\!\!\!\!
\end{minipage}} 

\begin{abstract}
We investigate the problem of agent-to-agent interaction in decentralized (federated) learning over time-varying directed graphs, and, in doing so, propose a consensus-based algorithm called \DSGTmTV. 
The proposed algorithm incorporates gradient tracking and heavy-ball momentum to distributively optimize a global objective function, while preserving local data privacy. Under \DSGTmTV, agents will update local model parameters and gradient estimates using information exchange with neighboring agents enabled through row- and column-stochastic mixing matrices, which we show guarantee both consensus and optimality. Our analysis establishes that \DSGTmTV~exhibits linear convergence to the exact global optimum when exact gradient information is available, and converges in expectation to a neighborhood of the global optimum when employing stochastic gradients.
Moreover, in contrast to existing methods, \DSGTmTV~preserves convergence for networks with uncoordinated stepsizes and momentum parameters, for which we provide explicit bounds. 
These results enable agents to operate in a fully decentralized manner, independently optimizing their local hyper-parameters.
% Finally, w
We demonstrate the efficacy of our approach via comparisons with state-of-the-art baselines on real-world image classification and natural language processing tasks.
\end{abstract}


\allowdisplaybreaks
\section{Introduction}
In an increasingly interconnected world, there is a growing demand for advanced optimization, learning, and inference models that efficiently leverage data from diverse and distributed sources~\cite{jordan2015machine,wang2023scientific}. 
Traditional machine learning (ML) approaches rely on centralizing data at a powerful server for ML model training~\cite{bkassiny2012survey}. 
However, with the growth at the edge of networks~\cite{jordan2015machine}, centralized ML has raised several concerns, such as device-to-server transmission bottlenecks, data privacy leakage, and, underutilized computing resources on edge devices, which have grown in local processing power. %resource idleness as edge devices have grown in local compute power.
% Traditional machine learning (ML) approaches rely on centralizing data at some powerful server for ML model training, and typically centralize data for training, which raises privacy concerns and introduces inefficiencies due to data movement and processing bottlenecks. 
% The exponential increase in data volume has further compounded the challenge of processing and storing data on a single device. 
On the other hand, federated learning (FL)~\cite{pmlr-v54-mcmahan17a} enables ML model training directly on the agents (e.g., devices, clients or learners), keeping data local and thereby addressing these challenges.
Specifically, FL operates in a three-step iterative process: (i) agents locally train their ML models, (ii) agents transmit their local ML models to a server, and (iii) the server aggregates and synchronizes them % agents' ML models
into a global ML model.
% Specifically, FL operates by aggregating local updates from multiple agents to build a global model without sharing raw data, thereby facilitating the collective harnessing of knowledge from distributed sources without compromising data security or privacy.}

\begin{figure}[t!]
\centering
\subfigure[Centralized FL]{
	\raisebox{0.18\height}{\includegraphics[width=0.225\textwidth]{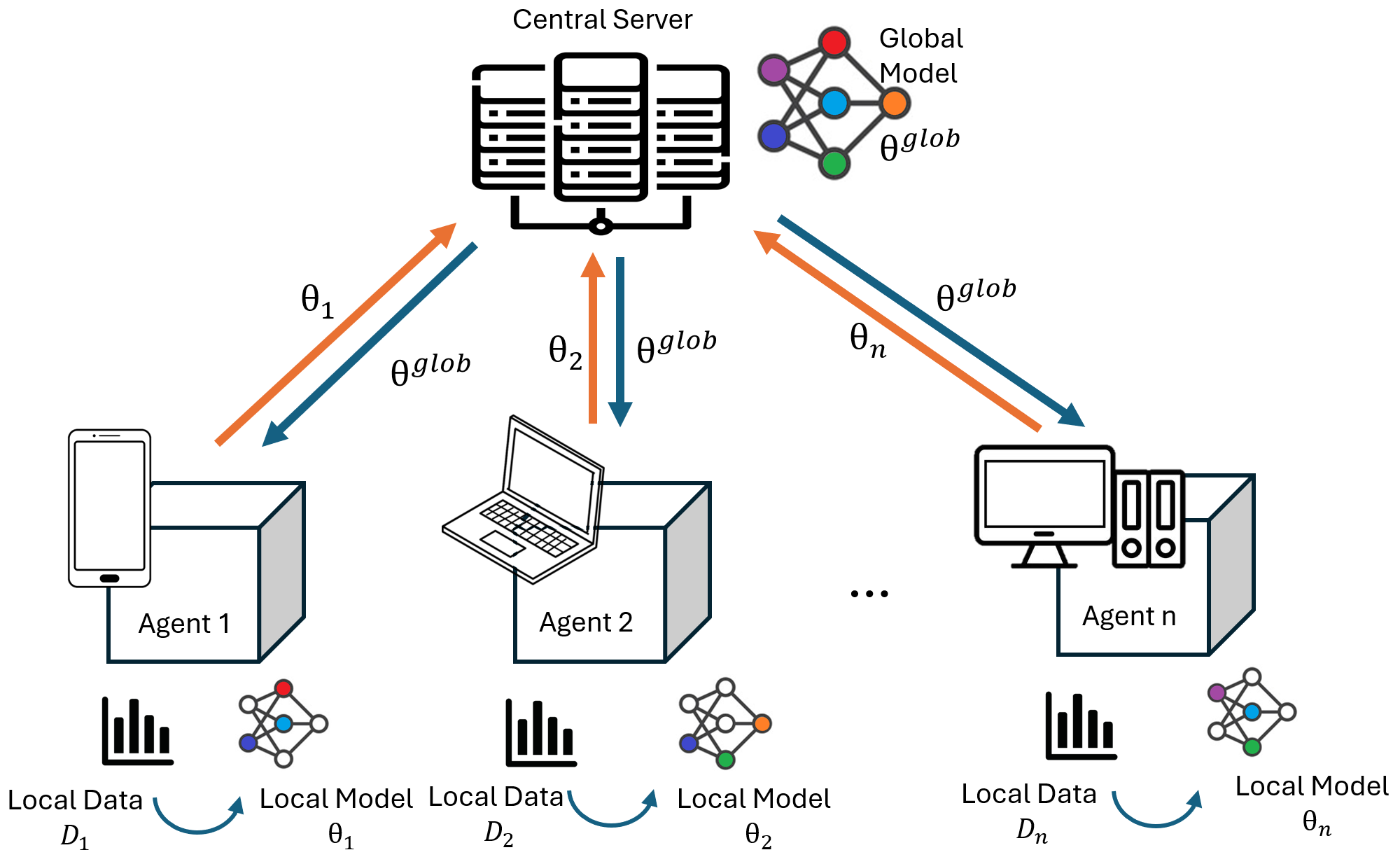}}\label{fig:FL_centralized}
	}
    \hspace*{-.7em}
\subfigure[Decentralized FL]{
	\includegraphics[width=0.23\textwidth]{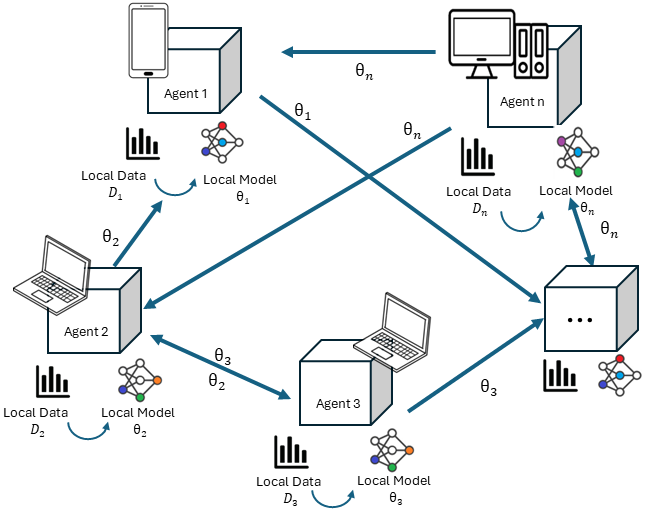}\label{fig:FL_decentralized}  }
 \vspace{-0.2cm}
 \caption{Traditional FL has a centralized server as in Fig.~\ref{fig:FL_centralized}. However, such a server is absent in decentralized FL, see Fig.~\ref{fig:FL_decentralized}. We focus on the more challenging case of decentralized FL, with directed and unbalanced links among agents.} \label{fig:FL}
\vspace{-0.4cm}
\end{figure}

Existing works in FL~\cite{ghosh2022efficient,avdiukhin2021federated,pmlr-v54-mcmahan17a,wang2019adaptive,wang2024multi} primarily focus on the centralized case, and, thus, rely on a central server.
However, this centralized architecture introduce several limitations: (i) it becomes a single point of failure, which can make the network vulnerable to attacks or server outages~\cite{lyu2022privacy}, (ii) it suffers from communication bottlenecks, exacerbated by agent failures and straggler effects~\cite{wang2021device}, and (iii) it leads to privacy risks as the server handles all ML models, which may be exploited to infer information about agents' local data~\cite{wang2019beyond}.
As they avoid the above problems and offer benefits to scalability, decentralized (federated) learning and optimization methodologies have thus gained substantial traction in many applications, e.g., big-data analytics~\cite{Daneshmand2015}, autonomous vehicles~\cite{NguyenAutonomousDriving2022}, crowdsourcing~\cite{wiopt23}, and UAV management~\cite{XiaoUAVFL2021}.
% satellite communications~\cite{YanSatellites2024}, and healthcare~\cite{XuDFL2021}. 
Fig.~\ref{fig:FL} depicts diagrams illustrating the two variants of FL. 
Decentralized FL eliminates the dependency on a central server by distributing the coordination and aggregation of model updates among participating agents, who communicate and collaborate directly in a peer-to-peer (P2P) manner.
However, fully decentralized FL has yet to be characterized for strongly heterogeneous network settings involving time-varying directed networks~\cite{yuan2024decentralized,li2024decentralized}.

\begin{figure}[t!]
\centering
\includegraphics[width=0.48\textwidth]{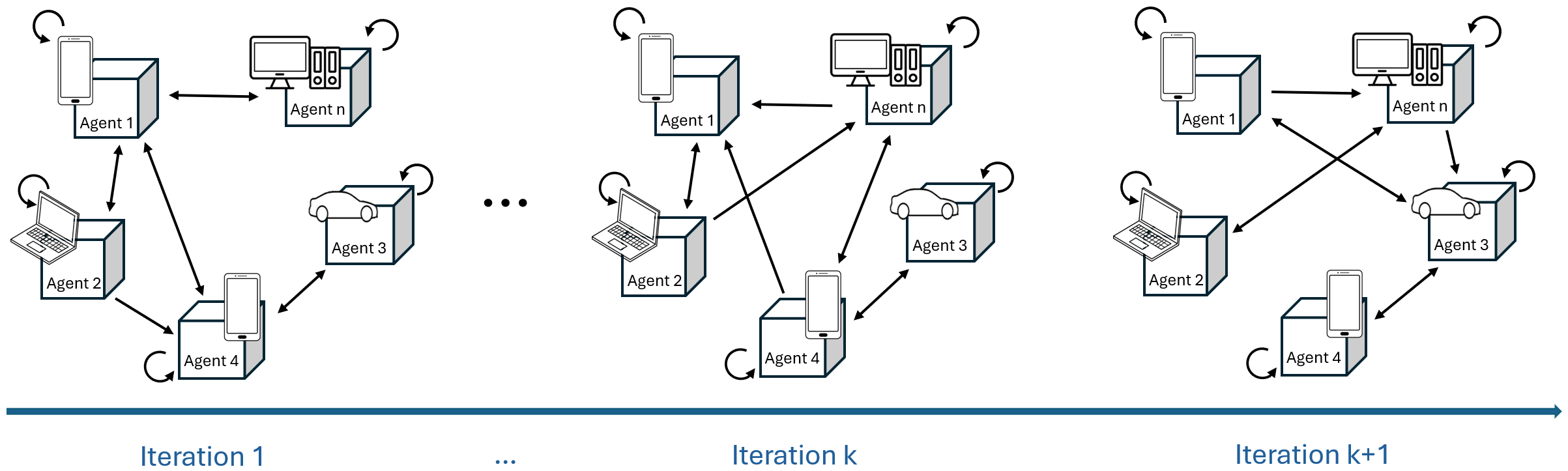}
\caption{Our problem further assumes that directed links in decentralized FL are time-varying. As the time iteration changes, directed links among network agents may shift as a result.}  \label{fig:Time-varying-Directed}
\vspace{-0.4cm}
\end{figure}

We thus aim to study the convergence of decentralized learning and optimization for heterogeneous agents over time-varying directed networks.
% which are essential in practical applications but have been largely overlooked in the literature~\cite{XXX,XXX}. 
In such networks, communication among agents can be dynamic (see Fig.~\ref{fig:Time-varying-Directed}) as a result of user mobility, limited bandwidth, failures and/or straggler effects, which themselves arise from heterogeneous device hardware and networking capabilities. 
Additionally, agent-to-agent links may be unilateral, i.e., directed edges, as observed in applications like UAV networks~\cite{XiaoUAVFL2021}, sensor networks~\cite{intanagonwiwat2003directed}, and heterogeneous wireless systems~\cite{jaradat2017capacity}, as links may have varying power levels and/or bandwidth, posing additional challenges to convergence. %further complicating convergence.
Moreover, we consider data-level heterogeneity, which is characterized by non-identically and independently distributed (non-i.i.d.) training data across agents, and batch gradients, as large datasets or complex ML models may prevent full batch gradient descent. 
In particular, the above forms of data-level heterogeneity may yield local ML models that diverge from the global optimum and noisy gradients, both of which hinder convergence of decentralized FL.  
To this end, we propose a decentralized stochastic gradient tracking optimization algorithm, \DSGTmTV, designed to effectively address these complexities and enhance convergence.

\subsection{Literature Review}

% The literature extensively covers decentralized optimization, emphasizing 
Existing literature on decentralized learning and optimization focuses on consensus-based protocols and gradient descent methods for addressing both deterministic problems \cite{Nedic2011, AngeliaCommunicationTopology2018} and stochastic problems \cite{Nedic2014, lian2017can}. 
Specifically for general distributed learning and FL type scenarios,
% Within the domain of ML, 
decentralized algorithms have garnered significant attention, primarily focusing on addressing undirected or weight-balanced communications~\cite{lian2017can, Tiancheng2022, LuDFL2023, Chellapandi2023, KoloskovaDML}, with results on directed communications being rather restrictive \cite{lalitha2018fully, NguyenAutonomousDriving2022}. 
% \cite{ParasnisDFL2023} explores semi-decentralized federated learning, integrating the traditional device-to-server (D2S) interaction model with multiple D2D clusters characterized as time-varying and directed communication graphs.
To address the challenge associated with data heterogeneity, \textit{gradient tracking} \cite{Tsianos2012,Qu2017, nedic2017achieving,xi2018add, nedic2015distributed,KoloskovaDML} has recently emerged as a promising approach.
% For directed communications, several gradient tracking methods have been proposed \cite{nedic2017achieving, xi2018add, Tsianos2012, nedic2015distributed}, some relying on constructing a column-stochastic weight matrix that requires knowledge of agents' out-degree, while others, such as \cite{xi2018linear, Xin2019FROSTFastRO}, utilize only row-stochastic weights.
For directed communications, several gradient tracking methods have been proposed, including Push-DIGing \cite{nedic2017achieving}, ADD-OPT \cite{xi2018add}, and subgradient-push \cite{Tsianos2012, nedic2015distributed}, which rely on constructing a column-stochastic weight matrix that requires knowledge of each agent's out-degree. Other algorithms, such as \cite{xi2018linear} and FROST \cite{Xin2019FROSTFastRO}, only use row-stochastic weights. 
However, these algorithms may have stability issues, as they typically require separate iterations for Perron eigenvector estimation \cite{pshi21}. 
% typically require separate iterations for Perron eigenvector estimation, which can cause stability issues. 

Recent works on the $AB$/Push-Pull algorithm \cite{xin2018linear, pshi21, Saadatniaki2020, Angelia2022AB} have eliminated the need for estimating the Perron eigenvector by utilizing both row- and column-stochastic weights.
% have addressed this stability issue by utilizing both row- and column-stochastic weights and eliminating the need for eigenvector estimation. 
These approaches demonstrate superior performance compared to existing algorithms that utilize only row- or column-stochastic weight matrices, especially when the problem is ill-conditioned and the communication networks are not well balanced, as discussed in \cite{pshi21} and \cite{Angelia2022AB}. The effectiveness of these methods is attributed to their ability to mitigate critical bottlenecks caused by nodes with minimal in-degree or out-degree neighbors in significantly unbalanced directed graphs. Follow-up works, such as \cite{Saadatniaki2020} and \cite{Angelia2022AB}, further establish linear convergence when the graphs are time-varying.
We further investigate such classes of algorithms in networks with data-level heterogenity and heterogeneous stochastic gradients.
% While reference \cite{Saadatniaki2020} derives linear convergence using Euclidean norms,  \cite{Angelia2022AB} provides an improved analysis using time-varying weighted average and weighted norms. 

Concurrently, the heavy-ball method \cite{POLYAK19641} and Nesterov momentum \cite{nesterov2003introductory} have been integrated into various distributed gradient-based algorithms. 
Notably, the linear convergence for the $AB$/Push-Pull method with a heavy-ball momentum term is proven in \cite{Xin2020Heavy}, while \cite{Huaqing2021} proposes a double-accelerated method by incorporating both momentum terms. For time-varying communication networks, \cite{Xiasheng2021} suggests utilizing heavy-ball and Nesterov techniques to accelerate the Push-DIGing algorithm \cite{nedic2017achieving}. However, the double-accelerated methods \cite{Huaqing2021, Xiasheng2021} require identical values for the heavy-ball and Nesterov acceleration parameters. 
Conversely, reference \cite{Nguyen2023AccAB} demonstrates linear convergence when the heavy-ball and Nesterov momentum parameters can be non-identical.

The above-mentioned works primarily focus on deterministic approaches with exact gradient information, while some other works have explored stochastic approaches to handle noisy gradients \cite{QureshiSADDOPT, pu2021stochastic, XinSahuKhanKar2019, HU2021107894}. In particular, the S-ADDOPT algorithm proposed in \cite{QureshiSADDOPT} is the stochastic version of ADD-OPT, while algorithms DSGT \cite{pu2021stochastic} and \SAB~\cite{XinSahuKhanKar2019} are stochastic versions of $AB$/Push-Pull. In \cite{HU2021107894}, two stochastic algorithms for convex problems, SAGA-UDN and SVRG-UDN, which utilize uncoordinated stepsizes, are studied. However, their convergence is demonstrated solely through numerical examples. %Reference \cite{Wang2021} integrates Nesterov momentum for accelerated optimization of non-convex and non-smooth problems.
Importantly, all these works consider static communication networks and do not address scenarios where the P2P communication links change over time. To the best of our knowledge, existing research has yet to simultaneously address all of the aforementioned challenges while accommodating time-varying directed networks as our work does.

% \textbf{Contribution.} 
% \ella{
\subsection{Outline and Summary of Contributions}
\subsubsection{\textbf{Formulation of \DSGTmTV}} 
Our proposed algorithm for decentralized learning and optimization considers stochastic gradients, uncoordinated (i.e., heterogeneous) model training hyper-parameters (e.g., stepsizes/learning rates and momentum parameters), and dynamic, time-varying directed networks, which, when combined together, results in a general but highly complex problem.
To address this, \DSGTmTV~leverages gradient tracking to mitigate global variance resulting from agent-level heterogeneity (see Section~\ref{sec-GradientTracking}), and moreover, improves convergence rates via heavy-ball momentum. 
Simultaneously, to overcome time-varying directed links and remove the need for Perron eigenvector estimation, we employ a two-phase update mechanism, pulling decision variables from neighbors and pushing gradients to neighbors using time-varying row- and column-stochastic weight matrices, respectively. 
As our method allows for agents to have heterogeneous local hyper-parameters,
% (i.e., heavy-ball momentum, learning rates, etc.), 
\DSGTmTV~is thus equipped to function in fully decentralized settings. 

% This paper proposes the \DSGTmTV~algorithm for decentralized \textit{stochastic} optimization, where computational tasks are distributed across agents, each of which typically possesses local data or information. Each agent has access to a stochastic first-order oracle to obtain gradient estimates with zero mean and bounded variance noise. 
% The \DSGTmTV~algorithm utilizes gradient tracking to mitigate the global variance resulting from the heterogeneity of the nodes' data (see also Section~\ref{sec-GradientTracking}) and incorporates heavy-ball momentum for accelerated performance. To facilitate directed communication networks, a two-phase update mechanism is employed, which involves pulling information about the decision variable from neighbors and pushing information about the gradients to neighbors using row- and column-stochastic weights, respectively. Remarkably, unlike deterministic algorithms  \cite{Qu2017,nedic2017achieving,xi2018add,Saadatniaki2020,Angelia2022AB,Nguyen2023AccAB}, as well as relevant works on stochastic algorithms \cite{XinSahuKhanKar2019,pu2021stochastic,Nguyen2023SABTV}, which assume identical constant stepsize values, agents' stepsizes (or learning rates) and momentum parameters in our proposed algorithm are deliberately made uncoordinated. This enables \DSGTmTV~to function in fully decentralized settings.}
% This approach enables each agent the autonomy to independently determine these values in a decentralized manner, thereby eliminating the need for centralized coordination.

\subsubsection{\textbf{Linear convergence}} 
We rigorously prove that \DSGTmTV~achieves \textit{linear convergence} to the global optimum with exact gradient information and to a neighborhood of the global optimum with stochastic gradients, even with uncoordinated hyper-parameters. Moreover, our convergence bounds 
% extend to stochastic and decentralized algorithms over time-varying directed networks, and they 
capture explicit relationships for training hyper-parameters in terms of cost functions, mixing matrices, and graph connectivity structures. The incorporation of heavy-ball momentum adds complexity to the analysis, which we navigate by tracking the expected state difference. 
To address the challenge of time-varying graphs, our analysis employs time-varying stochastic weights (associated with the weight matrices) to construct Lyapunov functions and establish \textit{one-step} contraction relations for the error terms.
Importantly, our bounds for uncoordinated stepsizes are independent of any notion of heterogeneity, allowing agents to optimize local learning rates independently without affecting global convergence guarantees.
In contrast, existing works~\cite{Xu2018,LU2018516,nedic2017geometrically}, which develop bounds for uncoordinated stepsizes, impose restrictions on the possible stepsizes across network agents.
% This contrasts with much of the existing works, which may employ uncoordinated stepsizes but sets bounds based on stepsize heterogeneity defined as the relative deviation of stepsizes from their average in \cite{Xu2018}, and as the ratio of the largest to the smallest stepsize in \cite{LU2018516} and \cite{nedic2017geometrically}. 
% Consequently, an inherent trade-off exists between the tolerance for stepsize heterogeneity and the achievable largest stepsize. 
% While our approach is inspired by \cite{Xin2019FROSTFastRO} and \cite{pshi21}, which focus on deterministic algorithms in static networks, our analysis further develops convergence for stochastic and decentralized algorithms over time-varying networks. 
% Additionally, we show that momentum parameters can also be uncoordinated and simultaneously provide explicit bounds for both stepsizes and momentum parameters based on the properties of cost functions, mixing matrices, and graph connectivity structures.}

\subsubsection{\textbf{Generality of \DSGTmTV}} 
The proposed methodology encompasses several existing methods in decentralized learning and/or optimization as special cases. By adjusting various control parameters of \DSGTmTV~and settings, we can recover algorithms from previous works \cite{xin2018linear, Saadatniaki2020, pu2021stochastic, Angelia2022AB}. These mechanics are discussed in detail in Remark~\ref{remark:SpecialCaseSABTV}.

% Our analysis improves and generalizes the analyses conducted in \cite{xin2018linear, pshi21, Saadatniaki2020,Angelia2022AB,XinSahuKhanKar2019, pu2021stochastic, Nguyen2023AccAB}, all of which are algorithms that represent special cases of \DSGTmTV. Refer to Remark~\ref{remark:SpecialCaseSABTV} for more details.}

\subsubsection{\textbf{Experimental evaluation}} We verify the generality and efficacy of \DSGTmTV~relative to baselines from literature on both image classification and natural language processing tasks in Section~\ref{sec:simulation}. Our proposed algorithm achieves both faster convergence and better final performance than existing methods, highlighting the veracity of our theoretical results.

% Some deterministic algorithms integrate the use of uncoordinated stepsizes, as demonstrated in works such as \cite{nedic2017geometrically} and \cite{Xu2018} for undirected graphs with doubly-stochastic weights, as well as in \cite{Xin2019FROSTFastRO} and \cite{LU2018516,pshi21} for directed graphs with row-stochastic and column-stochastic weights, respectively. 

% This distinguishes our work from prior literature on stochastic distributed algorithms \cite{QureshiSADDOPT,Wang2021,pu2021stochastic,XinSahuKhanKar2019,Nguyen2023SABTV,HU2021107894}.

This paper is structured as follows. Section~\ref{sec:formu} presents the problem formulation, with the decentralized algorithm introduced in Section~\ref{sec:algo}. Section~\ref{sec:conv_results} provides auxiliary results and convergence analysis. Section~\ref{sec:simulation} numerically evaluates the performance of the proposed algorithm. Section~\ref{sec:conc} concludes with key points. 
% The Appendix outlines notational conventions, additional foundational results, and proofs.

\subsection{Notational Conventions}
\label{ss:notation}
All vectors are viewed as column vectors unless otherwise stated. The transpose of a vector $u\in \re^n$ is denoted by $u^\T$. The vector with all entries as $0$ is denoted by $\zero$, and with all entries as $1$ by $\one$. The $i$-th entry of a vector $u$ is denoted by $u_i$, while it is denoted by $[u_k]_i$ for a time-varying vector $u_k$. We denote $\min(u)=\min_i u_i$ and $\max(u)=\max_i u_i$. A non-negative vector is called stochastic if its entries sum to $1$.  

We use $A_{ij}$ to denote the $ij$-th entry of a matrix $A$, and $[A_k]_{ij}$ when the matrix is time-dependent. The notation $A \!\le\! B$ is used when $A_{ij} \!\le\! B_{ij}$ for all $i, j$, where $A$ and $B$ are matrices of the same dimension. A matrix $A$ is nonnegative if all its entries are nonnegative and ${\min}^{+}(A)$ denotes the smallest positive entry of $A$. A nonnegative matrix $A\in\mathbb{R}^{n\times n}$ is considered to be row-stochastic if $A\mathbf{1}=\mathbf{1}$, and a nonnegative matrix $B\in\mathbb{R}^{n\times n}$ is considered to be column-stochastic if $\mathbf{1}^{\T} B=\mathbf{1}^{\T}$. $\Diag(u)$ denotes the diagonal matrix with vector $u$ on its diagonal. The identity matrix is denoted by $\mbi$. 

Given a positive vector $\ba=(a_1,\ldots,a_n)\in\re^n$, we denote:
\begin{center}
    $\la \bu,\bv\ra_{\ba}=\sum_{i=1}^m a_i\la u_i,v_i \ra $ and $\|\bu\|_{\ba}=\sqrt{\sum_{i=1}^m a_i\|u_i\|^2},$
\end{center}
where $\bu\!:=\![u_1,\ldots,u_n]^\T, \bv \!:=\![v_1,\ldots,v_n]^\T \!\!\in\!\! \re^{n\times p}$, and $u_i,v_i\!\in\!\re^p$. When $\ba = \one$, we write $\la \bu,\bv\ra$ and $\|\bu\|$. We have
\begin{align}\label{eq-NormIneq}  \|\bu\| \le \tfrac{1}{\sqrt{\min(\ba)}}\|\bu\|_{\ba} \text{and} |\la \bu,\bv \ra|_{\ba} \le\|\bu\|_{\ba}\|\bv\|_{\ba}.
\end{align}

We let $[n]=\{1,\ldots,n\}$ for an integer $n\ge 1$. Consider a directed graph $\G=([n],\mathcal{E})$, specified by the set of edges $\mathcal{E}\subseteq [n]\times[n]$ of ordered pairs of nodes. Let $\mathsf{D}(\bbG)$ denote the diameter and $\mathsf{K}(\bbG)$ denote the maximal edge-utility of the graph $\bbG$, as defined in Definition 2.1 and Definition 2.2 in \cite{Angelia2022AB}. Additionally, the in-neighbor and out-neighbor set for every agent $i$ are defined, as follows:
\[\cNini=\{j\in[n]|(j,i)\!\in\!\mathcal{E}\} \!\!\!\text{and}\!\!\! \cNouti=\{l\in[n]|(i,l)\!\in\!\mathcal{E}\}.\] 
A directed graph $\G$ is {\it strongly connected} if there is a directed path from any node to all other nodes in $\G$. 
% We use $\mathsf{D}(\G)$ and $\mathsf{K}(\G)$ to denote the diameter and the maximal edge-utility of a strongly connected directed graph $\G$, respectively, as defined in Definition 2.1 and Definition 2.2 in \cite{Angelia2022AB}.

When the graph varies over time, we use a subscript $k$ to indicate the time instance. For example, $\mathcal{E}_k$  denotes the edge-set of a graph $\G_k$, $\Nikin$ and $\Nikout$ denote the in-neighbors and the out-neighbors of a node $i$, respectively, at time $k$.

\section{Problem Formulation}
\label{sec:formu}
Consider a distributed or decentralized (federated) learning task performed across a network of $n$ agents (e.g., devices, clients or learners). Each agent $i \in [n]$ possesses a local dataset $\mathcal{D}_i$, consisting of  $m_i$ data samples $\{\xi_{ij} \!=\! (\mathscr{b}_{ij},\ell_{ij})\}_{j=1}^{m_i}$, where $\mathscr{b}_{ij} \!\in\!\re^p$ represents the feature vector, and $\ell_{ij} \!\in\! \re$ is the corresponding label. Agents independently process their local data and communicate with neighboring agents to collaboratively train a global model.

For a global model parameterized by $x\in \re^p$, let $\mathcal{L}(x;\xi)$ denote the loss incurred when predicting the sample $\xi$ using the model $x$. The average loss for $x$ over the local dataset $\mathcal{D}_i$ of agent $i$ is given by $f_i: \re^p \to \re$, as follows:
\[f_i(x) = \E_{\xi_i \sim \mathcal{D}_i} \mathcal{L}(x,\xi_i),\] 
where $\xi_i$ represents a data sample drawn from $\mathcal{D}_i$.

The overall objective in decentralized FL is to determine the model parameters $x$ that minimizes the global loss function, which aggregates the local losses from all agents. This corresponds to solving the following optimization problem: 
\begin{equation} \label{eq-problem}
\min_{x\in \re^p}~ f(x)=\frac{1}{n}\sum\limits_{i=1}^n f_i(x).
\end{equation}
% Therefore, the problem is formulated as a collaborative optimization task as described in \eqref{eq-problem}. 

% \ella{The form of the loss function $\mathcal{L}(x;\xi)$  varies depending on the specific ML task. Some common examples include Mean Squared Error (MSE), which is used in regression tasks, and hinge loss, employed to optimize Support Vector Machines (SVMs) by penalizing predictions that fall within a specific margin. In classification, cross-entropy loss evaluates the model's prediction probabilities across different classes, while logistic loss is specifically for binary classification tasks.}

\subsection{General Setting: Decentralized Stochastic Optimization} 
\label{sec:DecentralizedOpt}
To put the above problem into a more general context, consider a system of $n$ agents connected via a communication network. These agents aim to collaboratively solve the stochastic optimization problem in \eqref{eq-problem} for a global decision $x$. Here, each function $f_i: \mathbb{R}^p \rightarrow \mathbb{R}$ represents a local cost function known only to agent $i$. Agent $i\in[n]$ obtains a gradient estimate $g_i(x, \xi_i)$ through a stochastic first-order oracle. The estimate is a noisy approximation of the true gradient $\nabla f_i(x)$, influenced by a random variable $\xi_i$ that captures the inherent variability in data sampling or measurement errors.

This problem is relevant to various applications beyond FL. For instance, in distributed sensor networks \cite{intanagonwiwat2003directed}, sensors aggregate local measurements to minimize global estimation error, with stochastic gradients reflecting random measurement errors or variability in sensor readings. In multi-agent systems, such as autonomous vehicles \cite{NguyenAutonomousDriving2022,XiaoUAVFL2021}, agents coordinate their actions to optimize overall system performance. Here, stochastic factors include uncertainties in sensor data, environmental conditions, and agent behaviors, all contributing to the variability in gradient estimates.

\subsection{Assumptions}
We make the following assumptions regarding the local objective functions: 
\begin{assumption} \label{asm-functions}
Each $f_i$ is continuously differentiable and has $L$-Lipschitz continuous gradients, i.e., for some $L>0$,
\begin{equation*}
\|\nabla f_i(x)-\nabla f_i(y)\|\le L \|x-y\|,
\quad\hbox{for all $x,y\in \re^p$}.
\end{equation*}
\end{assumption}
	
\begin{assumption}\label{asm-strconv}
The average-sum function $f=\frac{1}{n}\sum_{i=1}^n f_i$ is $\mu$-strongly convex, i.e., for some $\mu>0$, 
\[\la \nabla f(x)-\nabla f(y),x-y\ra \ge \mu\|x-y\|^2\quad\hbox{for all $x,y\in \re^p$}.\]
\end{assumption}
Assumption~\ref{asm-strconv} implies that problem~\eqref{eq-problem} has a unique optimal solution $x^*$ such that $x^*=\argmin_{x\in\re^p} f(x)$.

\begin{remark}
Assumption~\ref{asm-strconv}, similar to~\cite{Saadatniaki2020, Angelia2022AB, Nguyen2023AccAB, Nguyen2023SABTV,ghosh2022efficient,avdiukhin2021federated}, is more general than requiring strong convexity for all local cost functions $f_i$, which is common in relevant works~\cite{xin2018linear, pshi21,pu2021stochastic,XinSahuKhanKar2019}. Assumption~\ref{asm-strconv} is satisfied, for instance, when $f_i$ is convex $\forall i$ and any $f_i$ has strong convexity.
\end{remark}

% \begin{remark}
% Assumption~\ref{asm-strconv}, as also employed in \cite{Saadatniaki2020, Angelia2022AB, Nguyen2023AccAB, Nguyen2023SABTV}, is less restrictive compared to the strong convexity requirement for every local cost function $f_i$ assumed in related works \cite{xin2018linear, pshi21,pu2021stochastic,XinSahuKhanKar2019}. Assumption~\ref{asm-strconv} is satisfied, for instance, when each $f_i$ is convex and at least one of them exhibits strong convexity.
% \end{remark}

The noisy gradient estimates satisfy the following standard assumptions \cite{pu2021stochastic}:
\begin{assumption} \label{asm-SFO}
For all $i\in[n]$ and $x\in\re^p$, each random vector $\xi_i\in\re^m$ is independent, and
\begin{flalign*}
&\text{(a)}\!\!\E[g_i(x,\xi_i)|x] = \nabla f_i(x),&&\\
&\text{(b)}\!\! \E[\|g_i(x,\xi_i)-\nabla f_i(x)\|^2|x] \le \sigma^2, \!\!\text{for some $\sigma>0$.}&&\!\!\!\!
\end{flalign*}
\end{assumption}

\begin{remark}
Assumption~\ref{asm-SFO} is applicable in several fields such as online learning, reinforcement learning, generative models, and parameter estimation in signal processing and communication systems \cite{Bottou2018,Jakovetic2018,pu2021stochastic}. The bounded variance assumption, however, can be relaxed, as demonstrated in \cite{Bottou2018}.
\end{remark}

% Assumption~\ref{asm-SFO} is a widely applicable assumption that is particularly pertinent  in the context of online distributed learning, reinforcement learning, generative models, and parameter estimation  in signal processing and communication systems. %Another example where 
% %This assumption is also applicable to
% It is also relevant in simulation-based optimization, where gradient estimation often produces noise owning to various sources such as simulation model discretization, measurement errors, incomplete convergence, and finite sample size for Monte-Carlo methods \cite{pu2018distributed}.
% Assumption~\ref{asm-SFO} holds true in many on-line distributed learning problems, where $f_i(x) = \E_{\xi_i} [F_i(x, \xi_i)]$ denotes the expected loss function agent $i$ wishes to minimize, while independent samples $\xi_i$ are gathered continuously over time. For another example, in simulation-based optimization, the gradient estimation often incurs noise that can be due to various sources, such as the discretization of the simulation model or measurement errors, incomplete convergence, and finite sample size for Monte-Carlo methods.

\subsection{Communication Networks} 
\label{sec:Communication}
The agents aim to collaboratively train the model $x^*$ through local computations and information exchange with neighboring agents. At each iteration $k$, the communication network is represented by the directed graph $\bbG_k=([n],\mathcal{E}_k)$, where the vertices correspond to the $n$ agents and the edges indicate their connectivity. Specifically, an edge $(j,i) \in \mathcal{E}_k$ signifies that agent $i$ receives information from agent $j$ at iteration $k$. Fig.~\ref{fig:FL_decentralized} illustrates the information flow among devices within a decentralized FL framework at a given iteration.

We consider the challenging scenario where agent communication occurs through a sequence of time-varying directed networks, governed by the sequence of directed graphs $\{\bbG_k\}$. In these networks, agent-to-agent links may be unilateral due to factors such as varying power levels, bandwidth limitations, or environmental conditions like interference and obstructions \cite{intanagonwiwat2003directed,jaradat2017capacity}. The communication among agents can be dynamic \cite{NguyenAutonomousDriving2022,XiaoUAVFL2021}, influenced by limited bandwidth, user mobility, intermittent connectivity, network failures, and straggler effects. These dynamics arise from heterogeneous device hardware, differing networking capabilities, and fluctuating network loads, all contributing to varying communication patterns.
 Fig.~\ref{fig:Time-varying-Directed} illustrates how the communication topology changes across consecutive iterations. Note that time-varying directed networks include static networks  or undirected networks, considered in prior studies \cite{lian2017can, Tiancheng2022, LuDFL2023, Chellapandi2023, KoloskovaDML}, as special cases.

Formally, we consider the following assumption regarding the communication graphs:
\begin{assumption} \label{asm-graphs}
For each $k \ge 0$, the directed graph $\G_k$ is strongly connected and has a self-loop at every node $i\in[n]$.
\end{assumption}

\begin{remark}
Assumption~\ref{asm-graphs} can be relaxed by considering 
% a sequence of graphs that are $C$-strongly connected
a sequence of $C$-strongly connected graphs, i.e., for every $k\ge 0$, there exists an integer $C\ge 1$ such that the graph formed by the edge set $\mathcal{E}^C_k=\bigcup_{i=kC}^{(k+1)C-1}\mathcal{E}_i$ is strongly connected.  
\end{remark}

We parameterize communication through mixing matrices. Specifically, the information exchange at time $k$ is governed by two non-negative weight matrices, $A_k$ and $B_k$, that {\it align} with the topology of the graph $\bbG_k$, in the following sense:
\begin{align}\label{eq-alignA}
\!\!\![A_k]_{ij}\! >0,~\forall j\!\in\!\Nikin\!\cup\!\{i\}&;~\,[A_k]_{ij}\!=0,~\forall j\!\not\in\!\Nikin\!\cup\!\{i\},\!\!\\
\label{eq-alignB}
\!\!\![B_k]_{ji}\! >0,~\forall j\!\in\!\Nikout\!\cup\!\{i\}&;~\, [B_k]_{ji}\!=0,~\forall j\!\not\in\!\Nikout\!\cup\!\{i\}.\!\!
\end{align}

Each agent $i$ independently decides the entries $[A_k]_{ij}$ for their in-neighbors $j\in\Nikin$, while agent $j\in\Nikin$ determines the value $[B_k]_{ij}$. We make the following assumptions regarding the matrices $A_k$ and $B_k$:

\begin{assumption} \label{asm-amatrices}
For each $k \ge 0$, the matrix $A_k$ is row-stochastic, i.e., $A_k\1=\1$, and compatible with the graph $\bbG_k$ in the sense of relation~\eqref{eq-alignA}. Moreover, there exists a scalar $a\!>\!0$ such that $\min^+(A_k)\ge a$ for all $k\ge0$.
\end{assumption}

\begin{assumption}\label{asm-bmatrices}
For each $k \ge 0$, the matrix $B_k$ is column-stochastic, i.e., ${\1}^{\T} B_k={\1}^{\T}$, and compatible with the graph $\bbG_k$ in the sense of relation~\eqref{eq-alignB}. Moreover, there exists a scalar $b>0$ such that $\min^+(B_k)\ge b$ for all $k\ge0$.
\end{assumption}

\section{The \DSGTmTV~Algorithm} \label{sec:algo}
This section introduces the Momentum-Accelerated \underline{D}ecentralized \underline{S}tochastic \underline{G}radient \underline{T}racking for \underline{T}ime-\underline{V}arying Directed Networks (\DSGTmTV) algorithm. This algorithm represents a consensus-based decentralized learning framework incorporating gradient tracking and momentum-based acceleration (``$m$"), to collaboratively solve problem~\eqref{eq-problem} and obtain the optimal global solution or model $x^*$, where agents communicating over time-varying directed networks.

\begin{algorithm}[t!]
\KwIn{Local cost function $f_i(x)$ for all agent $i \in [n]$.}
\KwOut{Optimal global decision $x^*$.}

\textbf{Initialize:} Every agent $i \in [n]$ sets arbitrary initial vectors $x_{-1}^i, x_0^i \in \mathbb{R}^p$, local gradient estimate $y_0^i = g_i(x_0^i, \xi_0^i)$, stepsize $\alpha_i$, and heavy-ball momentum parameter $\beta_i$.

\For{$k=0, 1, 2, \ldots$}
  {
    \SetKwBlock{Begin}{Each agent $i \in [n]$ do}{end}
    \Begin
    {
        \!\!\!\textit{Communication Step:}\\
        Choose the weights $[B_k]_{ji},j\in\cNoutik$;\\
        Send $x_k^i$ and $[B_k]_{ji}y_k^i$ to  out-neighbors $j\in\cNoutik$;\\
        Receive $x_{k}^j$ and $[B_k]_{ij}y_k^j$ from in-neighbors $j\!\in\!\cNinik$;\\
        Choose the weights $[A_k]_{ij},j\in\cNinik$;

        \!\!\!\textit{Local State Update:}
        \begin{equation*}
        x_{k+1}^i = \sum_{j=1}^n [A_k]_{ij} x_k^j - \alpha_i y_k^i + \beta_i (x_k^i - x_{k-1}^i),
        \end{equation*}
        
        \!\!\!\textit{Obtain Gradient Estimate:}\\
        Query the stochastic first-order oracle for a gradient estimate $g_i(x_{k+1}^i, \xi_{k+1}^i)$;

        \!\!\!\textit{Gradient Tracking Update:}
        \begin{equation*}
        y_{k+1}^i = \sum_{j=1}^n [B_k]_{ij} y_k^j + g_i(x_{k+1}^i, \xi_{k+1}^i) - g_i(x_k^i, \xi_k^i).
        \end{equation*}
    }
  }
\caption{The \DSGTmTV~Algorithm} \label{alg:DSGTmTV}
\end{algorithm}

Each agent $i\in[n]$ maintains a local model $x_k^i\in \re^p$, which is an estimate of the ``global model", and a local update direction $y_k^i\in\mathbb{R}^p$ at iteration $k$. These variables are updated over time. Each agent $i$ initializes their updates using arbitrary vectors $x_{-1}^i,x_0^i$ and $y_0^i=g_i(x_0^i,\xi_0^i)$, without the need for coordination among the agents. At time $k\ge 0$, every agent $i$ sends its vector $x_k^i$ and a scaled direction $[B_k]_{ji}y_k^i$ to its out-neighbors $j\in\Nikout$ and receives these vectors sent by its in-neighbors $j\in\Nikin$.
Upon the information exchange, every agent $i$ performs a local state update and a gradient tracking update, as follows:
\begin{subequations}\label{eq-met}
	\begin{align}
    % & \text{Local state update:} \nonumber\\
	&x_{k+1}^i =  \sum_{j=1}^n[A_k]_{ij}x_{k}^j  - \a_i y_k^i+\!\!\!\!\!\!\!\!\underbrace{\b_i\left(x_{k}^i-x_{k-1}^i\right),}_{\text{heavy-ball momentum}} \!\!\!\!\label{eq-x}\\
    % & \text{Gradient tracking update:} \nonumber\\
	&y_{k+1}^i =  \sum_{j=1}^n [B_k]_{ij}y_{k}^j + \underbrace{g_i(x_{k+1}^i,\xi_{k+1}^i) -g_i(x_k^i,\xi_k^i)}_{\text{correction}},\label{eq-y}
	% &\hbox{with arbitrary $x_0^i\in\mathbb{R}^p$ and $y_0^i=g_i(x_0^i,\xi_0^i)$}.
	% \label{eq-init}
	\end{align}
\end{subequations}
where $\a_i>0$ is a constant stepsize, and $\b_i\ge 0$ is the heavy-ball momentum parameter chosen by agent $i$ for all $i\in [n]$. The procedure is outlined in Algorithm~\ref{alg:DSGTmTV}.

\begin{remark}[Generality of \DSGTmTV] \label{remark:SpecialCaseSABTV}
In our proposed algorithm, model training hyper-parameters such as learning rates/stepsizes $\a_i$ and momentum parameters $\b_i$ (for all $i\in[n]$), are fully uncoordinated (i.e., heterogeneous across agents). By tuning these hyper-parameters and considering special cases of the communication graphs $\G_k$, we can recover several existing algorithms from previous works. Specifically, when $\a_i=\a$ and $\b_i=0$, for all $i\in [n]$, we recover the DSGT algorithm  \cite{pu2021stochastic} for a static (i.e., $\G_k = \G$ for all $k\ge 0)$ undirected communication graph, and the \SAB~algorithm \cite{XinSahuKhanKar2019} for a static directed communication graph. Moreover, the algorithms presented in \cite{xin2018linear, pshi21, Saadatniaki2020, Angelia2022AB, Xin2020Heavy, Nguyen2023AccAB} represent deterministic versions of \SABTV, assuming exact gradient information ($\sigma^2 = 0$). Specifically, references \cite{xin2018linear, pshi21} focus on undirected graphs, while references \cite{Saadatniaki2020, Angelia2022AB,Xin2020Heavy, Nguyen2023AccAB} consider directed graphs, with the algorithms in \cite{Xin2020Heavy, Nguyen2023AccAB} further incorporating momentum terms ($\b > 0$).
\end{remark}

\subsection{Compact Form}
To express the updates from \eqref{eq-met} in a compact form, let $\bx_k=[x_k^1,\ldots,x_k^n]^\T\in\re^{n\times p}$ be the matrix that stacks all the local models of the agents at iteration $k$. Similarly, let $\by_{k}=[y_k^1,\ldots,y_k^n]^\T\in\re^{n\times p}$ represent the matrix stacking all local update directions. The stochastic gradient information from all agent is represented by $\bg(\bx_k,\bxi_{k})=[g_1(x_{k}^1,\xi_{k}^1),\ldots,g_n(x_{k}^n,\xi_{k}^n)]^\T\in\re^{n\times p}$, where $\bxi_k=[\xi_k^1,\ldots,\xi_k^n]^\T\in\re^{n\times m}$ is a matrix containing the random variables associated with each agent's gradient estimate at iteration $k$. With this notation, \eqref{eq-met} can be written, as follows:
\begin{subequations}\label{eq-met-comp}
	\begin{align}
	&\bx_{k+1} =  A_k\bx_k  - \Da \by_{k}+\Db(\bx_{k}-\bx_{k-1}),\label{eq-x-comp}\\
	&\by_{k+1} =  B_k\by_{k} + \bg(\bx_{k+1},\bxi_{k+1}) -\bg(\bx_k,\bxi_{k}),\label{eq-y-comp}
	\end{align}
\end{subequations}
where $\balpha=(\a_1,\ldots,\a_n)^\T$ and $\bbeta=(\beta_1,\ldots,\beta_n)^\T$ 
define the vectors of the stepsizes and the momentum parameters, respectively, and $\Da=\Diag(\balpha)$ and $\Db=\Diag(\bbeta)$ are the corresponding diagonal matrices. We will also employ the following notation in our subsequent analysis:
\begin{align}\label{eq-quant-avg1}
&\bby_k=\frac{1}{n}\one_n^\T\by_{k},~~ \bar{g}_k(\bx_k,\bxi_k)=\frac{1}{n}\one_n^\T\bg(\bx_k,\bxi_k),\\ 
&h(\bx_k)=\frac{1}{n}\one^\T_n\nabla F(\bx_k) = \frac{1}{n}\sum_{i=1}^n \nabla f_i(x_k^i), \label{eq-quant-avg2}
\end{align}
where $\nabla F(\bx_k)=[\nabla f_1(x_k^1),\ldots,\nabla f_n(x_k^n)]^\T\!\in\!\re^{n\times p}$.

Finally, we denote by $\F_k$ the $\sigma$-algebra generated by the set of random vectors $\{\bxi_0,\ldots,\bxi_{k-1}\}$, and denote by $\E[\cdot|\F_k]$ the conditional expectation given $\F_k$.

We have the following results regarding stochastic gradients.
\begin{lemma}[\!\!\cite{pu2021stochastic}] \label{lem-sumgrad}
Consider the iterates $\{y_k^i\}_{i\in [n], k\ge 0}$ generated by the \DSGTmTV~method in~\eqref{eq-met}. Let Assumption~\ref{asm-functions}, Assumption~\ref{asm-SFO} and Assumption~\ref{asm-bmatrices} hold. For all $k\ge0$, we have
\begin{flalign*}
&\text{(a)} \bby_k=\bar{g}_k(\bx_k,\bxi_k),&& \\
&\text{(b)} \E[\bby_k|\F_k] = h(\bx_k),&& \\
&\text{(c)} \E[\|\bby_k-h(\bx_k)\|^2|\F_k] \le \tfrac{\sigma^2}{n}.&& 
\end{flalign*}
% where $\F_k$ is the $\sigma$-algebra generated by the set of random vectors $\{\bxi_0,\ldots,\bxi_{k-1}\}$ and $\E[\cdot|\F_k]$ is the conditional expectation given $\F_k$.
\end{lemma}
Lemma~\ref{lem-sumgrad}(a) establishes that the average of agents' local directions, $\bby_k$, effectively tracks the estimate of the global direction $\bar{g}_k(\bx_k,\bxi_k)$ through the gradient tracking technique in \eqref{eq-y}. 
The proof uses mathematical induction on $k$, utilizing the column-stochasticity of $B_k$ and the initialization of the $y$-variable. For a comprehensive proof of all statements in Lemma~\ref{lem-sumgrad}, refer to \cite{pu2021stochastic,XinSahuKhanKar2019}.

\subsection{Gradient Tracking} \label{sec-GradientTracking}
One of the challenges in the decentralized setting is ensuring that the updates made by each agent are representative of the global objective, especially when data is distributed heterogeneously across agents. Below, we discuss the gradient tracking technique employed to address this issue.

% To understand the rationale behind incorporating gradient tracking, w
We first revisit the Decentralized Stochastic Gradient Descent (DSGD) algorithm, a widely adopted method in the DFL setting (see \cite{Tiancheng2022} and references therein). Consider a simple case where the communication graph $\G_k=\G$ is static and undirected, and the corresponding mixing matrix $W$ is doubly-stochastic.  In this setting, the update rule for DSGD at each agent $i\in[n]$ is given by
\begin{align}\label{eq:DSGD}
x_{k+1}^i =  \sum_{j=1}^n W_{ij}x_{k}^j  - \a_k g_i(x_k^i,\xi_k^i), \text{for all $k\ge 0.$} 
\end{align}

Reference \cite{Yuan2020} establishes that under a sufficiently small constant stepsize $\a_k = \a$ for all $k\ge 0$, DSGD achieves linear convergence at a rate of $\rho^k$ with $\rho = 1-\mathcal{O}(\alpha\mu)$ towards a neighborhood of the global optimum $x^*$. Furthermore, the steady-state error is characterized by
\begin{align}\label{eq-DSGD-steady-state-error}
\lim\sup_{k\to\infty} &\frac{1}{n}\sum_{i=1}^n\E[\|x^i_k-x^*\|_{2}^2] \\
&= \mathcal{O}\left(\frac{\a\sigma^2}{n\mu} + \frac{L^2}{\mu^2}\frac{\a^2\lambda_M^2\sigma^2}{1-\lambda_M} + \frac{L^2}{\mu^2}\frac{\a^2\lambda_M^2 \mathscr{e}^2}{(1-\lambda_M)^2}\right),\nonumber
\end{align}
where $\lambda_M$ denotes the second largest singular value of $W$, and $\mathscr{e}=\frac{1}{n}\sum_{i=1}^n\|\nabla f_i(x^*)\|^2_2$. Notably, DSGD introduces an additional error term $\frac{L^2}{\mu^2}\frac{\a^2\lambda_M^2 \mathscr{e}^2}{(1-\lambda_M)^2}$. This error becomes significant when $\mathscr{e}$ is large, often due to heterogeneous data distributions among agents. Addressing this term is crucial for improving steady-state performance, especially when it dominates the overall error.

The discussion above highlights the motivation for employing the gradient tracking technique in DSGD to derive DSGT \cite{pu2021stochastic}, a specific instance of \DSGTmTV~under static undirected communication graphs (i.e., $\G_k=\G$ and $A_k = B_k =W$), with $\b_i=0$ for all $i\in [n]$. Specifically, the updates in \eqref{eq-met} reduce to the following form, which defines DSGT:
\begin{subequations}\label{eq-met-DSGT}
	\begin{align}
    &x_{k+1}^i =  \sum_{j=1}^nW_{ij}x_{k}^j  - \a_k y_k^i,\label{eq-x-DSGT}\\
	&y_{k+1}^i =  \sum_{j=1}^n W_{ij}y_{k}^j + g_i(x_{k+1}^i,\xi_{k+1}^i) -g_i(x_k^i,\xi_k^i).\label{eq-y-DSGT}
	\end{align}
\end{subequations}
It has been demonstrated in \cite{pu2021stochastic} that for a sufficiently small constant stepsize $\a_k = \a$ for all $k\ge 0$, DSGT achieves linear convergence at the rate of $\rho^k$ with $\rho = 1-\mathcal{O}(\alpha\mu)$, converging to a neighborhood of $x^*$. The steady-state error is given by
\begin{align}\label{eq-DSGT-steady-state-error}
\lim\sup_{k\to\infty} \frac{1}{n}\sum_{i=1}^n &\E[\|x^i_k-x^*\|_{2}^2] \nonumber\\
&= \mathcal{O}\!\left(\!\frac{\a\sigma^2}{n\mu} \!+ \frac{L^2}{\mu^2}\frac{\a^2\lambda_M^2\sigma^2}{(1-\lambda_M)^3}\!\right)\!.
\end{align}

Unlike DSGD, the steady-state error of DSGT does not depend on $\mathscr{e}$; it completely eliminates the additional error term $\frac{L^2}{\mu^2}\frac{\a^2\lambda_M^2 \mathscr{e}^2}{(1-\lambda_M)^2}$. This improvement is attributed to the use of gradient tracking
% , where the local descent direction $y_k^i$ in DSGT approximates the global stochastic gradient; 
(see Lemma~\ref{lem-sumgrad}(a)). 
  
Consequently, in the deterministic case (i.e., $\sigma^2=0$) and with a constant stepsize, it can be observed from \eqref{eq-DSGD-steady-state-error} and \eqref{eq-DSGT-steady-state-error} that DSGT achieves linear convergence to the global minimum $x^*$, while DSGD converges only to a neighborhood of $x^*$.

In the stochastic case, it should be noted that the majority of DFL algorithms commonly employ DSGD with decaying stepsizes \cite{Tiancheng2022}, for which gradient tracking might not provide a clear advantage over standard DSGD. Specifically, with stepsizes decaying as $\a_k = \mathcal{O}(1/k)$, both DSGD and DSGT converge sublinearly to the exact global minimum $x^*$, and,
\begin{align}\label{eq-DSG-decaying}
    \frac{1}{n}\sum_{i=1}^n\E[\|x^i_k-x^*\|_{2}^2] = \mathcal{O}\left(\frac{1}{k}\right), \text{for all $k\ge 0.$}
\end{align}
However, when the stepsize is constant and in scenarios with significantly heterogeneous data distributions among agents (i.e.,  large $\mathscr{e}$), gradient tracking is preferable. Moreover, while gradient tracking alone does not guarantee exact linear convergence in this context, its combination with variance reduction techniques like SAGA \cite{SAGADefazio2014} can achieve this goal. However, this is beyond the scope of the current work.
% Although this is beyond the scope of the current work, it represents a promising direction for future research.
% It is worth noting that the network dependence of steady-state error in DSGD is better, being $\mathcal{O}((1-\lambda_M)^{-2})$, compared to $\mathcal{O}((1-\lambda_M)^{-3})$ in DSGT. 

\section{Convergence Analysis }\label{sec:conv_results}
This section provides a detailed analysis of the convergence properties of the \DSGTmTV~algorithm over a time-varying directed communication network. We first present preliminary results, followed by the establishment of contraction relations for four critical error quantities. We demonstrate that \DSGTmTV~achieves linear convergence to the exact global optimum when exact gradient information is available ($\sigma^2=0$), and converges linearly in expectation to a neighborhood of the optimal solution when employing stochastic gradient estimates ($\sigma^2>0$), provided that sufficiently small constant stepsizes and momentum parameters are appropriately chosen.

\subsection{Preliminaries}\label{sec:preliminaries}
This section derives several preliminary results necessary for the main convergence analysis, while additional fundamental results from existing works that support the analysis are provided in Appendix~\ref{Appendix:basicre}. Specifically, the following lemmas define the stochastic vector sequences associated with the sequences of stochastic mixing matrices $\{A_k\}$ and $\{B_k\}$. These stochastic vectors will be employed to construct suitable Lyapunov functions associated with the proposed method.

\begin{lemma}[\!\!\cite{nguyen2022distributed}, Lemma 5.4] \label{lem-amatrices}
Let Assumption~\ref{asm-graphs} hold. Consider a row-stochastic matrix sequence $\{\!A_k\!\}$ satisfying Assumption \ref{asm-amatrices}, there exists a sequence $\{\phi_k\}$ of stochastic vectors such that
\begin{align}\label{eq-phik}
\phi_{k+1}^{\T}A_k=\phi_k^{\T},\qquad\hbox{for all $k\ge0$}.
\end{align}
Moreover, $[\phi_k]_i\ge \tfrac{a^n}{n}$ for all $i\in[n]$ and for all $k\ge0$.
\end{lemma}

\begin{lemma}[\!\!\cite{Angelia2022AB}, Lemma 3.4] \label{lem-bmatrices}
Let Assumption~\ref{asm-graphs} hold. Given a column-stochastic matrix sequence $\{B_k\}$ satisfying Assumption~\ref{asm-bmatrices}, define the vector sequence ${\pi_k}$, for $k\ge 0$, as follows:
\begin{align}\label{eq-pik}
\pi_{k+1}=B_k\pi_k,\qquad\hbox{initialized with }\ \pi_0=\tfrac{1}{n}\1.
\end{align} 
Then, the vectors $\pi_k$ are stochastic and  $[\pi_k]_i\ge \tfrac{b^n}{n}$ for all $i\in[n]$ and $k\ge0$.
\end{lemma}

Let $\hat x_k= \sum_{i=1}^n[\phi_k]_i x_k^i$ be the $\phi_k$-weighted average of the iterates $x_k^i$ where $\{\phi_k\}$ is the sequence of stochastic vectors satisfying $\phi_{k+1}^{\T}A_k=\phi_k^{\T}$ (see Lemma~\ref{lem-amatrices}) and $\hat{\bx}_k=\one_n \hat x_k$. The recursive relation for the weighted average ${\hat x_k}$, presented below, will be utilized to track the consensus among agents and the optimality gap.

\begin{lemma}\label{lem-weightedavg}
The weighted average sequence $\{\hat x_k\}$ satisfies, for all $k\ge0$,
\begin{align*}%\label{eq-hatx}
	\hat x_{k+1}=\hat x_{k} -  \sum_{i=1}^n [\phi_{k+1}]_i\a_i y_k^i+\sum_{i=1}^n [\phi_{k+1}]_i\b_i(x_k^i-x_{k-1}^i).
\end{align*}
\end{lemma}

\begin{proof}
    See Appendix~\ref{Appendix-lem-weightedavg}.
\end{proof}

We proceed to define the constants to be used in the subsequent analysis. Define $\mathsf{D}(\bbG_k)$ as the diameter and $\mathsf{K}(\bbG_k)$ as the maximal edge-utility of the graph $\bbG_k$.
% (see Appendix~\ref{Appendix:basicre}). 
The constants $\chi_k\ge 1$, $\varphi_k\ge 1$, $\gamma_k\in(0,1]$, $\psi_k>0$, $\tau_k\in(0,1)$, $c_k\in(0,1)$, $\nu_k>0$ and $\zeta_k>0$, are defined as follows:
\begin{align}\label{eq-const-all}
&\chi_k\!=\!\sqrt{\tfrac{1}{\min(\pi_{k})}},~\varphi_k\!=\!\sqrt{\!\tfrac{1}{\min(\phi_{k})}},~\g_k\!=\!\!\sqrt{\max_{i\in[n]} ([\phi_{k+1}]_i[\pi_k]_i)},
 \nonumber\\
&\psi_k\!=\!n(\chi_{k+1}^2\!\!-\!1),~ \tau_k \!=\!\sqrt{\!1 \!-\!   \tfrac{\min^2(\pi_k)\,b^2}
{\max^2(\pi_k\!) \max(\pi_{k+1}\!) \mathsf{D}(\bbG_k\!)\mathsf{K}(\bbG_k\!)}},\nonumber\\
& c_k=\sqrt{1- \tfrac{\min(\phi_{k+1})\, a^2} {\max^2(\phi_k)\,\mathsf{D}(\bbG_k)\mathsf{K}(\bbG_k)}},~\nu_k=\tfrac{6L^2\chi_{k+1}^2\tau^2}{1-\tau^2}+3\psi_kL^2\!, \nonumber\\
&\zeta_k = \left(c\varphi_{k+1}+\varphi_k \right)^2\nu_k, \!\!\!\!
% &\zeta_k = \tfrac{6L^2\left(c\varphi_{k+1}+\varphi_k \right)^2\chi_{k+1}^2\tau^2}{1-\tau^2}+3L^2\psi_k\left(c\varphi_{k+1}+\varphi_k \right)^2,
\end{align} 
where $c\in(0,1)$ and $\tau\in(0,1)$ are upper bounds for $c_k$ and $\tau_k$, respectively. Furthermore, we denote $\eta \in (0,1)$ as the lower bound for $\phi_{k+1}^\T\pi_k$. Specifically, we have
\begin{align}\label{eq-const-max}
\max_{k\ge0}c_k\le c,~ \max_{k\ge0}\tau_k\le \tau, ~ \min_{k\ge0}\phi_{k+1}^\T\pi_k\ge \eta.
\end{align}

We conclude this section by providing the bound for the conditional expectation of $\left\| \sum_{i=1}^n y_k^i\right\|^2$ given $\F_k$, which will soon be employed to analyze the behavior of the $y$-iterates.
\begin{lemma}\label{lem-ysum}
Let  Assumption~\ref{asm-functions}, Assumption~\ref{asm-strconv}, Assumption~\ref{asm-SFO}, and Assumption~\ref{asm-bmatrices} hold, then for all $k\!\ge\!0$ we have
\begin{align*}
\E\!\!\left[\!\left\| \sum_{i=1}^n y_k^i\right\|^2\!\!\bigg|\F_k\!\right] \!\!\le\! 2nL^2\varphi_k^2 \!\left(\|\hat x_k \!-\! x^*\|^2 \!\!+\! \|\bx_k \!-\!\hat{\bx}_k\|_{\phi_k}^2\!\right)\!+\!2n\sigma^2\!.
\end{align*}
\end{lemma}
\begin{proof}
See Appendix~\ref{Appendix-Lemmalem-ysum}.
\end{proof}

\subsection{Main Results}\label{sec:conv}
In the following analysis, we examine the behavior of the local models $x_k^i$ and local directions $y_k^i$ for $i\in [n]$ and $k\ge 0$, generated by the method in \eqref{eq-met} and derive bounds for four critical quantities in expectation, namely, (i) the optimality gap $\E[\|\hat{x}_k -x^*\|^2]$, (ii) the consensus error $\E[\|\bx_k-\hat{\bx}_k\|_{\phi_k}^2]$, (iii) the state difference $\E[\left\|\bx_{k}-\bx_{k-1}\right]\|^2]$ and (iv) the gradient tracking error $\E[S^2(\by_{k},\pi_k)]$, with $S(\by_{k},\pi_k)$ defined as:
\begin{align}
S(\by_{k},\pi_k) =\sqrt{\sum_{i=1}^n[\pi_k]_i\Bigg\|\frac{y_k^i}{[\pi_k]_i} -\sum_{j=1}^n y_{k}^j\Bigg\|^2}, \label{eq-x-S-quants}
\end{align}
where $\bx_k$, $\by_{k}$, and $\hat \bx_k$ are defined as in the previous sections and $x^*$ is the optimal global solution of problem~\eqref{eq-problem}. These analyses establish a recursive relation for each quantity. Collectively, they form a linear system of inequalities for a vector composed of all four error quantities, governing the convergence behavior of the proposed algorithm.

% We then establish a linear system of inequalities for a vector composed of these four quantities and show that the vector converges linearly to a neighborhood around the optimal solution for sufficiently small constant stepsizes and momentum parameters. 

% The bounds on  $\E\left[\|\hat{x}_k -x^*\|^2|\F_k\right]$, $\E\left[\|\bx_k-\hat{\bx}_k\|_{\phi_k}^2|\F_k\right]$, $\|\bx_{k+1}-\bx_k\|^2$ and $\E\left[S^2(\by_{k},\pi_k)|\F_k\right]$ are given respectively as follows:

Let $\bar{\a} = \max(\balpha)$ and $\bar{\b} = \max(\bbeta)$ represent the largest stepsize and momentum parameter, respectively. We start with establishing a bound on the expected squared distance between $\hat x_k$, the weighted average of the local models, and the global model $x^*$.

\begin{proposition} \label{prop-waverx}
% Let Assumption~\ref{asm-functions}, Assumption~\ref{asm-strconv}, Assumption~\ref{asm-graphs}, Assumption~\ref{asm-SFO} and Assumption~\ref{asm-amatrices} hold. 
Let Assumptions~\ref{asm-functions}--\ref{asm-bmatrices} hold. 
% Let the stepsize in method~\eqref{eq-met} be such that 
For $0\!<\!\bar{\a}\!<\!\tfrac{2}{n\eta(L+\mu)}$, we have for all $k\ge0$, 
\begin{align*}
\E[\|\hat x_{k+1}-x^*\|^2|&\F_k]
\le \left(1- \tfrac{\bar{\a} n\mu \eta}{2}\right)\|\hat x_k - x^*\|^2\\
+&\tfrac{3\bar{\a} L^2\varphi_k^2}{\mu}\|\bx_k-\hat{\bx}_k\|_{\phi_k}^2+\tfrac{6\bar{\b}^2}{\bar{\a} n\mu \eta}\|\bx_k-\bx_{k-1}\|^2\\
+&\tfrac{6\bar{\a}}{ n\mu \eta}\E\left[S^2(\by_{k},\pi_k)|\F_k\right]+\tfrac{3\bar{\a}^2 n\sigma^2}{2}.
\end{align*}
\end{proposition}
\begin{proof}
See Appendix~\ref{Appendix-Propositionprop-waverx}.
\end{proof}

We proceed to analyze the $\phi_{k}$-weighted dispersion of the iterates $x_k^i$ for $i \in [n]$, relative to their weighted average $\hat x_k$, which quantifies consensus among the agents.

\begin{proposition}\label{prop-xcontract}
Let Assumption~\ref{asm-functions}, Assumption~\ref{asm-SFO},
Assumption~\ref{asm-graphs}, and Assumption~\ref{asm-amatrices} hold. We have for all $k\ge0$,
\begin{align*}
&\E\left[\left\|\bx_{k+1}-\hat \bx_{k+1}\right\|_{\phi_{k+1}}^2|\F_k\right]\le\tfrac{1+c^2}{2}\left\|\bx_k - \hat \bx_k\right\|_{\phi_{k}}^2\\
+&\tfrac{2\bar{\a}^2 (1+c^2)\g_k^2}{1-c^2}~ \E\left[\|\by_{k}\|_{\pi_k^{-1}}^2|\F_k\right]+\tfrac{2\bar{\b}^2 (1+c^2)}{1-c^2} \|\bx_k-\bx_{k-1} \|^2,
\end{align*}
where $\g_k$ and $c$ are as defined in \eqref{eq-const-max}.
\end{proposition}
\begin{proof}
See Appendix~\ref{Appendix-Propositionprop-xcontract}.
\end{proof}

The next proposition provides the recursive relation for the state difference between successive iterations of the algorithm.
\begin{proposition}\label{prop-xdiff}
Let Assumption~\ref{asm-graphs} and Assumption~\ref{asm-amatrices} hold.
Then, for all $k\ge0$, we have 
\begin{align*}
\left\|\bx_{k+1}-\bx_k\right\|^2
&\le 3\left(c\varphi_{k+1}+\varphi_k \right)^2 \|\bx_k-\hat{\bx}_k\|_{\phi_k}^2 \\
&+ 3\bar{\a}^2 \|\by_{k}\|_{\pi_k^{-1}}^2+3\bar{\b}^2\left\|\bx_{k}-\bx_{k-1}\right\|^2.
\end{align*}
\end{proposition}
\begin{proof}
See Appendix~\ref{Appendix-Propositionprop-xdiff}.
\end{proof}

The last quantity used to analyze the performance of the \DSGTmTV~algorithm involves measuring the $\pi_k$-weighted dispersion of the scaled directions $\frac{y_k^i}{[\pi_k]_i}$ for $i\in [n]$, relative to the global stochastic direction $\sum_{j=1}^n y_{k}^j$. 
\begin{proposition}\label{prop-ycontract}
% Let Assumption~\ref{asm-functions}, Assumption~\ref{asm-graphs}, Assumption~\ref{asm-SFO} and Assumption~\ref{asm-bmatrices} hold.
Let Assumptions~\ref{asm-functions}--\ref{asm-bmatrices} hold.
Then, for the quantity $S(\by_{k},\pi_k)$ defined in~\eqref{eq-x-S-quants}, we have for all $k\ge0$,
\begin{align*}
\E\left[S^2(\by_{k+1},\pi_{k+1})|\F_k\right]
\le \tfrac{1+\tau^{2}}{2}\E\left[S^2(\by_{k},\pi_k)|\F_k\right]& \nonumber\\
+\zeta_k\|\bx_k-\hat{\bx}_k\|_{\phi_k}^2+\nu_k \bar{\a}^2\E\left[\|\by_{k}\|_{\pi_k^{-1}}^2|\F_k\right]&\nonumber\\
+\nu_k \bar{\b}^2\|\bx_k-\bx_{k-1} \|^2+2(2+\bar{\a} L)n\psi_k\sigma^{2}&,
\end{align*}
where the scalars $\zeta_k$, $\nu_k$ and $\psi_k$ are given in \eqref{eq-const-all} while $\tau\in(0,1)$ is defined as in \eqref{eq-const-max}.
\end{proposition}
\begin{proof}
See Appendix~\ref{Appendix-Propositionprop-ycontract}.
\end{proof}

The bounds derived in the preceding propositions depend on $\E\left[\|\by_{k}\|_{\pi_k^{-1}}^2|\F_k\right]$. To establish a composite relation, we introduce the following lemma, which bounds this quantity in terms of the four considered error quantities.
\begin{lemma}\label{lem-ypik}
Let Assumption~\ref{asm-graphs} and Assumption~\ref{asm-bmatrices} hold. Then
\begin{align*}
\E\left[\|\by_{k}\|_{\pi_k^{-1}}^2|\F_k\right] \le&~ 2nL^2\varphi_k^2 \left(\|\hat x_k - x^*\|^2+ \|\bx_k -\hat{\bx}_k\|_{\phi_k}^2\right)\\
&+\E[S^2(\by_{k},\pi_k)|\F_k]+2n\sigma^2.
\end{align*}
\end{lemma}
\begin{proof}
See Appendix~\ref{Appendix-Lemmalem-ypik}.
\end{proof}

\subsection{Composite Relation} \label{subsec-comrel}
% Building on the preceding results, w
We now present a composite relation for the four error quantities. We define the vector $V_k$ as follows:
\begin{align}\label{eq-vk}
V_k=[&\E[\|\hat{x}_k -x^*\|^2],~\E[\|\bx_k-\hat{\bx}_k\|_{\phi_k}^2], \nonumber\\
&\E[\left\|\bx_{k}-\bx_{k-1}\right]\|^2],~\E[S^2(\by_{k},\pi_k)]]^\T.
\end{align}
% \begin{align}\label{eq-vk}
% \!\!V_k=\begin{bmatrix}
%            \E[\|\hat{x}_k -x^*\|^2] \\
%            \E[\|\bx_k-\hat{\bx}_k\|_{\phi_k}^2] \\
%         ~\E[\left\|\bx_{k}-\bx_{k-1}\right]\|^2]~\\
%            \E[S^2(\by_{k},\pi_k)]
%          \end{bmatrix}.
% \end{align}
By substituting the bound on $\E[\|\by_{k}\|_{\pi_k^{-1}}^2|\F_k]$ from Lemma~\ref{lem-ypik} in  Proposition~\ref{prop-xcontract}, Proposition~\ref{prop-xdiff}, and Proposition~\ref{prop-ycontract}, and taking the full expectation of both sides, it can be verified that for $0\!<\!\bar{\a}\!<\!\tfrac{2}{n\eta(L+\mu)}$ then $V_k$ follows the dynamical system below:
\begin{equation*}
V_{k+1}\le
M_k(\bar{\a},\bar{\b})V_k+\bb_k(\bar{\a}),
\end{equation*}
where $\bb_k(\bar{\a})=\left[b_k^{1}\bar{\a}^2,b_k^{2}\bar{\a}^2,b_k^{3}\bar{\a}^2,b_k^{4}+b_k^{5}\bar{\a}+b_k^{6}\bar{\a}^2\right]^\T$, and $M_k(\bar{\a},\bar{\b})$ is given by
\begin{align*}
\left[\begin{array}{cccc}
1-m_k^{1}\bar{\a} & m_k^{2}\bar{\a} & m_k^{3}\frac{\bar{\b}^2}{\bar{\a}} & m_k^{3}\bar{\a} \cr
% \hbox{}\cr
m_k^{4}\bar{\a}^2 &\tfrac{1+c^2}{2} + m_k^{4}\bar{\a}^2 & m_k^{5}\bar{\b}^2 & m_k^{6}\bar{\a}^2\cr
% \hbox{}\cr
m_k^{7}\bar{\a}^2 & m_k^{8}+m_k^{7}\bar{\a}^2 & 3\bar{\b}^2& 3\bar{\a}^2 \cr
% \hbox{}\cr
m_k^{9}\bar{\a}^2 & m_k^{10}+m_k^{9}\bar{\a}^2 & m_k^{11}\bar{\b}^2 & \tfrac{1+\tau^2}{2}+m_k^{11}\bar{\a}^2 \end{array}\right].
\end{align*}
% \begin{align*}
% \left[\begin{array}{cccc}
% 1-m_k^{1}\bar{\a} & m_k^{2}\bar{\a} & m_k^{3}\frac{\bar{\b}^2}{\bar{\a}} & m_k^{3}\bar{\a} \cr
% m_k^{4}\bar{\a}^2 &\tfrac{1+c^2}{2} + m_k^{4}\bar{\a}^2 & m_k^{5}\bar{\b}^2 & m_k^{6}\bar{\a}^2\cr
% m_k^{7}\bar{\a}^2 & m_k^{8}+m_k^{7}\bar{\a}^2 & 3\bar{\b}^2& 3\bar{\a}^2 \cr
% m_k^{9}\bar{\a}^2 & m_k^{10}+m_k^{9}\bar{\a}^2 & m_k^{11}\bar{\b}^2 & \tfrac{1+\tau^2}{2}+m_k^{11}\bar{\a}^2 \end{array}\right]
% \end{align*}
Here, the constants are as follows:
\begin{align*}
&m^{1}_{k}=\tfrac{n\mu \eta}{2},~m^{2}_{k}=\tfrac{3 L^2\varphi_k^2}{\mu},~m^{3}_{k}=\tfrac{6}{ n \mu \eta},\\
&m^{4}_{k}=\tfrac{4nL^2\varphi_k^2(1+c^2)\g_k^2}{1-c^2},~m^{5}_{k}=\tfrac{2(1+c^2)}{1-c^2},~m^{6}_{k}=m^{5}_{k}\g_k^2,\\
&m^{7}_{k}=6nL^2\varphi_k^2,~m^{8}_k=3\left(c\varphi_{k+1}+\varphi_k \right)^2,~m^{9}_{k}=2nL^2\varphi_k^2\nu_k,\\
&m^{10}_{k}=\zeta_k,~m^{11}_{k}=\nu_k,~b_k^{1}=\tfrac{3 n \sigma^2}{2},~b_k^{2}=\tfrac{4n(1+c^2)\g_k^2\sigma^2}{1-c^2},\\
&b_k^{3}=6n\sigma^2,~b_k^{4}=4n\psi_k\sigma^2,~b_k^{5}=2Ln\psi_k\sigma^2,~b_k^{6}=2n\nu_k\sigma^2.
\end{align*}

We now determine an upper bound matrix $M(\bar{\a},\bar{\b})$ for $M_k(\bar{\a},\bar{\b})$. For the constants defined in \eqref{eq-const-all}, we note that $\gamma_k \le 1$ and $\eta <\phi_{k+1}^\T\pi_k \le 1$.
% , for all $k$. 
Let $\psi>0$, $\chi>1$, and $\varphi>1$, be upper bounds for $\psi_k$, $\chi_k$, and $\varphi_k$, respectively, i.e.,
\begin{align}\label{eq-const-max2}
\max_{k\ge0}\psi_k\le \psi,~\max_{k\ge0} \chi_k \le \chi,~ \max_{k\ge0} \varphi_k \le \varphi.
\end{align}
We will also define the following constants:
% for our subsequent analysis:
\begin{align}\label{eq-const-nuvarsigma}
\nu ~&=\tfrac{6L^2\chi^2\tau^2}{1-\tau^2}+3\psi L^2, \\
\varsigma^2 &= \tfrac{L^2(n\eta\mu^2+12L^2\varphi^2)}{\eta\mu^2} + \tfrac{8\nu(n\eta^2\mu^2+48L^2\varphi^2)}{n\eta\mu^2(1-\tau^2)}. \nonumber
% \varsigma^2 = \tfrac{(8\nu+L^2n\eta)(n\eta\mu^2+48L^2\varphi^2)}{n\eta^2\mu^2(1-\tau^2)}. \!
% \varsigma^2 = \tfrac{L^2(\eta+48\chi^2)(n\mu^2\eta+48L^2\varphi^2)}{\mu^2\eta^2(1-\tau^2)^2}
\end{align}

Using these upper bounds, we have $M_k(\bar{\a},\bar{\b})\le M(\bar{\a},\bar{\b})$ and $\bb_k(\bar{\a})\le \bb(\bar{\a})$, for all $k \ge 0$. It is straightforward to obtain the following proposition:

\begin{proposition}\label{prop-composite-relation}
Let Assumptions~\ref{asm-functions}--\ref{asm-bmatrices} hold. For $0\!<\!\bar{\a}\!<\!\tfrac{2}{n\eta(L+\mu)}$, $V_k$ evolves according to the following system:
\begin{equation}\label{eq-vkrel}
V_{k+1}\le M(\bar{\a},\bar{\b})V_k+\bb(\bar{\a}), \qquad \hbox{for all $k\ge0$},
\end{equation}
where $V_k$ is as defined in \eqref{eq-vk}, and $M(\bar{\a},\bar{\b})$ and $\bb(\bar{\a})$ are given respectively by
\begin{align*}
&\!\!M(\bar{\a},\bar{\b})\!=\!\!\left[\!\begin{array}{cccc}
\!1\!-\!m_{1}\bar{\a}\!\! & m_{2}\bar{\a} & \!\!m_{3}\frac{\bar{\b}^2}{\bar{\a}}\!\! & m_{3}\bar{\a} \cr
% \hbox{}\cr
\!\!m_{4}\bar{\a}^2\!\! &\!\tfrac{1+c^2}{2} \!+\! m_{4}\bar{\a}^2 & \!\!m_{5}\bar{\b}^2\!\! & m_{5}\bar{\a}^2\cr
% \hbox{}\cr
\!\!m_{7}\bar{\a}^2\!\! & \!\!m_{8}\!+\!m_{7}\bar{\a}^2 & \!\!3\bar{\b}^2\!\!& 3\bar{\a}^2 \cr
% \hbox{}\cr
\!\!m_{9}\bar{\a}^2\!\! & \!\!m_{10}\!+\!m_{9}\bar{\a}^2 & \!\!m_{11}\bar{\b}^2\! & \!\tfrac{1+\tau^2}{2}\!+\!m_{11}\bar{\a}^2\!\!\! \end{array}\right],\!\!\!\! \nonumber\\
&%\!\!\!\text{and} 
\bb(\bar{\a})=\left[b_{1}\bar{\a}^2,b_{2}\bar{\a}^2,b_{3}\bar{\a}^2,b_{4}+b_{5}\bar{\a}+b_{6}\bar{\a}^2\right]^\T,%\label{eq-vecg}
\end{align*}
where the constants in these expressions are given by
% $M(\bar{\a},\bar{\b})$ and $\bb(\bar{\a})$ are 
\begin{align*}
&m_{1}=\tfrac{n\mu\eta}{2},~m_{2}=\tfrac{3L^2\varphi^2}{\mu},~m_{3}=\tfrac{6}{ n \mu \eta},~m_{4}=\tfrac{4nL^2\varphi^2(1+c^2)}{1-c^2},\\
&m_{5}=\tfrac{2(1+c^2)}{1-c^2}\!,~m_{7}=6nL^2\varphi^2\!,~m_{8}=12\varphi^2\!,~m_{9}=2nL^2\!\varphi^2\!\nu,\\
&m_{10}=4\varphi^2\nu,~m_{11}=\nu,~b_{1}=\tfrac{3 n\sigma^2}{2},~b_{2}=\tfrac{4n(1+c^2)\sigma^2}{1-c^2},\\
&b_{3}=6n\sigma^2,~b_{4}=4n\psi\sigma^2,~b_{5}=2Ln\psi\sigma^2,~b_{6}=2n\nu\sigma^2,
\end{align*}
with
$\nu$ as defined in \eqref{eq-const-nuvarsigma}.
\end{proposition}

\subsection{Convergence Theorem} \label{subsec-conv-theo}
In this section, we establish the main result regarding the convergence properties of the \DSGTmTV~algorithm. Specifically, with sufficiently small constant stepsizes and momentum parameters, the algorithm converges linearly to a neighborhood surrounding the optimal global model. Additionally, in the deterministic setting where exact gradients are available without stochastic noise, the \DSGTmTV~algorithm achieves linear convergence directly to the global model.

\def\limsup{\mathop{\rm limsup}}

\begin{theorem}\label{theo-main}
% Let Assumption~\ref{asm-functions}, Assumption~\ref{asm-strconv}, Assumption~\ref{asm-graphs}, Assumption~\ref{asm-SFO}, Assumption~\ref{asm-amatrices} and Assumption~\ref{asm-bmatrices} hold.
Let Assumptions~\ref{asm-functions}--\ref{asm-bmatrices} hold. Consider the iterates produced by Algorithm 1, the notations in \eqref{eq-const-all}, \eqref{eq-const-max}, \eqref{eq-const-max2} and \eqref{eq-const-nuvarsigma}. Assume we choose the largest stepsize $\bar{\a}>0$ such that 
\begin{align}\label{eq:alpha-range}
\bar{\a} < \min \left\{\tfrac{2}{n\eta(L\!+\!\mu)},\tfrac{1\!-\!c^2}{2\varphi\varsigma\sqrt{2(1+c^2)}}\right\},
\end{align}
and the largest momentum parameter $\bar{\b}\ge 0$ satisfies
\begin{align}\label{eq:beta-range}
\!\!\bar{\b} \!<\! \min \!\left\{\!\bar{\a}\sqrt{\tfrac{n\eta L^2(1 \!-\! \tau^2) \!+\! 32\nu}{48(1 \!-\! \tau^2)}},\sqrt{\!\tfrac{(1-c^2)^2}{96(1\!+\!c^2)\varphi^2}\!-\!\tfrac{\varsigma^2\bar{\a}^2}{12}}\right\}.\!\!\!
\end{align}
Then, $\rho_{M}<1$ where $\rho_{M}$ is the spectral radius of $M(\bar{\a},\bar{\b})$, the vector $(\mbi-M(\bar{\a},\bar{\b}))^{-1}\bb(\bar{\a})$ has non-negative components, and we have that
\begin{align*}
&\limsup_{k\to\infty} V_k \le (\mbi-M(\bar{\a},\bar{\b}))^{-1}\bb(\bar{\a}),\\
&\limsup_{k\to\infty} \E[\|\hat{x}_k -x^*\|^2] \le [(\mbi-M(\bar{\a},\bar{\b}))^{-1}\bb(\bar{\a})]_1,\\
&\limsup_{k\to\infty} \E[\|\bx_k-\hat{\bx}_k\|_{\phi_k}^2] \le [(\mbi-M(\bar{\a},\bar{\b}))^{-1}\bb(\bar{\a})]_2,
\end{align*}
with a linear convergence rate of the order of $\mathcal{O}\left(\rho_M^k\right)$, where the matrix $M(\bar{\a},\bar{\b})$ and the vector $\bb(\bar{\a})$ are defined in Proposition~\ref{prop-composite-relation}. Here, $[(\mbi-M(\bar{\a},\bar{\b}))^{-1}\bb(\bar{\a})]_j$ is the $j$-th element of the vector $(\mbi-M(\bar{\a},\bar{\b}))^{-1}\bb(\bar{\a})$.
\end{theorem}

\begin{proof}
Recall from \eqref{eq-vkrel} that we have
\[V_{k+1}\le M(\bar{\a},\bar{\b})V_k+\bb(\bar{\a}), \qquad \hbox{for all $k\ge0$}.\]
Thus,
\[V_{k}\le [M(\bar{\a},\bar{\b})]^kV_0+\sum_{t=0}^{k-1}[M(\bar{\a},\bar{\b})]^{t}\bb(\bar{\a}), ~~ \hbox{for all $k\ge 1$}.\]
The matrix $[M(\bar{\a},\bar{\b})]^k$ converges to $\zero$, as $k \to \infty$, linearly at rate $\mathcal{O}\left(\rho_M^k\right)$ if the spectral radius $\rho_{M}$ of $M(\bar{\a},\bar{\b})$ satisfies $\rho_{M}<1$. Furthermore, $\mbi-M(\bar{\a},\bar{\b})$ is nonsingular \cite{horn1990matrix}, and
\[\sum_{k=0}^{\infty}[M(\bar{\a},\bar{\b})]^{k}=(\mbi-M(\bar{\a},\bar{\b}))^{-1}.\]
The goal is to find the range of the stepsizes and momentum parameters such that $\rho_{M}<1$. In light of Lemma~\ref{lem-SpectralRadius}, it suffices to solve for the range of $\bar{\a}$ and $\bar{\b}$ such that $M(\bar{\a},\bar{\b})\bdelta < \bdelta$ holds for some positive vector $\bdelta=[\delta_1,\delta_2,\delta_3,\delta_4]^\T$. Expanding the $M(\bar{\a},\bar{\b})\bdelta < \bdelta$ inequality element-wise, we obtain
\begin{align*}
(1-m_1\bar{\a})\delta_1+m_2\bar{\a}\delta_2+\tfrac{m_3\bar{\b}^2}{\bar{\a}}\delta_3+m_3\bar{\a}\delta_4&\!<\!\delta_1,\\
m_{4}\bar{\a}^2\delta_1+\left(\tfrac{1+c^2}{2}+ m_{4}\bar{\a}^2\right)\delta_2+m_{5}\bar{\b}^2\delta_3+m_{5}\bar{\a}^2\delta_4&\!<\!\delta_2,\\
m_{7}\bar{\a}^2\delta_1+(m_{8}+m_{7}\bar{\a}^2)\delta_2+3\bar{\b}^2\delta_3+3\bar{\a}^2\delta_4&\!<\!\delta_3,\\
\!m_{9}\bar{\a}^2\!\delta_1\!\!+\!(m_{10}\!+\!m_{9}\bar{\a}^2)\delta_2\!\!+\!m_{11}\bar{\b}^2\!\delta_3\!\!+\!\!\left(\!\tfrac{1+\tau^2}{2}\!\!+\!m_{11}\bar{\a}^2\!\right)\!\delta_4&\!<\!\delta_4,
\end{align*}
which are equivalent to 
\begin{align}
&\!\!\!m_3\delta_3\bar{\b}^2<(m_1\delta_1-m_2\delta_2-m_3\delta_4)\bar{\a}^2,\label{eq-alpha1}\\
&\!\!\!(m_{4}\delta_1+m_{4}\delta_2+m_{5}\delta_4)\bar{\a}^2+m_5\delta_3\bar{\b}^2<\tfrac{1-c^2}{2}\delta_2,\label{eq-alpha2}\\
&\!\!\!(m_{7}\delta_1+m_{7}\delta_2+3\delta_4)\bar{\a}^2+3\delta_3\bar{\b}^2<\delta_3-m_{8}\delta_2,\label{eq-alpha3}\\
&\!\!\!(m_{9}\delta_1\!+\!m_{9}\delta_2\!+\!m_{11}\delta_4)\bar{\a}^2\!+\!m_{11}\delta_3\bar{\b}^2\!<\!\tfrac{1-\tau^2}{2}\delta_4\!-\!m_{10}\delta_2.\!\!\!\label{eq-alpha4}
\end{align}
For the inequality in \eqref{eq-alpha1} to hold, the following is required:
\begin{align}\label{eq-delta1}
\delta_1>\frac{m_2\delta_2+m_3\delta_4}{m_1}.
\end{align}
Additionally, to ensure that the RHSs of both \eqref{eq-alpha3} and \eqref{eq-alpha4} are positive, we need the following conditions to be fulfilled:
\begin{align}\label{eq-delta23}
\delta_3>m_{8}\delta_2 \text{and} 0<\delta_2<\frac{1-\tau^2}{2m_{10}}\delta_4.
\end{align}
It is straightforward to verify that the following values of $\delta_1$, $\delta_2$, $\delta_3$, and $\delta_4$ satisfy both \eqref{eq-delta1} and \eqref{eq-delta23}:
\begin{align*}
\delta_1\!=\!\frac{2}{m_1}\!\!\left(\!m_2+\frac{4m_3m_{10}}{1-\tau^2}\!\right)\!, \delta_2=1,\delta_3=2m_{8}, \delta_4=\frac{4m_{10}}{1-\tau^2}.
\end{align*}
Plugging in these values of $\bdelta$ into \eqref{eq-alpha1}--\eqref{eq-alpha4}, we have
$M(\bar{\a},\bar{\b})\bdelta < \bdelta$ when the following inequalities hold:
\begin{align*}
&\!\!\!\left(\!m_2\!+\!\tfrac{4m_3m_{10}}{1-\tau^2}\!\right)\bar{\a}^2>2m_3m_8\bar{\b}^2,\!\!\!\\
&\!\!\!\left(\!m_4\!+\!\tfrac{2m_4}{m_1}\!\!\left(\!m_2\!+\!\tfrac{4m_3m_{10}}{1-\tau^2}\!\right)\!\!+\!\tfrac{4m_5m_{10}}{1-\tau^2}\!\right)\bar{\a}^2+2m_5m_8\bar{\b}^2\!<\!\tfrac{1-c^2}{2},\!\!\!\nonumber\\
&\!\!\!\left(\!\tfrac{m_7}{m_8}\!+\!\tfrac{2m_7}{m_1m_8}\!\left(\!m_2\!+\!\tfrac{4m_3m_{10}}{1-\tau^2}\!\right)\!+\!\tfrac{12m_{10}}{m_8(1-\tau^2)}\!\right)\bar{\a}^2+6\bar{\b}^2\!<1,\!\!\!\nonumber\\
&\!\!\!\left(\!\tfrac{m_9}{m_{10}}\!+\!\tfrac{2m_9}{m_1m_{10}}\!\left(\!m_2\!+\!\tfrac{4m_3m_{10}}{1-\tau^2}\!\right)\!+\!\tfrac{4m_{11}}{1-\tau^2}\!\right)\bar{\a}^2+\tfrac{2m_8m_{11}}{m_{10}}\bar{\b}^2\!<1.\!\!\nonumber
\end{align*}
We can further verify that 
\begin{align*}
&m_4+\tfrac{2m_4}{m_1}\left(m_2+\tfrac{4m_3m_{10}}{1-\tau^2}\right)\!+\tfrac{4m_5m_{10}}{1-\tau^2}<\tfrac{4(1+c^2)\varphi^2\varsigma^2}{1-c^2},\\
&\tfrac{m_7}{m_8}+\tfrac{2m_7}{m_1m_8}\left(m_2+\tfrac{4m_3m_{10}}{1-\tau^2}\right)+\tfrac{12m_{10}}{m_8(1-\tau^2)}< \tfrac{\varsigma^2}{2},\\
&\tfrac{m_9}{m_{10}}+\tfrac{2m_9}{m_1m_{10}}\left(m_2+\tfrac{4m_3m_{10}}{1-\tau^2}\right)+\tfrac{4m_{11}}{1-\tau^2}<\tfrac{\varsigma^2}{2},
\end{align*}
where $\varsigma^2$ is as defined in \eqref{eq-const-nuvarsigma}.
% where $\varsigma^2 = \tfrac{L^2(\eta+48\chi^2)(n\mu^2\eta+48L^2\varphi^2)}{\mu^2\eta^2(1-\tau^2)^2}$ as defined in \eqref{eq-const-nuvarsigma}.

Using the above bounds, we obtain the upper bounds for $\bar{\a}$ and $\bar{\b}$ as shown in \eqref{eq:alpha-range} and \eqref{eq:beta-range}, respectively, thereby establishing the proof.
% \begin{align}
% &\!\!\b^2<\left(\frac{n\eta L^2}{48}+\frac{4L^2\chi^2\tau^2}{(1-\tau^2)^2}+\frac{2l^2n(\chi^2-1)}{1-\tau^2}\right)\a\!\!\!\label{eq-solve-beta-range1}\\
% &\!\!\b^2\!<\!\frac{(1-c^2)^2}{192\varphi^2}-\frac{L^2(n\mu^2+48L^2\varphi^2)(\eta+48\chi^2)}{24\mu^2\eta(1-\tau^2)^2}\a^2\!\!\!\label{eq-solve-beta-range2}\\
% &\!\!\b^2\!<\!\frac{1}{6}-\frac{L^2\left(\eta+48nL^2\chi^2\right)\left(n^2\mu^2\eta+72L^2\varphi^2\right)}{12\mu^2\eta(1-\tau^2)^2}\a^2\label{eq-solve-beta-range3}\\
% &\!\!\b^2\!<\!\frac{1}{6}-\frac{n\chi^2L^2\left(n\mu^2+24L^2\varphi^2\right)\left(\eta+96\chi^2\right)}{2\mu^2\eta(1-\tau^2)^4(\chi^2-1)(2\tau^2+n)}\a^2\label{eq-solve-beta-range4}
% \end{align}
\end{proof}

\begin{remark}
The error bounds in Theorem~\ref{theo-main} go to zero as the largest value of the stepsizes $\bar{\a}$ gets smaller and the variance $\sigma$ on the gradient noise decreases. 
% Consequently, the \DSGTmTV~algorithm exhibits linear convergence towards the global minimizer.
\end{remark}

\begin{remark} \label{rem-CgraphsAB}
When Assumption~\ref{asm-graphs} is relaxed and the graph sequence $\{\mathbb{G}_k\}$ is strongly connected over a period $C>1$, the more general results of Lemma~\ref{lem-amatrices} and Lemma~\ref{lem-bmatrices} indicate the existence of stochastic vector sequences $\{\phi_k\}$ and $\{\pi_k\}$, such that for all $k \ge 0$ and $i\in[n]$, it holds that $[\phi_k]_i\ge\frac{a^{nC}}{n}$ and $[\pi_k]_i\ge\frac{b^{nC}}{n}$. Furthermore,
\begin{align*}
&\phi^{\T}_{k+C}\left(A_{k+C-1}\ldots A_{k+1}A_{k}\right)=\phi^{\T}_k, \\
\!\!\!\text{and} &\pi_{k+C}=\left(B_{k+C-1}\ldots B_{k+1}B_{k}\right)\pi_k,  \text{for all} k \ge 0.
\end{align*}
In this case, the convergence analysis follows similarly to the provided analysis for strongly connected graphs, where the contraction due to stochastic matrices occurs after time $k \!=\! C$.
\end{remark}

\begin{remark}
% In Theorem~\ref{theo-main}, we have established a linear convergence rate of the \DSGTmTV~algorithm to a neighborhood of the global minimizer when the largest stepsize, $\bar{\a}$, and the largest momentum parameter, $\bar{\b}$, follow \eqref{eq:alpha-range} and \eqref{eq:beta-range}. However,
Computing the bounds for the largest stepsize, $\bar{\a}$, and the largest momentum parameter, $\bar{\b}$, using \eqref{eq:alpha-range} and \eqref{eq:beta-range} in Theorem~\ref{theo-main}, may be impractical. This difficulty arises from the intricate dependencies of contraction coefficients, $c$ and $\tau$, and other global constants on the connectivity structure of the communication graph, which may be unknown -- especially in the decentralized settings considered in this work. Additionally, it is worth noting that the theoretical bounds in Theorem~\ref{theo-main} tend to be conservative, as has also been noted in prior research \cite{Xin2020Heavy,pshi21}. In practical applications, achieving optimal performance often requires a manual optimization of stepsizes and momentum parameters \cite{YANG2020_hyperparameter}.
\end{remark}

\begin{remark}
In light of the bounds provided in Theorem~\ref{theo-main}, the largest stepsize and largest momentum parameter are independent of their heterogeneity. This eliminates the need to balance the trade-off between the tolerance for stepsize heterogeneity and the achievable largest stepsize, as required in some prior works \cite{Xu2018,LU2018516,nedic2017geometrically}. Consequently, agents have the flexibility to employ different nonnegative values for the stepsizes and momentum parameters, with only one of these stepsizes needing to be positive. 
\end{remark}

% \begin{remark}
% In light of the bounds provided in Theorem~\ref{theo-main}, there exists a trade-off between the tolerance for stepsize heterogeneity, determined by the condition number $\kappa$ of the stepsize matrix $\Da$, and the achievable largest stepsize $\bar{\alpha}$.
% \end{remark}

\section{Numerical results} \label{sec:simulation}
In this section, we evaluate the effectiveness and convergence properties of \DSGTmTV~through empirical studies on two tasks: image classification and natural language processing (NLP). Specifically, we analyze the performance of \DSGTmTV~and compare it with established benchmarks on the MNIST dataset for image classification and a sentiment analysis task for NLP, using the SMS Spam Collection dataset.

\subsection{Simulation Settings}
We consider $n=10$ agents embedded in a communication network that is directed and time-varying \cite{NguyenAutonomousDriving2022,XiaoUAVFL2021}.
Each agent $i$ possesses a local dataset $\mathcal{D}_i$, comprising $m_i$ training samples $\{\xi_{ij} \!=\! (\mathscr{b}_{ij},\ell_{ij})\}_{j=1}^{m_i}$. There is a total of $m=\sum_{i=1}^nm_i$ labeled data points for training and an additional $m^{test}$ labeled data points for testing.

The $j$-th sample at node $i$ is a tuple $\{\mathscr{b}_{ij},\ell_{ij}\} \!\subseteq \!\re^p \!\times \!\{
1,\!-1\}$. To construct an estimate of the coefficients $x\!\!=\!\![x_0,x_{1:}^\T]^\T\!\!\in\! \re^{p+1}$, where $x_{1:}\!\!=\!\![x_1,\ldots,x_p]^\T\!$, the local logistic regression loss function for the $j$-\textit{th} image at node $i$ is given by:
\[\mathcal{L}(x,\xi_{ij}) = \ln \Big[1+\exp\Big\{-(x_{1:}^\T \mathscr{b}_{ij}+x_0)\ell_{ij}\Big\} \Big] \!+\dfrac{\lambda}{2}\|x\|^2\!,\]
which is smooth and strongly convex due to the inclusion of the L2--regularization.

% \begin{table}[t!]
% \begin{tabular}{|l|c|c|c|c|c|c|c|c|}
% \hline
% \textbf{Algorithms}          & GD       & SGD      & CFL--GD   & CFL--SGD  & DFL--GD   & DFL--SGD  & \DGTmTV  & \DSGTmTV~\\ \hline
% \textbf{Exact Gradient}      & \cmark & \xmark & \cmark & \xmark & \cmark & \xmark & \cmark & \xmark \\ \hline
% \textbf{Stochastic Gradient} & \xmark & \cmark & \xmark & \cmark & \xmark & \cmark & \xmark & \cmark \\ \hline
% \textbf{Centralized ML}      & \cmark & \cmark & \xmark & \xmark & \xmark & \xmark & \xmark & \xmark \\ \hline
% \textbf{Centralized FL}      & \xmark & \xmark & \cmark & \cmark & \xmark & \xmark & \xmark & \xmark \\ \hline
% \textbf{Decentralized FL}    & \xmark & \xmark & \xmark & \xmark & \cmark & \cmark & \cmark & \cmark \\ \hline
% \end{tabular}
% \caption{Properties of Algorithms}
% \label{tab:algorithms_properties}
% \end{table}

We compare the performance of the proposed algorithm with state-of-the-art methods across three different settings:

\begin{itemize}
\item \textit{Traditional ML:} Data from multiple sources or agents is aggregated and processed centrally at a single server.
% using a gradient descent (GD) algorithm.
\item \textit{Centralized Federated Learning (CFL) (see Fig.~\ref{fig:FL_centralized}):} Agents locally train models using their own data, and a centralized server aggregates the models from the agents. 
% Each agent employs the GD algorithm.
\item \textit{Decentralized Federated Learning (DFL) (see Fig.~\ref{fig:FL_decentralized}):} There is no centralized server; instead, agents aggregate models through P2P communication.
% We consider two algorithms : (1) the conventional GD algorithm \cite{Tiancheng2022} and (2) the proposed algorithm with gradient tracking and heavy-ball momentum (\DSGTmTV).
\end{itemize}

For each algorithm, we consider two variants: (1) \textit{deterministic} gradient descent (GD) approach, where the entire dataset is used 
% at either the central server or each agent 
to compute exact gradients ($\sigma^2 = 0$), and (2) \textit{stochastic} GD (SGD) approach, which uses small, randomly sampled batches of the dataset to compute the gradients ($\sigma^2 > 0$). Table~\ref{tab:algorithms_properties} provides a comparative overview of the algorithms. 
% based on their gradient computation methods and learning settings.

\begin{table}[t!]
\begin{tabular}{|l|c|c|c|c|c|}
\hline
\!\!\!\textnormal{Algorithms} \!\!\!\!         & \textnormal{\begin{tabular}[c]{@{}c@{}}\!\!\!\!Exact\!\!\!\!\\  \!\!\!\!\! Gradient \!\!\!\!\!\end{tabular}} & \textnormal{\begin{tabular}[c]{@{}c@{}}\!\!\!\!\! Stochastic \!\!\!\!\!\\ \!\!\!\!Gradient\!\!\!\!\end{tabular}} & \textnormal{\begin{tabular}[c]{@{}c@{}}\!\!\!Centralized\!\!\!\\ ML\end{tabular}} & \textnormal{\begin{tabular}[c]{@{}c@{}}\!\!\!Centralized\!\!\!\\ FL\end{tabular}} & \textnormal{\begin{tabular}[c]{@{}c@{}}\!\!\!Decentralized\!\!\!\\ FL\end{tabular}} \\ \hline
\!\!\!\textnormal{GD}   \!\!\!               & \cmark                &                     & \cmark               &                 &                 \\ \hline
\!\!\!\textnormal{SGD}\!\!\!                 &                 & \cmark                    & \cmark               &                 &                 \\ \hline
\!\!\!\textnormal{CFL--GD}\!\!\!              & \cmark                &                     &                & \cmark                &                 \\ \hline
\!\!\!\textnormal{CFL--SGD}\!\!\!             &                 & \cmark                    &                & \cmark                &                 \\ \hline
\!\!\!\textnormal{DFL--GD}\!\!\!              & \cmark                &                     &                &                 & \cmark                \\ \hline
\!\!\!\textnormal{DFL--SGD}\!\!\!             &                 & \cmark                    &                &                 & \cmark                \\ \hline
\!\!\!\textnormal{\DGTmTV}\!\!\!             & \cmark                &                     &                &                 & \cmark                \\ \hline
\!\!\!\textnormal{\DSGTmTV}\!\!\!             &                 & \cmark                    &                &                 & \cmark                \\ \hline
\end{tabular}
\caption{Algorithm Properties and Settings.}
\label{tab:algorithms_properties}
\vspace{-0.6cm}
\end{table}

% In contrast, the centralized stochastic gradient descent (CSGD) approach uses a small batch of data point, randomly sampled from the entire batch, for each iteration. We also compare with the non-stochastic $AB$/Push-Pull (ABPP) \cite{Saadatniaki2020, Angelia2022AB} and $AB$/Push-Pull with heavy-ball momentum (ABPPm)  \cite{Nguyen2023AccAB} methods for decentralized algorithms, where each agent $i$ processes a local batch of $m_i$ labeled data points. On the other hand, in the proposed  \DSGTmTV~algorithm, each agent selects one data point uniformly from its local batch, for each iteration. The \SABTV~algorithm is the variant without heavy-ball momentum \cite{Nguyen2023SABTV}.

\subsection{Application 1: Image Classification}
We perform digit classification on the widely used MNIST dataset, with the goal of  distinguishing between handwritten digits $3$ and $5$.
Each image is represented as a $784$-dimensional vector. We use $2000$ samples for training and $500$ samples for testing purposes, with a regularization parameter $\lambda=0.001$.

\begin{figure}[t!]
\centering
~\includegraphics[width=0.2235\textwidth]{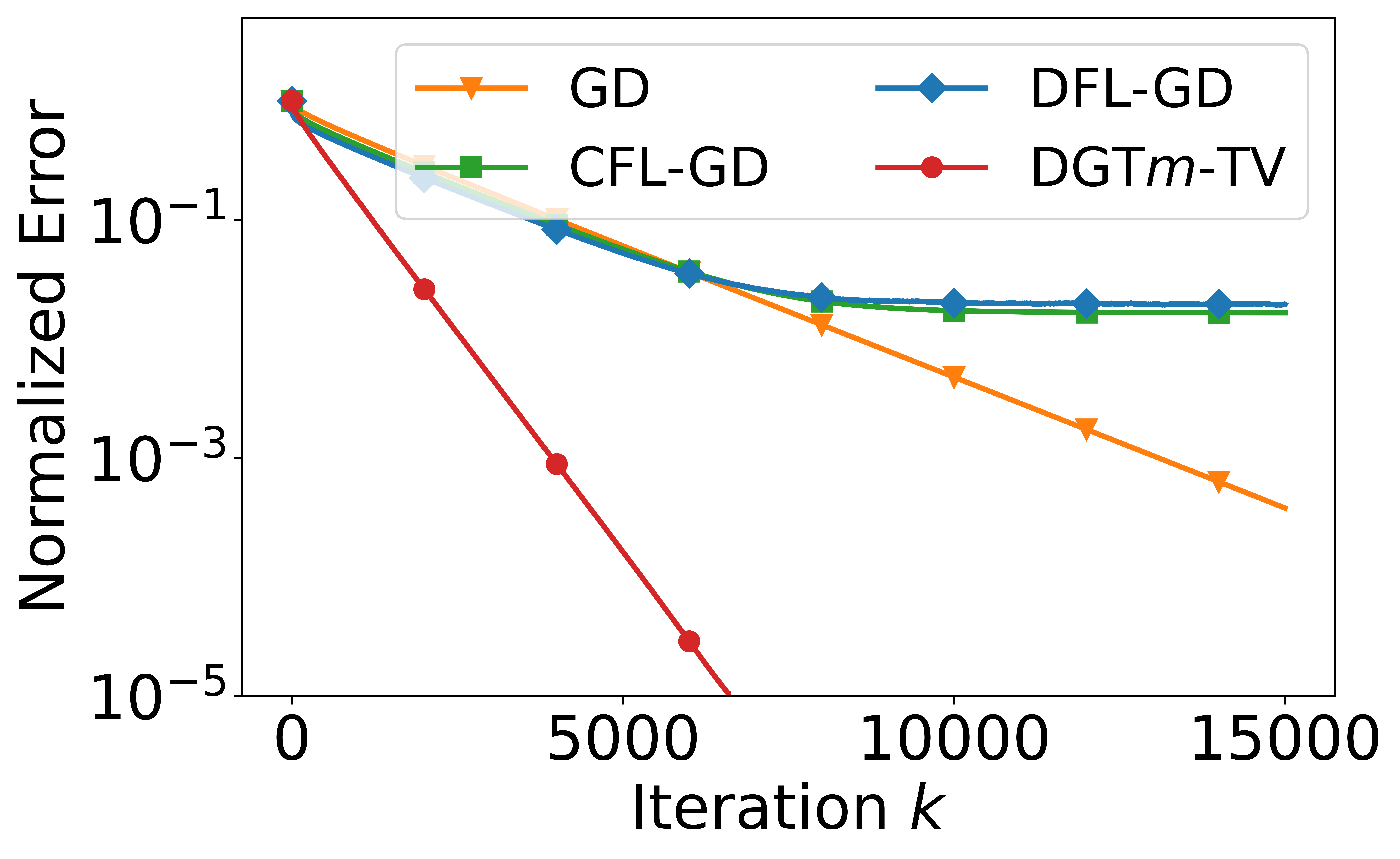}
\hspace*{-.4em}
\includegraphics[width=0.23\textwidth]{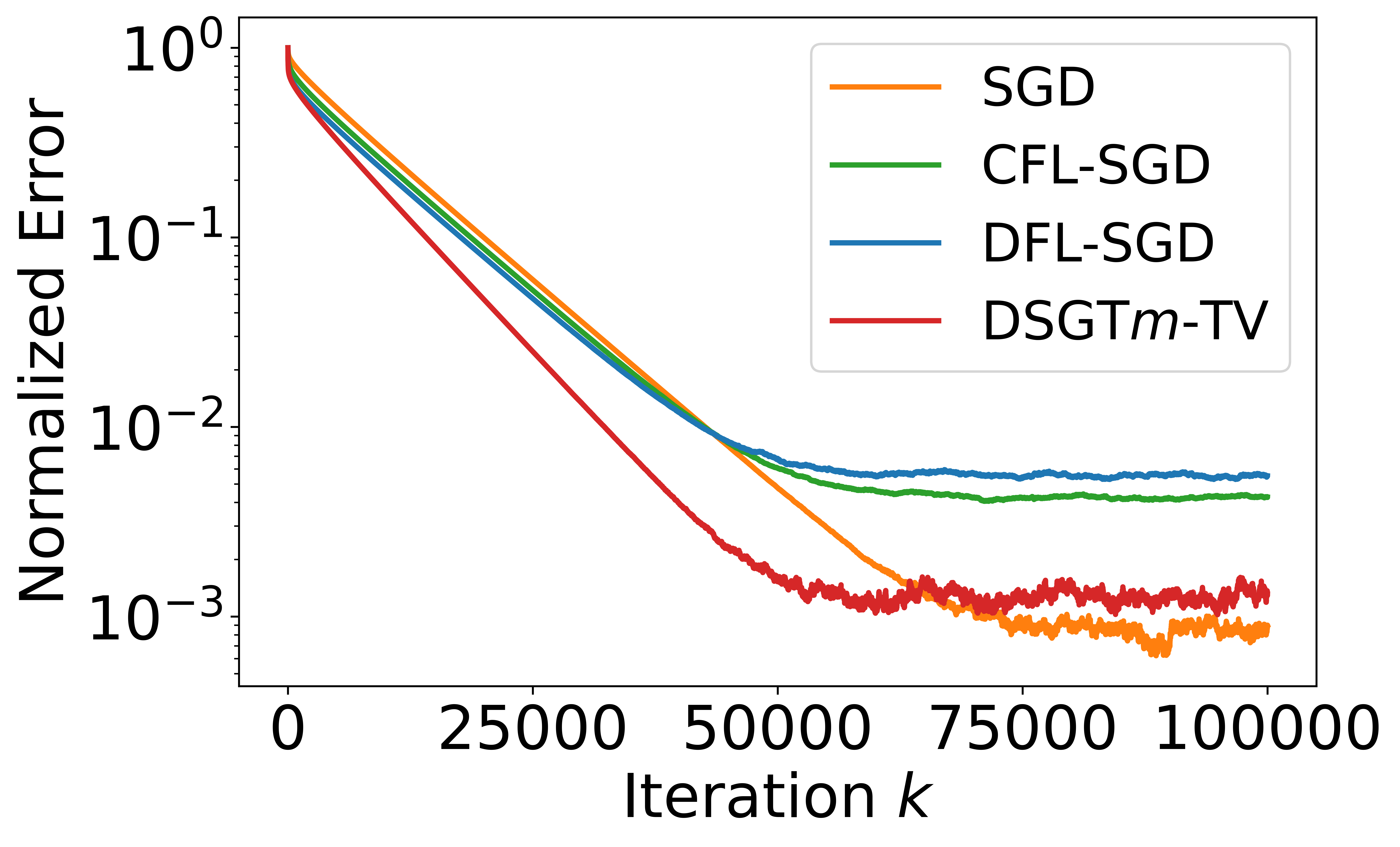}

\hspace*{-.3em}\includegraphics[width=0.22\textwidth]{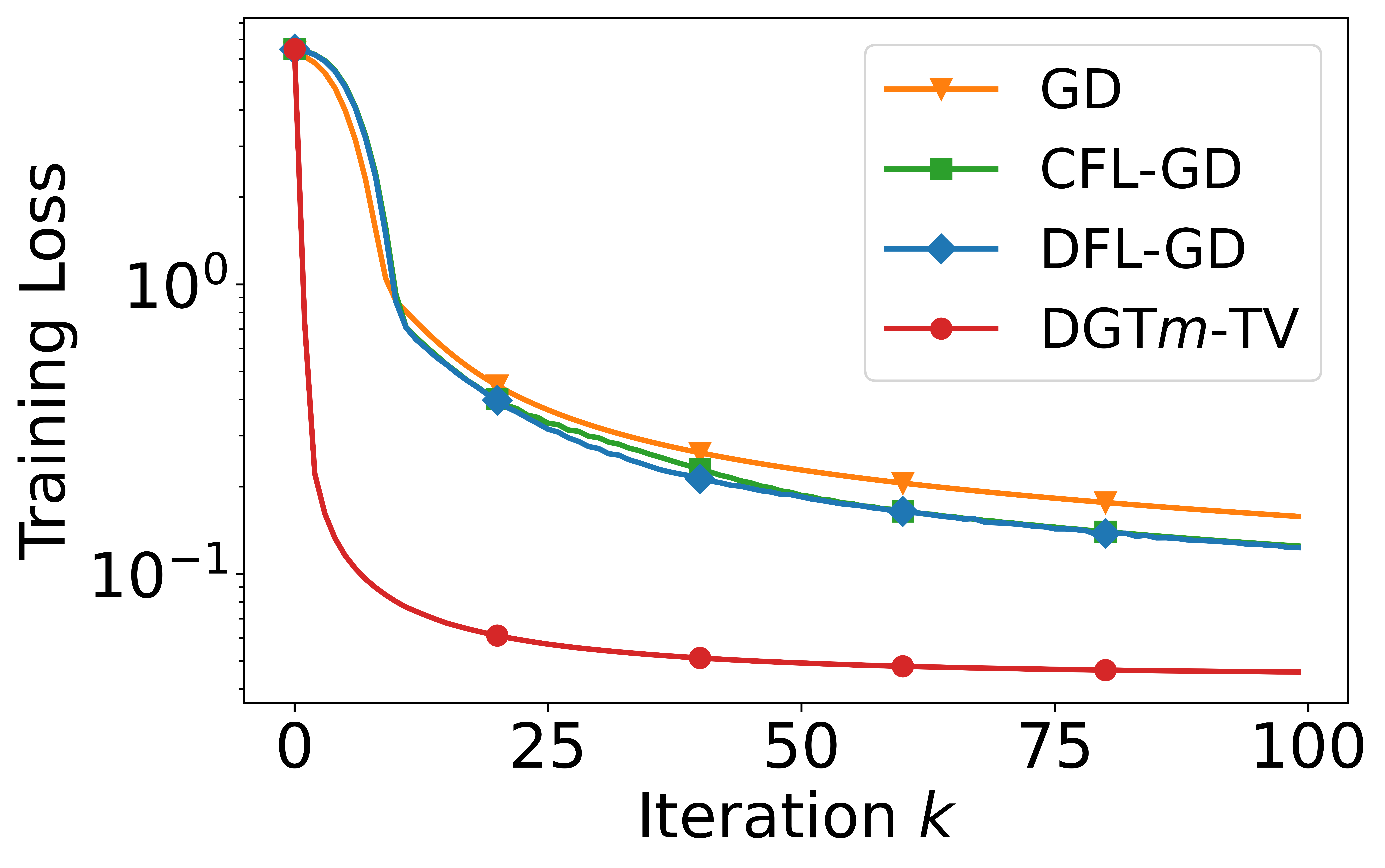}
\hspace*{-.2em}
\includegraphics[width=0.22\textwidth]{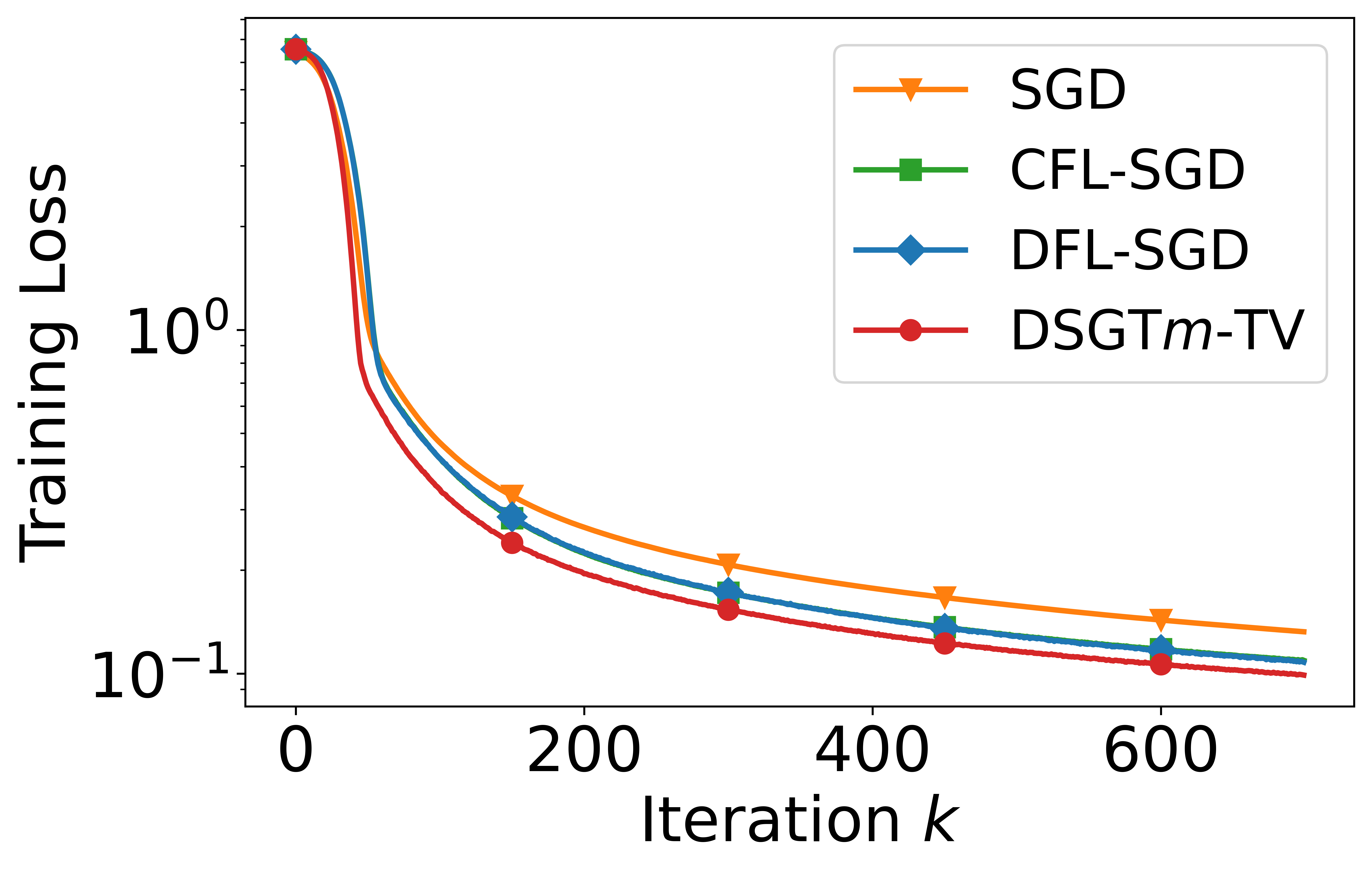}

\subfigure[Exact gradients ($\sigma^2=0$)]{\includegraphics[width=0.222\textwidth]{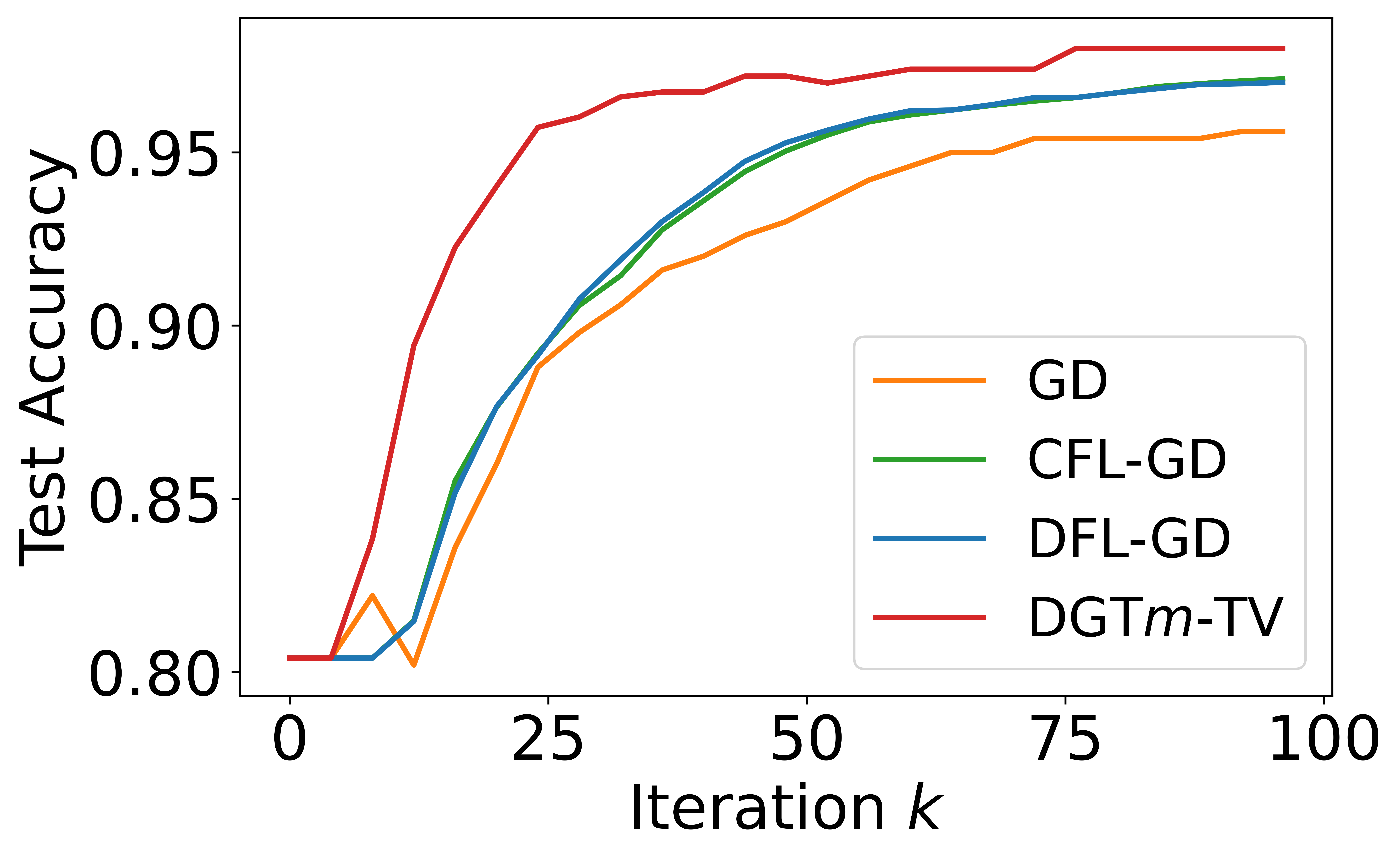}
}
\hspace*{-.4em}
\subfigure[Stochastic gradients ($\sigma^2>0$)]{\includegraphics[width=0.21\textwidth]{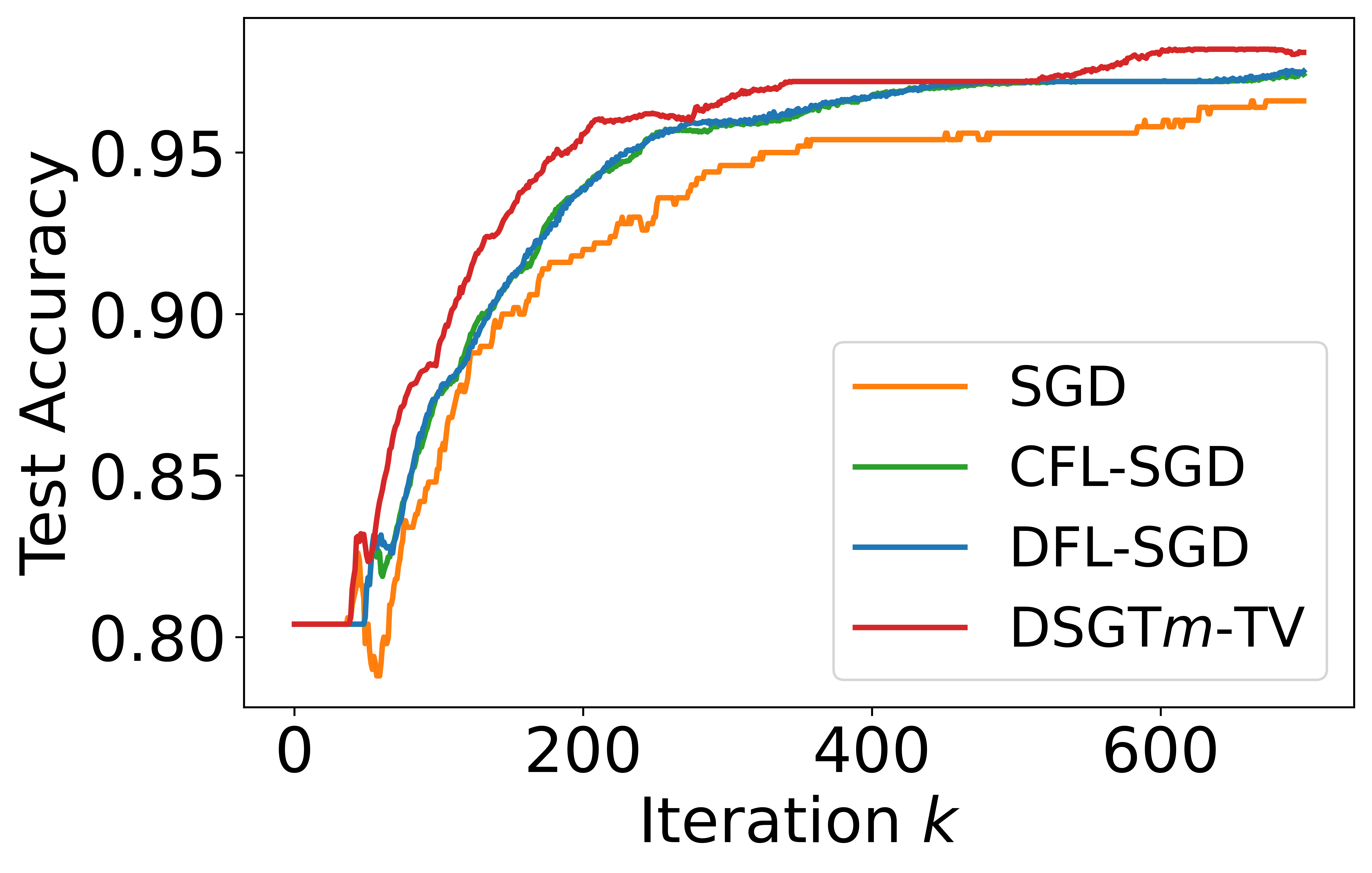}
}
\vspace{-0.2cm}
\caption{Performance on MNIST dataset. Top: Normalized error, Middle: Training loss,  Bottom: Accuracy on the test set.}
\label{fig:MNIST}
\vspace{-0.4cm}
\end{figure}

Fig.~\ref{fig:MNIST} compares the performance of different  algorithms on the MNIST dataset under both exact gradients ($\sigma^2=0$) and stochastic gradients ($\sigma^2>0$). We evaluate three performance metrics: (1) The normalized error, given by $\frac{\|\bx^k-\bx^*\|^2_2}{\|\bx^0-\bx^*\|^2_2}$, measures the relative error of the models $\bx^k$ produced by the algorithm at iteration $k$, compared to the optimal global model $\bx^*$, normalized by the initial error; (2) The training loss, computed as $\frac{1}{n}\sum_{i=1}^n \E_{\xi_{i} \sim \mathcal{D}_i} \mathcal{L}(x^i_k,\xi_{i})$, assesses the model's performance on the training data; and (3) The accuracy on the test set, calculated as $\textstyle \tfrac{\textnormal{Number of correct predictions}}{\textnormal{Total number of test samples}}$, provides insight into the model's effectiveness on unseen data.

\textbf{The Effect of Gradient Tracking:}
As expected, Fig.~\ref{fig:MNIST} (top left) shows that the \DGTmTV~algorithm converges to the exact optimal solution at a linear rate when using exact gradients,  similar to the performance of the GD algorithm. However, GD relies on centralized data aggregation at a single server, while the data in CFL--GD, DFL--GD and \DGTmTV~algorithms is distributed across multiple agents. In such distributed settings, other benchmarks like CFL--GD and DFL--GD only achieve convergence to a neighborhood around the optimal solution. When stochastic gradients are employed (Fig.~\ref{fig:MNIST} top right), all algorithms converge to a neighborhood of the optimum. However, the proposed \DSGTmTV~algorithm exhibits a smaller error in this neighborhood. As discussed in Section\ref{sec-GradientTracking}, this improvement is due to \DSGTmTV's ability to track the global gradient, ensuring that each agent's updates are representative of the global objective. In contrast, conventional DFL algorithms \cite{Tiancheng2022} only use local gradient information based on each agent’s local dataset. This approach can introduce significant bias, particularly with highly heterogeneous data distributions, leading to updates that may not accurately represent the global objective and increasing error.

\textbf{Performance Comparison:} The proposed algorithm, indicated by the red curve in the plots, demonstrates strong performance on the test set, achieving over $95\%$ accuracy after $25$ iterations when $\sigma^2=0$ (Fig.~\ref{fig:MNIST} bottom left) and reaching a similar level of accuracy after approximately $200$ iterations when $\sigma^2>0$ (Fig.~\ref{fig:MNIST} bottom right). Compared to other benchmarks, the proposed algorithm not only attains higher accuracy but also reduces training loss more rapidly. It is important to note that the same value of stepsizes are employed across all algorithms for a fair comparison. The observed acceleration in convergence can be attributed to the incorporation of the heavy-ball momentum. 
% Specifically, we use momentum parameters $\beta = 0.7$ and $\beta = 0.2$ for the two respective plots. A higher value of $\beta$, while adhering to the bounds established in Theorem~\ref{theo-main}, generally results in faster convergence.

\subsection{Application 2: Natural Language Processing}
We perform text classification by combining the proposed \DSGTmTV~algorithm with the pre-trained BERT (Bidirectional Encoder Representations from Transformers) model, aiming to effectively classify text data across various labels. Specifically, we use the SMS Spam Collection dataset, which consists of $5574$ SMS messages in English labeled as either legitimate or spam. Each message is tokenized and transformed into a high-dimensional vector representation using the pre-trained BERT-base model, which provides $768$-dimensional embeddings. We then augment the model with an additional classification layer for text classification tasks. This layer maps the high-dimensional embeddings produced by BERT into the label space of the classification problem. The classification layer is trained using either the proposed \DSGTmTV~algorithm or one of the benchmark algorithms, with the training dataset consisting of $5000$ samples and $550$ samples reserved for testing. We set the regularization parameter $\lambda$ to $0.001$.

\begin{figure}[t!]
	\centering
	\includegraphics[width=0.23\textwidth]{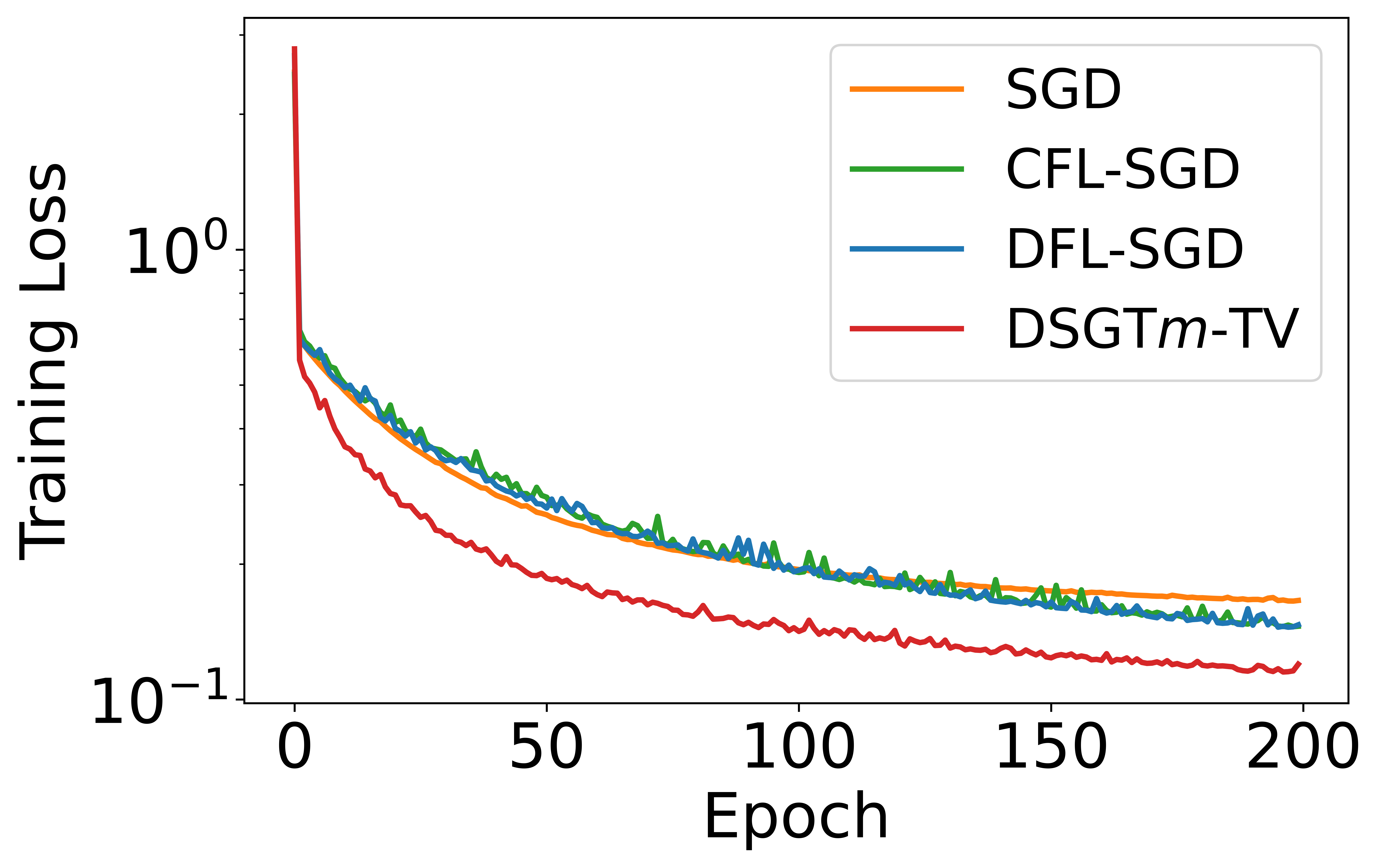}\label{fig:SpamLossStochastic} 
	   \hspace*{-.5em} 
	\includegraphics[width=0.23\textwidth]{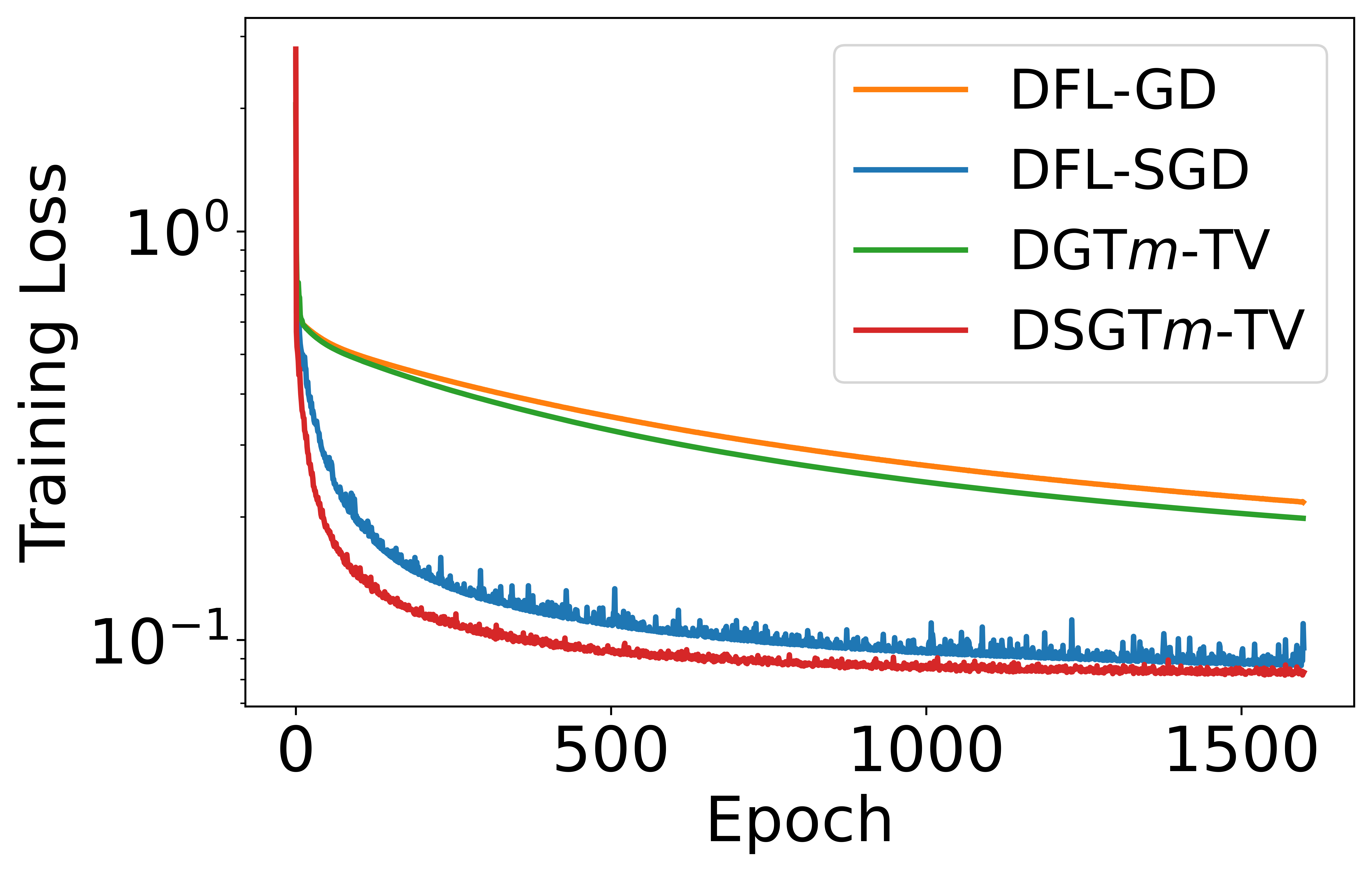}\label{fig:SpamLossSD}
 
\subfigure{
	\includegraphics[width=0.23\textwidth]{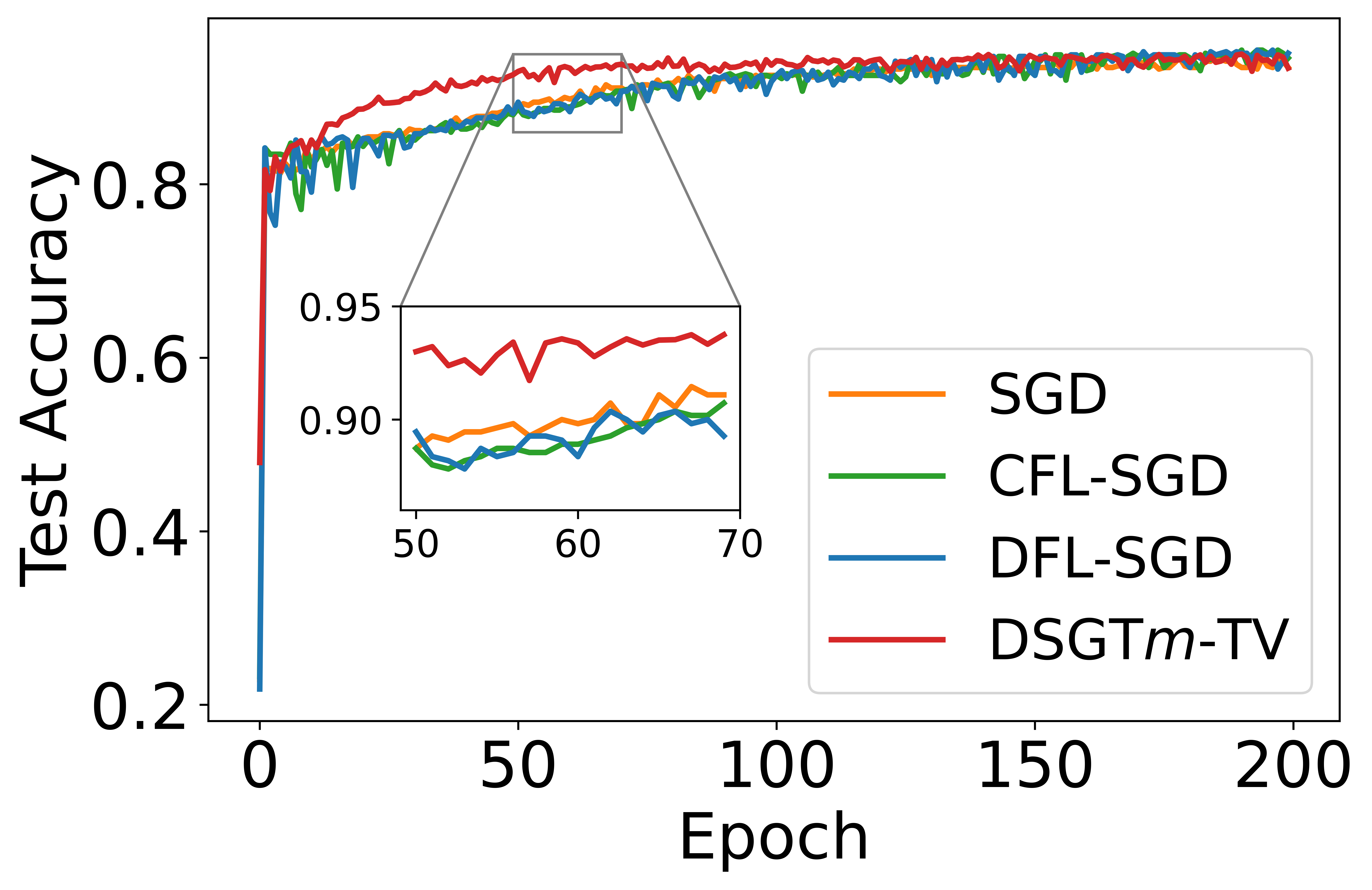}\label{fig:SpamStochastic}} 
	   \hspace*{-.5em} 
\subfigure{
	\includegraphics[width=0.23\textwidth]{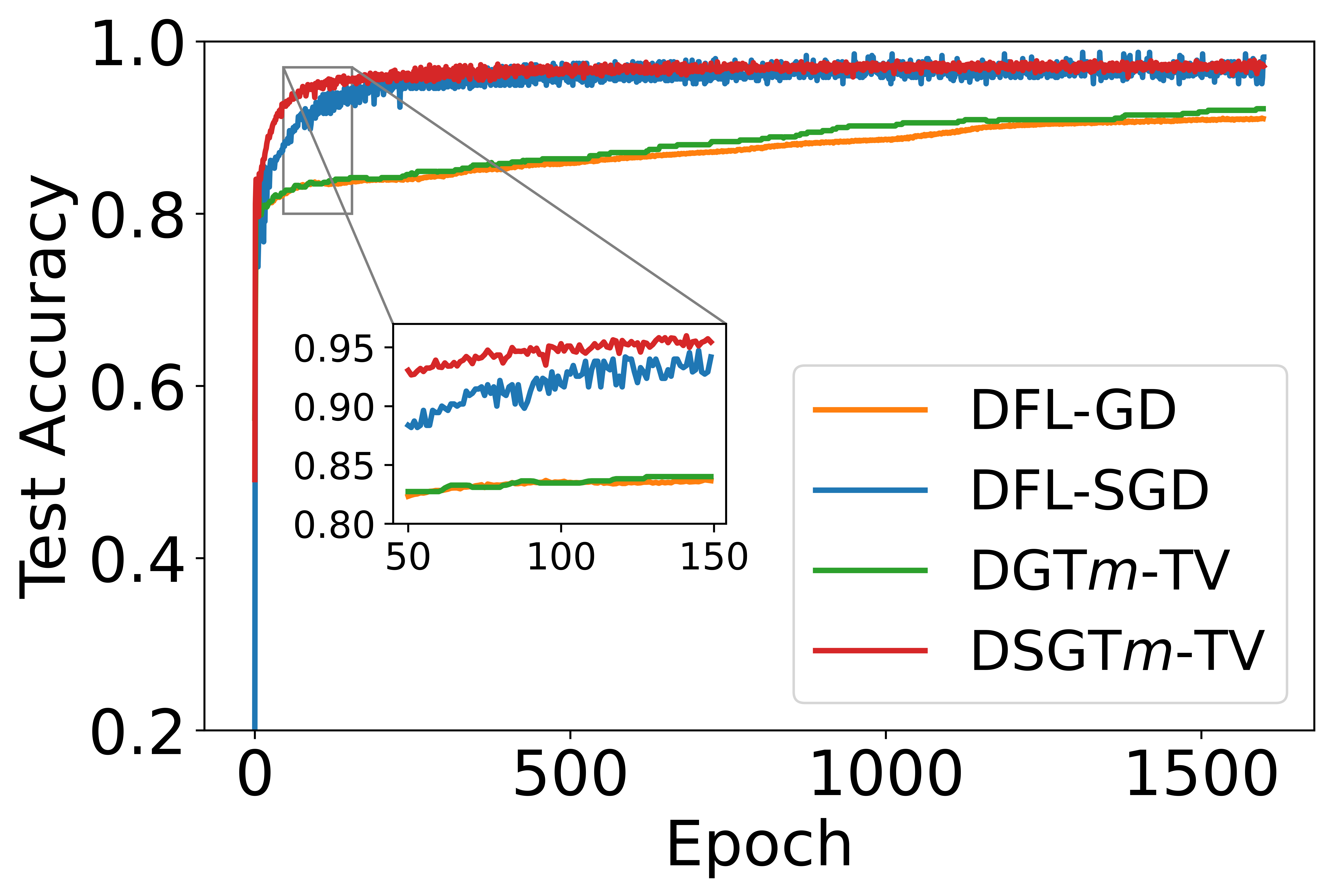}\label{fig:SpamSD}
	} 
 
	\caption{Performance on SMS Spam Dataset. Left: Performance comparison of \DSGTmTV~versus SGD algorithms accross different settings. Right: Comparison of decentralized FL algorithms: benchmark versus \DSGTmTV~with deterministic and stochastic gradient variants. }\label{fig:SMS Spam}
\vspace{-0.4cm}
\end{figure}

Fig.~\ref{fig:SMS Spam} presents the algorithms' performance in terms of training loss and test accuracy. 
The left panel compares the performance of the proposed \DSGTmTV~algorithm with that of several SGD variants under different settings. It is evident that \DSGTmTV~achieves higher accuracy and lower loss more rapidly compared to the benchmark algorithms. Specifically, \DSGTmTV~reaches an accuracy of $93\%$ after $50$ epochs, whereas the benchmark algorithms require over $100$ epochs to attain comparable results. Furthermore, among the two algorithms in the decentralized FL setting, \DSGTmTV~exhibits more stable performance compared to the state-of-the-art DFL--SGD. This increased stability can be attributed to the gradient tracking mechanism employed by \DSGTmTV, whereas DFL--SGD shows considerable oscillations in performance during the initial $25$ epochs.

The right panel further focuses on decentralized FL settings, comparing the proposed \DSGTmTV~algorithm with the DFL--SGD algorithm in both exact and stochastic gradient variants. The results confirm that the \DSGTmTV~algorithm outperforms the DFL--SGD benchmark. In addition, stochastic variants achieve over $95\%$ accuracy after $200$ epochs, whereas their deterministic counterparts require more than $1500$ epochs to attain comparable accuracy.
% This difference is primarily due to the stochastic gradient methods' ability to introduce randomness into the gradient estimation process, which accelerates convergence by allowing the algorithm to explore the parameter space more effectively. 
The inherent noise in the gradients can enable the algorithm to make more substantial updates and potentially escape from suboptimal regions, leading to faster convergence. In contrast, deterministic methods rely on precise gradient calculations, resulting in slower convergence as they make more conservative updates.

\section{Conclusions}
\label{sec:conc}
This paper has proposed \DSGTmTV, an accelerated consensus-based algorithm with gradient tracking for decentralized collaborative learning and optimization tasks among networked agents communicating over a time-varying directed network. Our theoretical analysis has established that \DSGTmTV~exhibits linear convergence to the exact global optimum when exact gradient information is available, and demonstrates linear convergence in expectation to a neighborhood of the global optimum under stochastic gradient estimations. Importantly, the stepsizes and momentum parameters can be fully uncoordinated, and we have derived explicit bounds for convergence guarantees that are independent of any heterogeneity, enabling each agent to determine these parameters autonomously in a fully decentralized manner. Consequently, this algorithm can be regarded as a unified framework that generalizes many existing methods and accommodates a wider range of communication connectivity structures. We have validated our theoretical results through image classification and natural language processing tasks. A promising avenue for future research is the investigation of variance-reduced techniques to address the inherent noise in gradient estimates and remove the steady-state error.

\bibliographystyle{IEEEtran}
\bibliography{references}

\newpage
\setcounter{page}{1}
\appendix
\subsection{Basic Results}\label{Appendix:basicre}
We leverage the following lemmas within our theoretical analysis:
% The following lemmas are crucial to our analysis: 
Lemma~\ref{lem-normlincomb} provides a relation for the squared norm of a linear combination of vectors.
% , while Lemma~\ref{lem-amatrices} and Lemma~\ref{lem-bmatrices} define and outline the fundamental properties of stochastic vectors associated with stochastic matrices.

\begin{lemma}[\!\!\cite{nguyen2022distributed}, Corollary~5.2]\label{lem-normlincomb}
Consider a vector collection $\{u_i, \, i\in[n]\}\subset\re^p$, and a scalar collection $\{\g_i,\, i\in[n]\}\subset\re$ such that $\sum_{i=1}^n \g_i=1$. For all $u\in \re^p$, we have
\begin{align*}
\Bigg\|\!\sum_{i=1}^n \!\g_i u_i - u \Bigg\|^2 \!\!\!= \!\sum_{i=1}^n \!\g_i \|u_i-u\|^2 \!-\!\!\sum_{i=1}^n \g_i \Bigg\|u_i \!-\! \Bigg(\sum_{j=1}^n \g_j u_j\Bigg)\!\Bigg \|^2\!\!\!.
\end{align*}
\end{lemma}

Consider a strongly connected directed graph $\bbG=([n],{\cal E})$, and weight matrices $A$ and $B$ that are aligned with the graph $\bbG$. The next lemmas establish the contraction of the stochastic matrices $A$ and $B$.
% with $A_{ij}>0$ when $j=i$ and $(j,i)\in {\cal E}$, and $A_{ij}=0$ otherwise
% $B_{ij}>0$ when $j=i$ and $(j,i)\in {\cal E}$, and $B_{ij}=0$ otherwise
\begin{lemma}[\!\!\cite{nguyen2022distributed}, Lemma 6.1]\label{lem-basic-xcontract}
Let $A$ be a row-stochastic matrix, $\phi$  be a stochastic vector and let $\pi$ be a nonnegative vector such that $\pi^{\T}A=\phi^{\T}$. Consider a collection of vectors $x_1,\ldots,x_n\in\re^p$. For $\hat x_\phi=\sum_{i=1}^n \phi_i x_i$, we have
\begin{align*}
\sqrt{\sum_{i=1}^n \pi_i\Bigg\|\sum_{j=1}^n A_{ij} x_j- \hat x_\phi\Bigg\|^2} &\le c\sqrt{\sum_{j=1}^n \phi_j \|x_j - \hat x_\phi \|^2},
\end{align*}
where $c=\sqrt{1 -   \frac{\min(\pi) (\min^+(A))^2}{\max^2(\phi)\mathsf{D}(\bbG)\mathsf{K}(\bbG)}}\in (0,1)$ is a scalar.
\end{lemma}

\begin{lemma}[\!\!\cite{Angelia2022AB}, Lemma 4.5]\label{lem-basic-ycontract}
Let $B$ be a column-stochastic matrix, $\nu$  be a stochastic vector with positive entries, i.e., $\nu_i>0,~\forall i\in[n]$, and let the vector $\pi$ be given by $\pi=B\nu$. Considering a collection of vectors $y_1,\ldots,y_n\in\re^p$, we have
\begin{align*}
	\sqrt{\!\sum_{i=1}^n\pi_i \Bigg\|\frac{\sum_{j=1}^n \!B_{ij} y_j}{\pi_i} \!-\! \sum_{j=1}^n y_j\Bigg\|^2}
	\!\le\!\tau
	\sqrt{\!\sum_{i=1}^n \nu_i \Bigg\|\frac{y_i}{\nu_i} \!- \!\sum_{j=1}^n y_j\Bigg\|^2}\!,
\end{align*}	
where $\tau=\sqrt{1 -   \frac{\min^2(\nu)\,(\min^+(B))^2}{\max^2(\nu) \max(\pi)\, \mathsf{D}(\bbG)\mathsf{K}(\bbG)}}\in (0,1)$ is a scalar.
\end{lemma}

The following result demonstrates the contraction property of the gradient mapping.
\begin{lemma}[\!\!\cite{pu2021stochastic}, Lemma 3]\label{lem-contraction} 
Under Assumption~\ref{asm-functions} and Assumption~\ref{asm-strconv}, we have for all $k\ge 0$,
\[\|h(\bx_k)-\nabla f(\hat{x}_k)\|\le \frac{L}{\sqrt{n}}\|\bx_k-\hat{\bx}_k\|.\]
Additionally if $0<\a<2/(\mu+L)$, then for all $x\in\re^p$,
\[\|x-x^*-\a \nabla f(x)\|\le (1-\a\mu) \|x-x^*\|.\]
\end{lemma}

Finally, we present a result from nonnegative matrix theory that will be used to prove the main convergence result.
\begin{lemma}[\!\!\cite{horn1990matrix}, Corollary 8.1.29] \label{lem-SpectralRadius}
Let $M\in\re^{n\times n}$ be a nonnegative matrix and $\delta\in\re^n$ be a positive vector. If $M\delta<\lambda \delta$ with $\lambda>0$, then $\rho(M)<\lambda$.
\end{lemma}

\subsection{Proof of Lemma~\ref{lem-weightedavg}} \label{Appendix-lem-weightedavg}
\begin{proof}
By taking a weighted average of $x_{k+1}^i$ using the $\phi_{k+1}$ weights from the update equation \eqref{eq-x} for $x_{k+1}^i$, we obtain
\begin{align*}
\hat x_{k+1} = \sum_{i=1}^n [\phi_{k+1}]_i x_{k+1}^i= \sum_{i=1}^n[\phi_{k+1}]_i\sum_{j=1}^n[A_k]_{ij}x_k^j \qquad\qquad\\
- \sum_{i=1}^n \a_i[\phi_{k+1}]_i y_k^i + \sum_{i=1}^n[\phi_{k+1}]_i \b_i(x_k^i-x_{k-1}^i).
\end{align*}
	Since $\phi_{k+1}^{\T}A_k=\phi_k^{\T}$, for the double-sum term we have
\begin{align*}
\sum_{i=1}^n [\phi_{k+1}]_i \sum_{j=1}^n [A_k]_{ij}x_k^j
&=\sum_{j=1}^n \Big(\sum_{i=1}^n [\phi_{k+1}]_i [A_k]_{ij}\Big) x_k^j\\
&=\sum_{j=1}^n [\phi_k]_j x_k^j = \hat x_{k},
\end{align*}
which completes the proof.
\end{proof}

\subsection{Proof of Lemma~\ref{lem-ysum}} \label{Appendix-Lemmalem-ysum}
\begin{proof}
Consider $\left\|\sum_{i=1}^n y_k^i\right\|$; we have 
\begin{align*}
&\left\| \sum_{i=1}^n y_k^i\right\|=\left\|n\left(\bby_k-h(\bx_k)\right)+n h(\bx_k)-\sum_{i=1}^n\nabla f_i(x^*) \right\|\\
\le&~n\left\|\bby_k-h(\bx_k)\right\|+\left\|\sum_{i=1}^n\left(\nabla f_i(x_k^i) -\nabla f_i(x^*)\right) \right\|,
\end{align*}
where we use \eqref{eq-quant-avg1}, \eqref{eq-quant-avg2}, and, the fact that $\sum_{i=1}^n\nabla f_i(x^*)=0$. By using the assumption that each $f_i$ has Lipschitz continuous gradients with a Lipschitz constant $L>0$, we obtain
\begin{align*}
&\left\|\sum_{i=1}^n\left(\nabla f_i(x_k^i) -\nabla f_i(x^*)\right) \right\|\le \sum_{i=1}^n\left\| \nabla f_i(x_k^i)-\nabla f_i(x^*)\right\|\\
\le&~ L\sum_{i=1}^n\| x_k^i -x^*\|\!=\! L\sqrt{n}\|\bx_k -\bx^*\|\!\le\! L\sqrt{n}\varphi_k\|\bx_k -\bx^*\|_{\phi_k}\!.
\end{align*}
Hence,
\begin{align*}
\left\| \sum_{i=1}^n y_k^i\right\|^2\le 2n^2\left\|\bby_k-h(\bx_k)\right\|^2+ 2nL^2\varphi_k^2\|\bx_k -\bx^*\|_{\phi_k}^2.
\end{align*}
Applying the Lemma~\ref{lem-normlincomb} with $u_i=x_k^i$, $\g_i=[\phi_k]_i$ for all $i$, $\hat x_k=\sum_{j=1}^n[\phi_k]_j x_{k}^j$ and $u=x^*$ yields
\[\|\bx_k -\bx^*\|_{\phi_k}^2=\|\hat x_k-x^*\|^2 + \|\bx_k -\hat{\bx}_k\|_{\phi_k}^2.\]
Combining the preceding two relations, taking the conditional expectation given $\F_k$ on both sides, and applying Lemma~\ref{lem-sumgrad}(c), we obtain the desired relation.
\end{proof}
\subsection{Proof of Proposition~\ref{prop-waverx}}\label{Appendix-Propositionprop-waverx}
\begin{proof}
Under Assumption~\ref{asm-strconv}, the unique minimizer $x^*$ of $f(x)$ over $x\in\re^p$ exists. 
% By subtracting the optimal point $x^*$ from both sides of the recursive relation for $\hat x_k$ in Lemma~\ref{lem-weightedavg}, and by adding and subtracting $n \phi_{k+1}^\T\Da\pi_k\bby_k =  \sum_{i=1}^n [\phi_{k+1}]_i [\pi_k]_i\a_i\sum_{j=1}^ny_{k}^j $, we obtain
From Lemma~\ref{lem-weightedavg},  noticing that $n \phi_{k+1}^\T\Da\pi_k\bby_k =  \sum_{i=1}^n [\phi_{k+1}]_i [\pi_k]_i\a_i\sum_{j=1}^ny_{k}^j $, we obtain
\begin{align*}
&\hat x_{k+1}-x^*\\
=&\left(\hat x_k - x^*  - n \phi_{k+1}^\T\Da\pi_k\bby_k\right) +\sum_{i=1}^n [\phi_{k+1}]_i\b_i(x_k^i-x_{k-1}^i)\\
&-\sum_{i=1}^n [\phi_{k+1}]_i [\pi_k]_i\a_i\Bigg(\frac{y_k^i}{[\pi_k]_i} - \sum_{j=1}^ny_{k}^j\Bigg).
\end{align*}
In order to make the notation more concise, we denote $\tilde{\a}_k=n \phi_{k+1}^\T\Da\pi_k$, 
% and $\bar{\a}_k=n \bar{\a}\phi_{k+1}^\T\pi_k$, 
for all $k \ge 0$. Taking norms, squaring both sides and using Young's inequality that $(a+b)^2\le (1+\varrho)a^2+\left(1+\frac{1}{\varrho}\right)b^2$ with $\varrho=\frac{\tilde{\a}_k\mu}{2}>0$, yields
\begin{align}\label{eq-xhatxstar1}
\|\hat x_{k+1}-x^*\|^2&\le \left(1+\tfrac{\tilde{\a}_k\mu}{2}\right)\|q_1\|^2\nonumber\\
&+\left(1+\tfrac{2}{\tilde{\a}_k\mu}\right)\left(2\bar{\b}^2\|q_2\|^2+2\bar{\a}^2\| q_3\|^2\right),
\end{align}
where $q_1 \!\triangleq\! \hat x_k \!-\! x^*  \!- \! \tilde{\a}_k\bby_k,~ q_2 \!\triangleq\! \sum_{i=1}^n [\phi_{k+1}]_i(x_k^i \!-\!x_{k-1}^i),$
and $q_3 \triangleq \sum_{i=1}^n [\phi_{k+1}]_i [\pi_k]_i\left(\frac{y_k^i}{[\pi_k]_i} - \sum_{j=1}^ny_{k}^j\right).$\\

\noindent
\textbf{Bounding $\E\left[\|q_1\|^2|\F_k\right]$}: We have
\begin{align*}
\|q_1\|^2 =& \left\|\hat x_k - x^* - \tilde{\a}_k f(\hat x_k) +\tilde{\a}_k \left(f(\hat x_k)-\bby_k\right)\right\|^2 \\
=& \|\hat x_k - x^* - \tilde{\a}_k f(\hat x_k)\|^2+\tilde{\a}_k^2 \|f(\hat x_k)-\bby_k\|^2\\
&+2\tilde{\a}_k \langle \hat x_k - x^* - \tilde{\a}_k f(\hat x_k),f(\hat x_k)-\bby_k\rangle.
\end{align*}
% while for regular graphs, for example, we do have $n[\pi_k]_i=1$ for all $i$ and $k$
Given Assumption~\ref{asm-functions} and Assumption~\ref{asm-strconv}, the function $f$ is $\mu$-strongly convex and its gradients are $L$-Lipschitz continuous.
Thus, if $\tilde{\a}_k\in (0,\tfrac{2}{\mu+L})$, according to Lemma~\ref{lem-contraction}, we obtain
\begin{align}\label{eq-StronglyConvRelation}
\left\|\hat x_k - x^* - \tilde{\a}_k f(\hat x_k) \right\| \le (1-\tilde{\a}_k\mu)\|\hat x_k - x^*\|.
\end{align}
% Hence,
% \begin{align*}
% \|q_1\|^2 \le& (1-\tilde{\a}_k\mu)^2\|\hat x_k - x^*\|^2+\tilde{\a}_k^2 \|f(\hat x_k)-\bby_k\|^2\\
% &+2\tilde{\a}_k \langle \hat x_k - x^* - \tilde{\a}_k f(\hat x_k),f(\hat x_k)-\bby_k\rangle.
% \end{align*}
% Notice that $\E[\bby_k|\F_k] = h(\bx_k)$ from Lemma~\ref{lem-sumgrad}(b), taking the conditional expectation on both sides of the preceding relation given $\F_k$, we obtain
Notice that $\E[\bby_k|\F_k] = h(\bx_k)$ from Lemma~\ref{lem-sumgrad}(b),  we obtain
\begin{align*}
\E\!\left[\|q_1\|^2|\F_k\!\right] \!\!\le& (1\!-\!\tilde{\a}_k\mu)^2\|\hat x_k \!-\! x^*\|^2\!\!+\!\tilde{\a}_k^2 \E\!\left[\|f(\hat x_k)\!-\!\bby_k\|^2|\F_k\!\right]\\
&+2\tilde{\a}_k \langle \hat x_k - x^* - \tilde{\a}_k f(\hat x_k),f(\hat x_k)-h(\bx_k)\rangle.
\end{align*}
To estimate the inner product, we use the Cauchy–Schwarz inequality and the relation in \eqref{eq-StronglyConvRelation} to derive
\begin{align*}
&2\tilde{\a}_k\langle \hat x_k - x^* - \tilde{\a}_k f(\hat x_k),f(\hat x_k)-h(\bx_k)\rangle\\ 
\le & 2\tilde{\a}_k(1-\tilde{\a}_k\mu)\|\hat x_k - x^*\| \|f(\hat x_k)-h(\bx_k)\|\\
\le &\tilde{\a}_k\mu(1-\tilde{\a}_k\mu)^2\|\hat x_k - x^*\|^2+\tfrac{\tilde{\a}_k}{\mu}\|f(\hat x_k)-h(\bx_k)\|^2.
\end{align*}
Therefore, 
\begin{align}\label{eq-q1rel}
\E\!\left[\|q_1\|^2|\F_k\right] \!\le& (1\!-\!\tilde{\a}_k\mu)\|\hat x_k \!-\! x^*\|^2\!\!+\!\tilde{\a}_k^2 \E\!\left[\|f(\hat x_k)\!-\!\bby_k\|^2|\F_k\right]\nonumber\\
&+\tfrac{\tilde{\a}_kL^2\varphi_k^2}{n\mu}\|\bx_k-\hat{\bx}_k\|_{\phi_k}^2,
\end{align}
where we use the fact that $(1+\tilde{\a}_k\mu)(1-\tilde{\a}_k\mu)<1$ to obtain the first term and Lemma~\ref{lem-contraction} together with the relation in \eqref{eq-NormIneq} to obtain the last term. We now bound $\E\left[\|f(\hat x_k)-\bby_k\|^2|\F_k\right]$: 
\begin{align*}
\E\!\left[\|f(\hat x_k)\!-\!\bby_k\|^2|\F_k\!\right]
\!=\E\!\left[\|f(\hat x_k)\!-\!h(\bx_k)\!+\!h(\bx_k)\!-\!\bby_k\|^2|\F_k\right]\\
=\|f(\hat x_k)-h(\bx_k)\|^2+\E\left[\|h(\bx_k)-\bby_k\|^2|\F_k\right],
\end{align*}
% where the inner product term in the last equality equals to $0$ since
where the inner product term equals to $0$ since $\E[\bby_k|\F_k] = h(\bx_k)$. Using Lemma~\ref{lem-contraction} and Lemma~\ref{lem-sumgrad}(c) yields
\begin{align*}
&\E\left[\|f(\hat x_k)-\bby_k\|^2|\F_k\right]
\le \tfrac{L^2\varphi_k^2}{n}\|\bx_k-\hat{\bx}_k\|_{\phi_k}^2+\tfrac{\sigma^2}{n}.
\end{align*}
Substituting the preceding relation into \eqref{eq-q1rel}, we further obtain
\begin{align}\label{eq-q1relbound}
\E\left[\|q_1\|^2|\F_k\right]  \le& (1-\tilde{\a}_k\mu)\|\hat x_k - x^*\|^2\\
+&\tfrac{\tilde{\a}_k^2L^2\varphi_k^2}{n}\left(\tfrac{1}{\tilde{\a}_k\mu}+1\right)\|\bx_k-\hat{\bx}_k\|_{\phi_k}^2+\tfrac{\tilde{\a}_k^2\sigma^2}{n}.\nonumber
\end{align}
\textbf{Bounding $\E\left[\|q_2\|^2|\F_k\right]$}: Since $\phi_{k+1}$ is stochastic, we have 
% \begin{align*}
% \!\!\!\|q_2\|\!=\!&\left\|\sum_{i=1}^n[\phi_{k+1}]_i\left(x_k^i\!-x_{k-1}^i\right)\right\|\!\le\!\sum_{i=1}^n[\phi_{k+1}]_i\|x_k^i\!-x_{k-1}^i\|\\
% \le& \sqrt{\!\sum_{i=1}^n[\phi_{k+1}]_i\|x_k^i-x^i_{k-1}\|^2} \le  \|\bx_k-\bx_{k-1}\|.
% \end{align*}
% \begin{align}\label{eq-q2relbound}
% \!\!\!\!\!\!\!\!\text{Thus,}~~~\quad \quad \E\left[\|q_2\|^2|\F_k\right] \le \|\bx_k-\bx_{k-1}\|^2.\qquad \quad 
% \end{align}
\begin{align*}
\!\!\!\|q_2\|\!=&\!\left\|\sum_{i=1}^n[\phi_{k+1}]_i\left(x_k^i\!-x_{k-1}^i\right)\right\|\!\le\!\sum_{i=1}^n[\phi_{k+1}]_i\|x_k^i\!-x_{k-1}^i\|\nonumber\\
&\le \sqrt{\!\sum_{i=1}^n[\phi_{k+1}]_i\|x_k^i-x^i_{k-1}\|^2} \le  \|\bx_k-\bx_{k-1}\|.
\end{align*}
Thus,
\begin{align}\E\left[\|q_2\|^2|\F_k\right] \le \|\bx_k-\bx_{k-1}\|^2.\label{eq-q2relbound}
\end{align}
\textbf{Bounding $\E\left[\|q_3\|^2|\F_k\right]$}: 
% Using the relation in Lemma~\ref{lem-normlincomb} with 
Using Lemma~\ref{lem-normlincomb} with $\g_i=[\phi_{k+1}]_i$, $u_i=[\pi_k]_i\left(\frac{y_k^i}{[\pi_k]_i} - \sum_{j=1}^ny_{k}^j\right)$ and $u=0$, it follows that
\begin{align*}
&\|q_3\|^2\le\sum_{i=1}^n[\phi_{k+1}]_i\Bigg\|[\pi_k]_i\Bigg(\frac{y_k^i}{[\pi_k]_i} - \sum_{j=1}^ny_{k}^j\Bigg)\Bigg\|^2\\
\le&  \sum_{i=1}^n\Bigg\|[\pi_k]_i\Bigg(\frac{y_k^i}{[\pi_k]_i} - \!\sum_{j=1}^ny_{k}^j\Bigg)\Bigg\|^2\!\!\le\! \sum_{i=1}^n[\pi_k]_i\Bigg\|\frac{y_k^i}{[\pi_k]_i} -\! \sum_{j=1}^ny_{k}^j\Bigg\|^2\!\!,
\end{align*}
where the third inequality follows from Lemma~\ref{lem-normlincomb}.

Using the notation for $S(\by_{k},\pi_k)$ in \eqref{eq-x-S-quants} and taking the conditional expectation given $\F_k$ yields
\begin{align} \label{eq-q3relbound}
\E[\|q_3\|^2|\F_k]\le \E[S^2(\by_{k},\pi_k)|\F_k].
\end{align}
\textbf{Bounding $\E\left[\|\hat x_{k+1}-x^*\|^2|\F_k\right]$:}
Combining the relations in \eqref{eq-xhatxstar1}, 
\eqref{eq-q1relbound}--\eqref{eq-q2relbound},
% \eqref{eq-q1relbound}, \eqref{eq-q3relbound}, \eqref{eq-q1q3relbound} and \eqref{eq-q2relbound}, 
we obtain
\begin{align}\label{eq-xhatxstarcomp}
&\E\left[\|\hat x_{k+1}-x^*\|^2|\F_k\right]\\
\le& \left(1+\tfrac{\tilde{\a}_k\mu}{2}\right)(1-\tilde{\a}_k\mu)\|\hat x_k - x^*\|^2+\left(1+\tfrac{\tilde{\a}_k\mu}{2}\right)\tfrac{\tilde{\a}_k^2\sigma^2}{n}\nonumber\\
&+\tfrac{\tilde{\a}_k^2L^2\varphi_k^2}{n}\left(1+\tfrac{\tilde{\a}_k\mu}{2}\right)\left(\tfrac{1}{\tilde{\a}_k\mu}+1\right)\|\bx_k-\hat{\bx}_k\|_{\phi_k}^2\nonumber\\
&+2\left(1+\tfrac{2}{\tilde{\a}_k\mu}\right)\left(\bar{\b}^2\|\bx_k-\bx_{k-1}\|^2+\bar{\a}^2\E[S^2(\by_{k},\pi_k)|\F_k]\right).\nonumber
\end{align}
For $\tilde{\a}_k\in (0,\tfrac{2}{\mu+L})$, we have $2\tilde{\a}_k\mu\le \tilde{\a}_k(\mu+L)<2$ since $\mu\le L$. Thus, $\tilde{\a}_k\mu\le 1$ and $1\le\frac{1}{\tilde{\a}_k\mu}$. Hence, we have 
\begin{align*}
1\!+\!\tfrac{\tilde{\a}_k\mu}{2}\!\le\!\tfrac{3}{2},~ \tfrac{1}{\tilde{\a}_k\mu}\!+\!1 \!\le\! \tfrac{2}{\tilde{\a}_k\mu},~\left(\!1\!+\!\tfrac{\tilde{\a}_k\mu}{2}\!\right)\!(1\!-\!\tilde{\a}_k\mu)\!\le\! 1\!-\! \tfrac{\tilde{\a}_k\mu}{2}.
\end{align*}
Using the preceding bounds in \eqref{eq-xhatxstarcomp}, we obtain
\begin{align*}
&\E\left[\|\hat x_{k+1}-x^*\|^2|\F_k\right]\\
\le& \left(1- \tfrac{\tilde{\a}_k\mu}{2}\right)\|\hat x_k - x^*\|^2+\tfrac{3\tilde{\a}_kL^2\varphi_k^2}{n\mu}\|\bx_k-\hat{\bx}_k\|_{\phi_k}^2\nonumber\\
+&\tfrac{6}{\tilde{\a}_k\mu}\left(\bar{\b}^2\|\bx_k-\bx_{k-1}\|^2+\bar{\a}^2\E[S^2(\by_{k},\pi_k)|\F_k]\right)+\tfrac{3\tilde{\a}_k^2\sigma^2}{2n}.
\end{align*}
Recalling that $ \tilde{\a}_k = n \phi_{k+1}^\T\Da\pi_k$, we have   
\[n\bar{\a}\eta\le\tilde{\a}_k \le n\bar{\a} \phi_{k+1}^\T\pi_k,\] 
where $\eta\in(0,1)$ is the lower bound for $\phi_{k+1}^\T\pi_k$ as defined in \eqref{eq-const-max} and $\phi_{k+1}^\T\pi_k \le 1$. Thus, we obtain the desired relation.
\end{proof}

\subsection{Proof of Proposition~\ref{prop-xcontract}}\label{Appendix-Propositionprop-xcontract}
\begin{proof}
We define $u_k^i=x_k^i-x_{k-1}^i$ and $\bu_k=\bx_{k}-\bx_{k-1}$.
% \begin{align*}
% &\hat\bu_k=\left(\sum_{j=1}^n[\phi_{k+1}]_j u_{k}^j,\ldots,\sum_{j=1}^n[\phi_{k+1}]_j u_{k}^j\right),\\
% \text{and} &\bv_k=\left(\sum_{j=1}^n[\phi_{k+1}]_j y_{k}^j,\ldots,\sum_{j=1}^n[\phi_{k+1}]_j y_{k}^j\right).
% \end{align*}
Let $\boldsymbol{\Phi}_{k}=\one_n \phi_{k}^\T$, for all $k \ge 0$. We can write the relation for the weighted averages in Lemma~\ref{lem-weightedavg} in compact representation:
\[\hat \bx_{k+1}=\hat \bx_k -\boldsymbol{\Phi}_{k+1}\Da \by_k + \boldsymbol{\Phi}_{k+1}\Db\bu_k.\]
Recall the compact form of $x$-iterate process (see~\eqref{eq-x-comp}),
\[\bx_{k+1} =  A_k\bx_k  - \Da \by_{k}+\Db\bu_k \qquad\hbox{for all $k\ge 0$}.\]
Upon subtracting the preceding relations.
% , we obtain
% \begin{align*}
% \bx_{k+1} \!- {\hat \bx}_{k+1} = \left(A_k\bx_k \!- \hat \bx_k\right) - &(\mbi -\boldsymbol{\Phi}_{k+1})\Da\by_{k}\\
% + &(\mbi -\boldsymbol{\Phi}_{k+1})\Db\bu_{k}.
% \end{align*}
Taking the $\phi_{k+1}$-norm, squaring both sides of the preceding relation and using Young's inequality with $\varrho>0$, we obtain
% $(1+\varrho)a^2+\left(1+\frac{1}{\varrho}\right)b^2$ with $\varrho=\frac{\tilde{\a}_k\mu}{2}>0$
\begin{align}\label{wq-Dquant1}
&\left\|\bx_{k+1}-\hat \bx_{k+1}\right\|_{\phi_{k+1}}^2 \le (1+\varrho)\left\|A_k\bx_k - \hat \bx_k\right\|_{\phi_{k+1}}^2\!+\\
&\!\!2\!\left(\!1\!+\!\tfrac{1}{\varrho}\!\right)\!\!\left(\!\left\|(\mbi \!-\!\boldsymbol{\Phi}_{k+1})\Da\by_{k}\right\|_{\phi_{k+1}}^2 \!\!+\!\left\|(\mbi \!-\!\boldsymbol{\Phi}_{k+1})\Db\bu_{k} \right\|_{\phi_{k+1}}^2\!\right)\!.\!\!\nonumber
\end{align}
% \begin{align}\label{wq-Dquant1}
% \!\!\left\|\bx_{k+1}-\hat \bx_{k+1}\right\|_{\phi_{k+1}}^2\!\!=&\left\|A_k\bx_k - \hat \bx_k\right\|_{\phi_{k+1}}^2\!\!+\!\a^2\!\left\|\by_{k} \!-\! \bv_k\right\|_{\phi_{k+1}}^2 \nonumber\\
% &\!\!-2\a\langle A_k\bx_k - \hat \bx_k,\by_{k} - \bv_k\rangle_{\phi_{k+1}}.\!\!
% \end{align}

For the first term on the right-hand-side (RHS), we can write explicitly in terms of the vector components with $z_k^i=\sum_{j=1}^n [A_k]_{ij} x_{k}^j$, and obtain 
% For the first term, let $z_k^i=\sum_{j=1}^n [A_k]_{ij} x_{k}^j$, we have 
\begin{align*}
\left\|A_k\bx_k - \hat \bx_k\right\|_{\phi_{k+1}}^2= \sum_{i=1}^n [\phi_{k+1}]_i \|z_k^i -\hat x_k\|^2. 
\end{align*}
Next, we note that the vectors $z_k^i$, $i\in[n]$, satisfy Lemma~\ref{lem-basic-xcontract}, with $A=A_k$, and $x_i=x_k^i$ for all $i\in[n]$. Moreover, since we have $\phi_{k+1}^{\T}A_k=\phi_k^{\T}$ by Lemma~\ref{lem-amatrices}, Lemma~\ref{lem-basic-xcontract} applies with $\pi=\phi_{k+1}$, $\phi=\phi_k$ and $\hat x_\phi=\hat x_k$, which yields
% \begin{equation*}
% \left\|A_k\bx_k - \hat \bx_k\right\|_{\phi_{k+1}}^2 
% \le c_k^2
% \left\|\bx_k - \hat \bx_k\right\|_{\phi_{k}}^2, 
% \end{equation*}
% where we use $\min(A_k^+)\ge a$ and $c_k$ as denoted in \eqref{eq-const-all}. Using the upper bound for $c_k$ as defined in \eqref{eq-const-max}, we have
\begin{equation}\label{eq-zxhat}
\left\|A_k\bx_k - \hat \bx_k\right\|_{\phi_{k+1}}^2 
\le c^2
\left\|\bx_k - \hat \bx_k\right\|_{\phi_{k}}^2.
\end{equation}

For the second term in \eqref{wq-Dquant1}, using relation in Lemma~\ref{lem-normlincomb} with $\g_i=[\phi_{k+1}]_i$, $u_i=y_k^i$ for all $i$ and $u=0$, it follows that
\begin{align}\label{eq-Drel2}
&\!\!\!\left\|(\mbi \!-\!\boldsymbol{\Phi}_{k+1})\Da\by_{k} \right\|_{\phi_{k+1}}^2\!\!=\!\!\sum_{i=1}^n [\phi_{k+1}]_i\a_i^2 \Bigg\|y_k^i \!-\!\! \sum_{j=1}^n[\phi_{k+1}]_j y_{k}^j\Bigg\|^2\nonumber\\
\le& \bar{\a}^2\sum_{i=1}^n [\phi_{k+1}]_i \left\|y_k^i \right\|^2
= \bar{\a}^2\sum_{i=1}^n [\phi_{k+1}]_i[\pi_k]_i \frac{\left\|y_k^i \right\|^2}{[\pi_k]_i}\nonumber\\ 
\le& \bar{\a}^2\max_{j\in[n]} ([\phi_{k+1}]_j[\pi_k]_j) \sum_{i=1}^n\frac{\left\|y_k^i \right\|^2}{[\pi_k]_i}=\g_k^2 \bar{\a}^2\|\by_{k}\|_{\pi_k^{-1}}^2,
\end{align}
with $\g_k\!=\!\!\displaystyle\sqrt{\max_{i\in[n]} [\phi_{k+1}]_i[\pi_k]_i}\!\le\! 1$ and  $\|\by_{k}\|_{\pi_k^{-1}}^2\!=\!\sum_{i=1}^n \!\!\frac{\|y_k^i\|^2}{[\pi_k]_i}$.

For the third term in \eqref{wq-Dquant1}, using the relation in Lemma~\ref{lem-normlincomb} with $\g_i=[\phi_{k+1}]_i$, $u_i=u_k^i$ and $u=0$, it follows that
\begin{align}\label{eq-Drel3}
&\left\|(\mbi \!-\!\boldsymbol{\Phi}_{k+1})\Db\bu_k \right\|_{\phi_{k+1}}^2 \!\!=\!\! \sum_{i=1}^n [\phi_{k+1}]_i\b_i \Bigg\|u_k^i \!-\!\! \sum_{j=1}^n[\phi_{k+1}]_j u_k^j\Bigg\|^2 \nonumber \\
&\le  \bar{\b}^2\sum_{i=1}^n [\phi_{k+1}]_i \|u_k^i \|^2 \le \bar{\b}^2\|\bx_k-\bx_{k-1} \|^2,
\end{align}

\noindent 
Finally, substituting the relations in \eqref{eq-zxhat}, \eqref{eq-Drel2} and \eqref{eq-Drel3} back into \eqref{wq-Dquant1} and choosing $\varrho=\frac{1-c^2}{2c^2}$, we have
\begin{align*}
\left\|\bx_{k+1}\!-\!\hat \bx_{k+1}\right\|_{\phi_{k+1}}^2
\!\le&\tfrac{1+c^2}{2}\left\|\bx_k \!-\! \hat \bx_k\right\|_{\phi_{k}}^2\!\!+\!\tfrac{2\bar{\a}^2 (1+c^2)\g_k^2}{1-c^2} \|\by_{k}\|_{\pi_k^{-1}}^2\\
&+\tfrac{2\bar{\b}^2 (1+c^2)}{1-c^2} \|\bx_k-\bx_{k-1} \|^2.
\end{align*}
Taking the expectation on both sides given $\F_k$ yields the desired relation.
\end{proof}

\subsection{Proof of Proposition~\ref{prop-xdiff}}\label{Appendix-Propositionprop-xdiff}
\begin{proof}
Adding and subtracting $\hat \bx_k=(\hat x_k,\ldots,\hat x_k)$, we have
\begin{align*}
&\left\|\bx_{k+1}-\bx_k\right\|=\left\|\bx_{k+1} - {\hat \bx}_k +{\hat \bx}_k-\bx_k\right\|\\
\le& \left\|A_k\bx_k-{\hat \bx}_k\right\|+\left\|\bx_k - {\hat \bx}_k\right\| + \bar{\a} \|\by_{k}\|+\bar{\b}\left\|\bx_{k}-\bx_{k-1}\right\|,
\end{align*}
where the last inequality follows from the compact representation of $x$-iterate process (see~\eqref{eq-x-comp}) and the triangle inequality. By the relation for norms in~\eqref{eq-NormIneq}, it follows that 
\begin{align*}
\left\|\bx_{k+1}-\bx_k\right\| \le& \varphi_{k+1} \left\|\bz_k-{\hat \bx}_k\right\|_{\phi_{k+1}}+\varphi_{k}\left\|\bx_k - {\hat \bx}_k\right\|_{\phi_k}\\
+& \bar{\a} \|\by_{k}\|+\bar{\b}\left\|\bx_{k}-\bx_{k-1}\right\|,
 \end{align*}
 where $\varphi_{k}$ is defined in \eqref{eq-const-all}. We notice that by the relation in \eqref{eq-zxhat}, we have
$\displaystyle \left\|\bz_k-{\hat \bx}_k\right\|_{\phi_{k+1}} \le c \left\|\bx_k - {\hat \bx}_k\right\|_{\phi_k}$. Thus, we obtain the desired relation upon using the inequality $(a+b+c)^2\le 3a^2+3b^2+3c^2$.
\end{proof}

\subsection{Proof of Proposition~\ref{prop-ycontract}}\label{Appendix-Propositionprop-ycontract}
\begin{proof}
To simplify the notations, we define
\begin{align*}
\bg_k=\bg(\bx_k,\bxi_{k}) &\text{ and } g_k^i=g_i(x_k^i,\xi_k^i),\\
\nabla_{k}=\nabla F(\bx_k) &\text{ and } \nabla_k^i=\nabla f_i(x_k^i).
\end{align*}
The alternative compact form of the $y$-update in \eqref{eq-y-comp} is
\[\by_{k+1} =  B_k\by_{k} + \bg_{k+1} -\bg_k.\]
By multiplying with $\diag^{-1}(\pi_{k+1})$ and subtracting the matrix $\bs_{k+1}=[\one_n^\T\by_{k+1},\ldots,\one_n^\T\by_{k+1}]^\T\in\re^{n\times p}$, from both sides of the preceding relation, we have for all $k\ge0,$
\begin{align*}
\diag^{-1}(\pi_{k+1})\by_{k+1}-\bs_{k+1} =\left(\diag^{-1}(\pi_{k+1})B_k\by_{k}-\bs_{k}\right)&\\
+\left(\diag^{-1}(\pi_{k+1})-\one_n\one_n^\T\right)\left(\bg_{k+1}-\bg_k\right) 
\triangleq r_1+r_2&,
\end{align*}
since $\bs_{k}\!=\!\one_n\one_n^\T\by_{k}\!=\!\one_n\one_n^\T\bg_k$ (see Lemma~\ref{lem-sumgrad}(a)), and $\diag(u)$ is the diagonal matrix with the vector $u$ on its diagonal.

Taking $\pi_{k+1}$-induced norm on both sides of the preceding equality and using the definition of  $S(\by_{k+1},\pi_{k+1})$ with $\one_n^\T\by_{k+1}=\sum_{j=1}^n y_{k+1}^j$, we have that
\begin{align}\label{eq-srel1}
\!\!S^2(\by_{k+1},\pi_{k+1})
\!=\! \|r_1\|_{\pi_{k+1}}^2 \!+\|r_2\|_{\pi_{k+1}}^2 \!+2\langle r_1, r_2\rangle_{\pi_{k+1}}\!.\!\!\!
\end{align}

\noindent
\textbf{Bounding $\E\!\left[\|r_1\|_{\pi_{k+1}}^2|\F_k\!\right]$}\!: From the definition of $\bs_{k}$, we have 
\begin{align*}
\|r_1\|_{\pi_{k+1}}^2 &= \|\diag^{-1}(\pi_{k+1})B_k\by_{k}-\bs_{k}\|_{\pi_{k+1}}^2\\
&=\sum_{i=1}^n[\pi_{k+1}]_i \Bigg\|\frac{\sum_{j=1}^n [B_k]_{ij} y_{k}^j}{[\pi_{k+1}]_i}-\sum_{j=1}^n y_{k}^j\Bigg\|^2\\
&\le\tau_k^2 \sum_{i=1}^n [\pi_k]_i \Bigg\|\frac{y_k^i}{[\pi_k]_i} - \sum_{j=1}^n y_{k}^j\Bigg\|^2 = \tau_k^2 S^2(\by_{k},\pi_k),
\end{align*}
where we apply Lemma~\ref{lem-basic-ycontract} with the following identifications $\bbG=\bbG_k$, $B=B_k$, $\pi=\pi_{k+1}$, and $\nu=\pi_k$, to obtain the inequality. Here, $\tau_k$ as defined in \eqref{eq-const-all} and the equality follows from the definition of $S(\by_{k},\pi_k)$ in~\eqref{eq-x-S-quants}. Taking the conditional expectation given $\F_k$ and using $\tau_k\le\tau$, yields
\begin{align} \label{eq-r1relbound}
\E[\|\|r_1\|_{\pi_{k+1}}^2|\F_k]\le \tau^2\E[S^2(\by_{k},\pi_k)|\F_k].
\end{align}

\noindent
\textbf{Bounding $\E\left[\|r_2\|_{\pi_{k+1}}^2|\F_k\right]$}: We have
\begin{align*}
\|r_2\|_{\pi_{k+1}}^2 &= \|\left(\diag^{-1}(\pi_{k+1})-\one_n\one_n^\T\right)(\bg_{k+1}-\bg_k)\|_{\pi_{k+1}}^2\\
&\le \|\diag^{-1}(\pi_{k+1})-\one_n\one_n^\T\|_{\pi_{k+1}}^2\|\bg_{k+1}-\bg_k\|^2\\
&\le \psi_k\|\bg_{k+1}-\bg_k\|^2,
\end{align*}
where $\psi_k=n\chi_{k+1}^2-n$. Furthermore,
\begin{align}\label{eq-r2quant}
&\E\left[\|\bg_{k+1}-\bg_k\|^2|\F_k\right]\nonumber\\
=&\E\left[\left\|\bg_{k+1}-\bg_k-\left(\nabla_{k+1}-\nabla_{k}\right)+\nabla_{k+1}-\nabla_{k}\right\|^2|\F_k\right]\nonumber\\
=&\E\!\left[\!\left\|\bg_{k+1}\!\!-\!\bg_k\!\!-\!\left(\nabla_{k+1}\!-\!\nabla_{k}\right)\!\right\|^2\!|\F_k\!\right]\!\!+\!\E\!\left[\!\left\|\nabla_{k+1}\!-\!\nabla_{k}\right\|^2\!|\F_k\!\right]\nonumber\\
+&2\E\left[\langle\bg_{k+1}-\bg_k-\left(\nabla_{k+1}-\nabla_{k}\right),\nabla_{k+1}-\nabla_{k}\rangle|\F_k\right].\!\!
\end{align}
For the first term in \eqref{eq-r2quant}, since $\F_k \subset \F_{k+1}$, we have
\begin{align*}
&\E\left[\langle\bg_{k+1}-\nabla_{k+1},\bg_k-\nabla_{k}\rangle|\F_k\right]\\
=&\E\left[\E\left[\langle\bg_{k+1}-\nabla_{k+1},\bg_k-\nabla_{k}\rangle|\F_{k+1}\right]\big|\F_k\right]\\
=&\E\left[\langle\E\left[\bg_{k+1}-\nabla_{k+1}|\F_{k+1}\right],\bg_k-\nabla_{k}\rangle\big|\F_k\right]=0,
\end{align*} 
where the last equality follows from Assumption~\ref{asm-SFO}(a). Thus,
\begin{align}\label{eq-r2firstterm}
&\E\left[\left\|\bg_{k+1}-\bg_k-\left(\nabla_{k+1}-\nabla_{k}\right)\right\|^2|\F_k\right]\\
=&\E\left[\left\|\bg_{k+1}-\nabla_{k+1}\right\|^2|\F_k\right]+\E\left[\left\|\bg_k-\nabla_{k}\right\|^2|\F_k\right]\le 2n\sigma^2\!\!,\nonumber
\end{align}
according to Assumption~\ref{asm-SFO}(b).

For the second term in \eqref{eq-r2quant}, we have
\begin{align*}
&\left\|\nabla_{k+1}-\nabla_{k}\right\|^2 = \sum_{i=1}^n \|\nabla f_i(x_{k+1}^i)- \nabla f_i(x_k^i)\|^2\\
\le&L^2\sum_{i=1}^n \|x_{k+1}^i - x_k^i\|^2=L^2\|\bx_{k+1}-\bx_k\|^2,
% \le & 2L^2\left(c\varphi_{k+1}+\varphi_k \right)^2 \!\|\bx_k-\hat{\bx}_k\|_{\phi_k}^2 \!+ 2\alpha^2L^2 \|\by_{k}\|_{\pi_k^{-1}}^2,
\end{align*}
by using the Lipschitz continuity of the gradients $\nabla f_i$. Hence, using Proposition~\ref{prop-xdiff}, we have
\begin{align}\label{eq-r2secondterm}
&\E\!\left[\left\|\nabla_{k+1}-\nabla_{k}\right\|^2|\F_k\right] 
\le 3L^2\!\left(c\varphi_{k+1}+\varphi_k \right)^2 \!\|\bx_k\!-\hat{\bx}_k\|_{\phi_k}^2 \! \nonumber\\
+&  3\bar{\a}^2L^2\E\left[\|\by_{k}\|_{\pi_k^{-1}}^2|\F_k\right] +3L^2\bar{\b}^2\left\|\bx_{k}-\bx_{k-1}\right\|^2.
\end{align}

For the last term in \eqref{eq-r2quant}, we first note
\begin{align*}
&\E\left[\langle\bg_{k+1}-\bg_k-\left(\nabla_{k+1}-\nabla_{k}\right),\nabla_{k}\rangle|\F_k\right]\\
=&\langle\E\left[\bg_{k+1}-\nabla_{k+1}|\F_k\right]-\E\left[\bg_k-\nabla_{k}|\F_k\right],\nabla_{k}\rangle=0,
\end{align*}
and
\begin{align*}
&\E\left[\langle\bg_{k+1}-\bg_k-\left(\nabla_{k+1}-\nabla_{k}\right),\nabla_{k+1}\rangle|\F_k\right]\\
=&\E\left[\E\left[\langle\bg_{k+1}-\bg_k-\left(\nabla_{k+1}-\nabla_{k}\right),\nabla_{k+1}\rangle|\F_{k+1}\right]\big|\F_k\right]\\
=&\E\left[\langle\E\left[\bg_{k+1}-\nabla_{k+1}|\F_{k+1}\right]-\left(\bg_k-\nabla_{k}\right),\nabla_{k+1}\rangle\big|\F_k\right]\\
=&-\E\left[\langle\bg_k-\nabla_{k},\nabla_{k+1}\rangle\big|\F_k\right].
\end{align*}
Combining the preceding two relations, the last term in \eqref{eq-r2quant} can be written as follows:
\begin{align*}
&2\E\left[\langle\bg_{k+1}-\bg_k-\left(\nabla_{k+1}-\nabla_{k}\right),\nabla_{k+1}-\nabla_{k}\rangle|\F_k\right]\\
=&2\E\left[\langle\nabla_{k+1},\nabla_{k}-\bg_k\rangle\big|\F_k\right].
\end{align*}
According to the update in \eqref{eq-x} and \eqref{eq-y}, $\nabla_{k+1}^i$ is equal to
\[\nabla\! f_i\!\!\left(\sum_{j=1}^n[A_k]_{ij}x_{k}^j  \!\!-\! \a_i \!\Bigg(\!\sum_{j=1}^n [B_{k-1}]_{ij}y_{k-1}^j \!+\! g_k^i \!-\!g_{k-1}^i\!\!\Bigg)\!\!+\!\b_i u_{k}^i\!\!\right)\!\!.\]
We also define $\tilde{\nabla}_{k+1}^i$ to be equal to the following expression:
\[\nabla\! f_i\!\left(\sum_{j=1}^n[A_k]_{ij}x_{k}^j  \!\!-\! \a_i\! \Bigg(\!\sum_{j=1}^n [B_{k-1}]_{ij}y_{k-1}^j \!\!+\! \nabla_k^i \!\!-\!g_{k-1}^i\Bigg)\!\!+\!\b_i u_{k}^i\!\!\right)\!\!.\]
By using the Lipschitz continuity of the gradients, we have
\[\|\nabla_{k+1}^i-\tilde{\nabla}_{k+1}^i\|\le \bar{\a} L\|\nabla_k^i-g_k^i\|.\]

We then proceed to bound the last term in \eqref{eq-r2quant} as follows:
\begin{align}\label{eq-r2thirdterm}
&\E\left[\langle\nabla_{k+1},\nabla_{k}-\bg_k\rangle\big|\F_k\right]\nonumber\\
= & \sum_{i=1}^n\E\left[\langle\nabla_{k+1}^i-\tilde{\nabla}_{k+1}^i,\nabla_k^i-g_k^i\rangle\big|\F_k\right]\nonumber\\
\le & \sum_{i=1}^n\E\left[\|\nabla_{k+1}^i-\tilde{\nabla}_{k+1}^i\|\|\nabla_k^i-g_k^i\|\big|\F_k\right]\nonumber\\
\le &~ \bar{\a} L\sum_{i=1}^n\E\left[\|\nabla_k^i-g_k^i\|^2\big|\F_k\right] ~~\le~~ \bar{\a} nL\sigma^2,
\end{align}
where we use Assumption~\ref{asm-SFO} and $\E\left[\langle\tilde{\nabla}_{k+1}^i,\nabla_k^i-g_k^i\rangle|\F_k\right] =\langle\tilde{\nabla}_{k+1}^i,\E\left[\nabla_k^i-g_k^i|\F_k\right]\rangle = 0$.

Substituting \eqref{eq-r2firstterm}, \eqref{eq-r2secondterm} and \eqref{eq-r2thirdterm} to \eqref{eq-r2quant}, we obtain 
\begin{align}\label{eq-r2bound}
\E&\left[\|r_2\|_{\pi_{k+1}}^2|\F_k\right] \le 3L^2\psi_k\left(c\varphi_{k+1}+\varphi_k \right)^2 \!\|\bx_k-\hat{\bx}_k\|_{\phi_k}^2 \nonumber\\
&+3\bar{\a}^2\psi_kL^2\E\left[\|\by_{k}\|_{\pi_k^{-1}}^2|\F_k\right]+2(1+\bar{\a} L)n\psi_k\sigma^2\nonumber\\
&+3\bar{\b}^2\psi_kL^2\left\|\bx_{k}-\bx_{k-1}\right\|^2.
\end{align}

\noindent
\textbf{Bounding $\E\left[2\langle r_1, r_2\rangle_{\pi_{k+1}}|\F_k\right]$}: We first note that
\begin{align*}
&\langle \diag^{-1}(\pi_{k+1})B_k\by_{k}-\bs_{k},\one_n\one_n^\T\left(\bg_{k+1}-\bg_k\right)\rangle_{\pi_{k+1}}\\
=&\langle \left(B_k-\diag(\pi_{k+1})\one_n\one_n^\T\right)\by_{k},\one_n\one_n^\T\left(\bg_{k+1}-\bg_k\right)\rangle=0,
\end{align*}
where the first equality follows from $\bs_{k}=\one_n\one_n^\T\by_{k}$. The last equality is obtained since $\pi_{k+1}$ is a stochastic vector, $B_k$ is a column stochastic matrix and the matrix $\one_n\one_n^\T\left(\bg_{k+1}-\bg_k\right)$ is consensus. Therefore,
\begin{align*}
&\langle r_1, r_2\rangle_{\pi_{k+1}} \\
=& \langle \diag^{-1}\!(\pi_{k+1})B_k\by_{k}\!-\bs_{k},\diag^{-1}\!(\pi_{k+1})\left(\bg_{k+1}\!-\bg_k\right)\rangle_{\pi_{k+1}}\\
=& \langle \diag^{-1}(\pi_{k+1})B_k\by_{k}-\bs_{k},\bg_{k+1}-\bg_k\rangle\\
=& \langle \diag^{-1}(\pi_{k+1})B_k\by_{k},\nabla_{k}-\bg_k\rangle
+\langle -\bs_{k},\nabla_{k}-\bg_k\rangle\\
+&\langle \diag^{-1}(\pi_{k+1})B_k\by_{k}-\bs_{k},\bg_{k+1}-\nabla_{k}\rangle\triangleq s_1+s_2+s_3.
\end{align*}
Hence,
\begin{align}\label{eq-r1r2}
&\E\left[2\langle r_1, r_2\rangle_{\pi_{k+1}}|\F_k\right] \nonumber\\
=~& 2\E\left[s_1|\F_k\right]+2\E\left[s_2|\F_k\right]+2\E\left[s_3|\F_k\right].
\end{align}
Using the $y$-update in compact form in \eqref{eq-y-comp}, we have
\begin{align*}
s_1=& \langle B_kB_{k-1}\by_{k-1} + B_k\left(\bg_k -\bg_{k-1}\right),\nabla_{k}-\bg_k\rangle_{\pi_{k+1}^{-1}}.
\end{align*}
Noting that 
\begin{align*}
&\E\left[\langle B_kB_{k-1}\by_{k-1} - B_k\bg_{k-1},\nabla_{k}-\bg_k\rangle_{\pi_{k+1}^{-1}}|\F_k\right]\\
=&\langle B_kB_{k-1}\by_{k-1} - B_k\bg_{k-1},\E\left[\nabla_{k}-\bg_k|\F_k\right]\rangle_{\pi_{k+1}^{-1}}=0,
\end{align*}
we further obtain
\begin{align}\label{eq-s1}
&\E\left[s_1|\F_k\right]=\E\left[\langle B_k\bg_k,\nabla_{k}-\bg_k\rangle_{\pi_{k+1}^{-1}}|\F_k\right]\nonumber\\
=&\sum_{i=1}^n\frac{1}{[\pi_{k+1}]_i} \E\left[\left\langle\sum_{j=1}^n [B_k]_{ij} g_{k}^j,\nabla_k^i-g_k^i\right\rangle\bigg|\F_k\right]\nonumber\\
=&\sum_{i=1}^n\frac{[B_k]_{ii}}{[\pi_{k+1}]_i} \E\left[\left\langle  g_k^i,\nabla_k^i-g_k^i\right\rangle|\F_k\right]\nonumber\\
=&\sum_{i=1}^n\!\frac{[B_k]_{ii}}{[\pi_{k+1}]_i} \!\!\left(\!-\E\!\left[\left\|g_k^i\!-\!\nabla_k^i\right\|^2 \!|\F_k\right]\!+\!\E\!\left[\left\langle  \nabla_k^i,\!\nabla_k^i\!-\!g_k^i\right\rangle\!|\F_k\!\right]\right)\nonumber\\
=&-\sum_{i=1}^n\!\frac{[B_k]_{ii}}{[\pi_{k+1}]_i} \E\!\left[\left\|g_k^i-\nabla_k^i\right\|^2 |\F_k\right]\le 0,
\end{align}
where we use the assumption that $g_k^i$’s are independent across the nodes and the fact that $\E\!\left[\left\langle  \nabla_k^i,\nabla_k^i\!-\!g_k^i\right\rangle|\F_k\!\right]=0$.

We proceed similarly for $s_2$. Since $\bs_{k}\!=\!\one_n\one_n^\T\bg_k$, we have
\begin{align}\label{eq-s2}
\E\left[s_2|\F_k\right]=&\E\left[\langle -\one_n\one_n^\T\bg_k,\nabla_{k}-\bg_k\rangle|\F_k\right]\nonumber\\
=&\sum_{i=1}^n\E\left[-\left\langle\sum_{j=1}^n g_{k}^j,\nabla_k^i-g_k^i\right\rangle\bigg|\F_k\right]\nonumber\\
=&\sum_{i=1}^n\E\left[\left\langle -g_k^i,\nabla_k^i-g_k^i\right\rangle|\F_k\right]\nonumber\\
=&\sum_{i=1}^n\E\left[\left\langle \nabla_k^i-g_k^i,\nabla_k^i-g_k^i\right\rangle|\F_k\right]\le n\sigma^2.
\end{align}

Finally, for $s_3$, we first note that
\begin{align*}
\E\left[s_3|\F_k\right]=&\E\left[\langle \diag^{-1}(\pi_{k+1})B_k\by_{k}-\bs_{k},\bg_{k+1}-\nabla_{k}\rangle|\F_k\right]\\
=&\E\left[\langle \diag^{-1}(\pi_{k+1})B_k\by_{k}-\bs_{k},\nabla_{k+1}-\nabla_{k}\rangle|\F_k\right]\!.
\end{align*}
By the Cauchy–Schwarz inequality, we have
\begin{align*}
&\langle \diag^{-1}(\pi_{k+1})B_k\by_{k}-\bs_{k},\nabla_{k+1}-\nabla_{k}\rangle\\
\le& \|\diag^{-1}(\pi_{k+1})B_k\by_{k}-\bs_{k}\| \|\nabla_{k+1}-\nabla_{k}\|\\
\le& \tfrac{1-\tau^2}{4\tau^2}\|r_1\|_{\pi_{k+1}}^2+\tfrac{\chi_{k+1}^2\tau^2}{1-\tau^2}\|\nabla_{k+1}-\nabla_{k}\|^2,
\end{align*}
recalling that $\|r_1\|_{\pi_{k+1}}^2 = \|\diag^{-1}(\pi_{k+1})B_k\by_{k}-\bs_{k}\|_{\pi_{k+1}}^2$.

Using the relations in \eqref{eq-r1relbound} and \eqref{eq-r2secondterm} and the preceding relation, we obtain the bound for $\E\left[s_3|\F_k\right]$ as follows:
\begin{align} \label{eq-s3}
&\E\left[s_3|\F_k\right] \le \tfrac{3L^2\left(c\varphi_{k+1}+\varphi_k \right)^2\chi_{k+1}^2\tau^2}{1-\tau^2}\|\bx_k-\hat{\bx}_k\|_{\phi_k}^2 \nonumber\\
&+\tfrac{1-\tau^2}{4}\E[S^2(\by_{k},\pi_k)|\F_k]+\tfrac{3\bar{\a}^2L^2\chi_{k+1}^2\tau^2}{1-\tau^2}\E\left[\|\by_{k}\|_{\pi_k^{-1}}^2|\F_k\right]\nonumber\\
&+\tfrac{3\bar{\b}^2L^2\chi_{k+1}^2\tau^2}{1-\tau^2}\left\|\bx_{k}-\bx_{k-1}\right\|^2.
\end{align}
Substituting the bounds in \eqref{eq-s1}, \eqref{eq-s2} and \eqref{eq-s3} into \eqref{eq-r1r2} yields
\begin{align}\label{eq-r1r2bound}
\E&\left[2\langle r_1, r_2\rangle_{\pi_{k+1}}|\F_k\right] \le \tfrac{6L^2(c\varphi_{k+1}+\varphi_k )^2\!\chi_{k+1}^2\tau^2}{1-\tau^2\!}\|\bx_k-\hat{\bx}_k\|_{\phi_k}^2 \nonumber\\
&+\tfrac{1-\tau^2}{2}\E[S^2(\by_{k},\pi_k)|\F_k]+\tfrac{6\bar{\a}^2L^2\chi_{k+1}^2\tau^2}{1-\tau^2}\E\left[\|\by_{k}\|_{\pi_k^{-1}}^2|\F_k\right]\nonumber\\
&+\tfrac{6\bar{\b}^2L^2\chi_{k+1}^2\tau^2}{1-\tau^2}\left\|\bx_{k}-\bx_{k-1}\right\|^2+2n\sigma^2.
\end{align}

Combining the bounds in \eqref{eq-r1relbound}, \eqref{eq-r2bound}, \eqref{eq-r1r2bound} with the relation in \eqref{eq-srel1} completes the proof.
\end{proof}

\subsection{Proof of Lemma~\ref{lem-ypik}}\label{Appendix-Lemmalem-ypik}
\begin{proof}
% Under given assumptions, by 
From Lemma~\ref{lem-bmatrices} we have that the stochastic vectors $\pi_k$, $k\ge0$, defined in~\eqref{eq-pik}, all have positive entries so 
% that the scaling with $[\pi_k]_i^{-1}$ is well defined for all $i\in[n]$ and $k\ge0$. In particular 
the $\pi_k^{-1}-$induced norm is well defined.
Noting that
\[\|\by_{k}\|^2_{\pi_k^{-1}}=\sum_{i=1}^n \frac{\|y_k^i\|^2}{[\pi_k]_i} =\sum_{i=1}^n [\pi_k]_i\,\left\|\frac{y_k^i}{[\pi_k]_i} \right\|^2,\]
and using Lemma~\ref{lem-normlincomb} for the weighted average of vectors, where $\g_i=[\pi_k]_i$ and $u_i=y_k^i/[\pi_k]_i$ for all $i$, and $u=0$, we obtain
\begin{align*}
\sum_{i=1}^n [\pi_k]_i\,\Bigg\|\frac{y_k^i}{[\pi_k]_i} \Bigg\|^2&=\sum_{i=1}^n [\pi_k]_i\Bigg\|\frac{y_k^i}{[\pi_k]_i} -\sum_{j=1}^n y_{k}^j\Bigg\|^2 +\Bigg\|\sum_{j=1}^n y_{k}^j\Bigg\|^2\\
&=S^2(\by_{k},\pi_k) +\Bigg\|\sum_{j=1}^n y_{k}^j\Bigg\|^2,
 \end{align*}
which establish the proof of the lemma upon taking the conditional expectation given $\F_k$ and using Lemma~\ref{lem-ysum}.
\end{proof}

\end{document}